\documentclass{article}
\usepackage{titlecaps}
\Addlcwords{an and the}
\usepackage[utf8]{inputenc}
\usepackage{amsmath}
\usepackage{amsfonts}
\usepackage{amssymb}
\usepackage{algorithm}
\usepackage[noend]{algpseudocode}
\usepackage{multirow}
\usepackage{float}
\usepackage{array}
\usepackage{lipsum}
\usepackage[left=0.75in,right=0.75in,top=0.75in,bottom=0.75in]{geometry}
\usepackage{comment}
\usepackage{graphicx}
\usepackage{url}
\usepackage{natbib}
\usepackage{color}
\usepackage{pgfplots}
\usepackage{mathtools}
\usepackage{caption}
\usepackage{subcaption}
\usepackage{pdflscape}
\usepackage{booktabs} 
\usepackage{enumerate}
\usepackage{authblk}
\usepackage{blkarray}
\usepackage{bm}
\usepackage{marvosym} 
\usepackage{soul}
\usepackage{textcomp, xspace}
\usepackage{lineno}

\newtheorem{definition}{Definition}

\usepgfplotslibrary{dateplot}
\usepackage[colorlinks]{hyperref}
\hypersetup{bookmarks=true, citecolor=blue, urlcolor=blue, linkcolor=blue}
\newcommand{\mys}{\hspace{1mm}}

\definecolor{greenxllite}{RGB}{196,215,155}
\definecolor{bluexllite}{RGB}{149,179,215}
\definecolor{redxllite}{RGB}{218,150,148}
\definecolor{greenxl}{RGB}{155, 187, 89}
\definecolor{bluexl}{RGB}{79,129,189}
\definecolor{redxl}{RGB}{192, 80, 77}
\definecolor{OliveGreen}{RGB}{70,136,52}

\usetikzlibrary{patterns}
\algrenewcommand\algorithmicrequire{\textbf{Input:}}
\algrenewcommand\algorithmicensure{\textbf{Output:}}
\newcommand{\say}[1]{``#1''}

\newcommand{\primalgraph}{\mathbb{G}}
\newcommand{\primalvertex}{\mathbb{V}}
\newcommand{\primaledge}{\mathbb{E}}
\newcommand{\primalterminal}{\mathbb{E}_T}

\newcommand{\dualgraph}{L}
\newcommand{\dualvertex}{V}
\newcommand{\dualedge}{E}
\newcommand{\dualterminal}{V_T}

\newcommand{\paertooptimalset}{Z^{*}}
\newcommand{\localsearchndset}{Z}
\newcommand{\setY}{Y}
\newcommand{\efffrontier}[1]{C(#1)}

\newcommand{\explorationset}{X}
\newcommand{\tourset}{P}
\newcommand{\violatedterminalset}[1]{V^{inf}_{T}(#1)}

\newcommand{\subpathcollection}{\hat{P}}
\newcommand{\feasibleinsertionset}[3]{V_{T}^{#1}(#2,#3)}

\newcommand{\pathbanksubset}[2]{P(#1, #2)}
\newcommand{\pathbank}[2]{P^{*}(#1, #2)}

\newcommand{\tourindex}[2]{#1_{#2}}

\newcommand{\objectivevector}[1]{\mathcal{C}({#1})}  
\newcommand{\adjterminalvector}[1]{\mathcal{S}_{#1}}

\newcommand{\temrinalPriority}[1]{\pi(#1)}

\newcommand{\subpath}[3]{#1_{#2, #3}}
\newcommand{\edgevar}[4]{x_{#1,#2,#3,#4}}
\newcommand{\arrivaltimevar}[2]{s_{#1,#2}}
\newcommand{\twnodevar}[2]{y_{#1,#2}}
\newcommand{\commodityflow}[4]{f_{#1,#2,#3,#4}}

\newcommand{\sizedualterminal}{{n_T}}

\newcommand{\energy}[1]{h_{#1}}
\newcommand{\turn}[1]{l_{#1}}
\newcommand{\meanenergy}[1]{\Bar{h}_{#1}}
\newcommand{\stdenergy}{\sigma_h}
\newcommand{\stdturns}{\sigma_l}
\newcommand{\meanturns}[1]{\Bar{l}_{#1}}

\newcommand{\hyperParameter}[1]{\xi_{#1}}

\newcommand{\earliesttime}[1]{a_{#1}}
\newcommand{\latesttime}[1]{b_{#1}}
\newcommand{\traveltime}[1]{t_{#1}}

\newcommand{\length}[1]{|{#1}|}
\newcommand{\penaltyfunction}[1]{\Phi(#1)}
\newcommand{\patharrivaltime}[2]{\Lambda(#1,#2)}
\newcommand{\gain}[2]{\Gamma{(#1,#2)}}

\newcommand{\candidates}{D}

\newcommand{\depot}{0}
\newcommand{\concatenate}{\oplus}

\DeclareMathOperator*{\argmax}{arg\,max}

\allowdisplaybreaks

\newcommand{\Keywords}[1]{\par\noindent
{\small{\em \textbf{Keywords}\/}: #1}}

\title{A Bi-criterion Steiner Traveling Salesperson Problem with Time Windows for Last-Mile Electric Vehicle Logistics}

\author[1] {Prateek Agarwal*}
\author[2] {Debojjal Bagchi*}
\author[1, 3] {Tarun Rambha} 
\author[4] {Venktesh Pandey} 

\affil[1]{\small Department of Civil Engineering, Indian Institute of Science (IISc), Bengaluru, India}
\affil[2]{\small Department of Civil, Architectural and Environmental Engineering, The University of Texas at Austin, Austin, USA}
\affil[3]{\small Center for infrastructure, Sustainable Transportation and Urban Planning (C\textit{i}STUP), Indian Institute of Science (IISc),
Bengaluru, India}
\affil[4]{\small Department of Civil, Architectural, and Environmental Engineering, North Carolina Agricultural and Technical State University, Greensboro, USA}

\date{ }
\newcommand\blfootnote[1]{%
  \begingroup
  \renewcommand\thefootnote{}\footnote{#1}%
  \addtocounter{footnote}{-1}%
  \endgroup
}

\begin{document}
\maketitle
\blfootnote{* These authors contributed equally to this manuscript.}
\vspace{-5mm}
\begin{abstract}
This paper addresses the problem of energy-efficient and safe routing of last-mile electric freight vehicles. With the rising environmental footprint of the transportation sector and the growing popularity of E-Commerce, freight companies are likely to benefit from optimal time-window-feasible tours that minimize energy usage while reducing traffic conflicts at intersections and thereby improving safety. We formulate this problem as a Bi-criterion Steiner Traveling Salesperson Problem with Time Windows (BSTSPTW) with energy consumed and the number of left turns at intersections as the two objectives while also considering regenerative braking capabilities. We first discuss an exact mixed-integer programming model with scalarization to enumerate points on the efficiency frontier for small instances. For larger networks, we develop an efficient local search-based heuristic, which uses several operators to intensify and diversify the search process. We demonstrate the utility of the proposed methods using benchmark data and real-world instances from Amazon delivery routes in Austin, US. Comparisons with state-of-the-art solvers shows that our heuristics can generate near-optimal solutions within reasonable time budgets, effectively balancing energy efficiency and safety under practical delivery constraints.

\vspace{4mm}
\end{abstract}

\section{Introduction}
\label{sec:intro}
According to \cite{ieq}, global CO$_2$ emissions from the transport sector were nearly eight gigatons in 2022. There is, thus, a pressing need for energy-efficient routing, especially for freight trips, which continues to grow rapidly with the widespread adoption of E-Commerce. For instance, Amazon delivered approximately 4.75 billion packages in 2021 with an estimated carbon footprint of 71.54 million metric tons of CO$_{2}$e \citep{amazon, statista2023}. Safety is another critical issue central to transportation systems. Approximately 1.19 million commuters die in traffic crashes annually, costing countries an average of 3\% of their gross domestic product \citep{who2023}. Many of these accidents directly involve logistics operations \citep{castillo2016exploring}. For example, research from the Reporting of Injuries, Diseases, and Dangerous Occurrences Regulations reveals that the transportation sector in the UK has a work-related fatal injury rate nearly twice the average across all industries \citep{hse_fatalinjuries}. 

In this context, a popular study by United Parcel Service \citep{holland2017ups} (UPS) on last-mile routing of their gasoline fleet showed that minimizing distance while avoiding left turns (and satisfying other domain-specific side constraints) can significantly reduce fuel use and accidents. They documented that tours discovered by their proprietary ORION software reduced carbon emissions by 20,000 tons, fuel consumption by 10 million gallons, and increased successful package deliveries by 350,000 annually. Extensive studies have also shown that turn movements at intersections significantly impact accident rates \citep{wang2008modeling, wood2020analyzing}. For example, based on the National Motor Vehicle Crash Causation Survey data, \cite{choi2010crash} found that 61\% of intersection-related accidents in the US occur during left-hand turns. In this paper, we address a related problem in the context of Electric Vehicles (EVs), employing bi-criterion algorithms to optimize both the energy consumed and the number of left turns in a tour while also considering regenerative braking capabilities.

EVs present an opportunity to reduce tail-pipe emissions and improve urban air quality. However, models for the energy consumption of EVs are nascent and evolving. They often include various factors such as vehicle type, road gradient, speed and driving style, traffic congestion, battery management system, State of Charge levels, thermal management, and weather \citep{keskin2016partial, roberti2016electric, berzi2016development, zhao2020development,lin2021electric}. Several studies have focused on estimating energy consumption for range prediction using physics-based and machine-learning models \citep{varga2019prediction, topic2019neural}. Many electric vehicles can also gain energy through regenerative braking, which recharges batteries by recovering some kinetic energy from vehicles. Consequently, from a modeling perspective, when we represent a network using nodes and edges and find paths that minimize energy usage, the underlying graph can have negative weights for road segments where energy is gained rather than expended. Eco-routing methods that account for the effect of road features have received little attention in logistics. The impact of turns is even more neglected. Although there is some research on reducing the number of turns through better network design \citep{eichler2013vortex, boyles2014equilibrium}, studies specific to delivery fleets have been sparse. Avoiding turns can significantly improve the quality and safety of a route, especially for large vehicles. However, turning movements in a roadway graph are determined by three nodes: the current node, a predecessor, and a successor. This makes it difficult to directly apply shortest path algorithms that optimize turns without additional graph transformations. In this paper, we address these challenges, with the formal problem statement defined next.

\subsection{Problem Statement}
In this paper, the problem of routing a single vehicle is considered and formulated as a variant of the well-known Traveling Salesperson Problem (TSP). The objective of the TSP is to find the shortest tour that visits each node in a graph exactly once.\footnote{Extensions to the more general vehicle routing problem can be carried out using a cluster-first, route-second approach, in which the customers are first grouped and a TSP is solved for each group \citep{fischer81}.} Last-mile logistics problems can be formulated as a Steiner Traveling Salesperson Problem (STSP) because it is enough to visit a subset of nodes or customer locations, which are also called \textit{terminals}. While optimizing the sequence of customers visited, evaluating the exact paths these vehicles take between terminals is crucial for problems with time windows and multiple objectives. Deliveries can be time-sensitive, and hence, visits must fall within specific time windows. Finally, since we seek tours that optimize energy consumption and left turns, the problem is treated as a Bi-criterion Steiner Traveling Salesperson Problem with Time Windows (BSTSPTW). Vehicles are allowed to wait to meet time limits in the event of early arrivals. 

Our work is motivated by a practical problem that logistics companies face in delivering goods while minimizing energy consumption and enhancing safety through routing. By considering energy consumption and the number of left turns as the two objectives, our BSTSPTW formulation finds Pareto-optimal tours that visit all terminals (at least once since revisits may be optimal in STSP) and return to the depot, providing a range of routing options for drivers and managers to choose from based on their preferences and trade-offs between the two objectives. Throughout the paper, we assume right-hand moving traffic and the term \textit{turns} refers specifically to the left turns. However, one can easily modify the methods presented in this paper to left-handed traffic scenarios.

Moreover, in the event of rerouting (due to delays from non-recurring forms of congestion), the ability to switch to alternative solutions on the Pareto frontier facilitates adaptive decision-making. Therefore, the proposed models provide valuable managerial insights and decision-making support to logistics firms, thereby reducing environmental impacts and improving safety. Although the paper is set against the backdrop of finding tours that minimize energy and left turns, the methods developed for the BSTSPTW are generic. They can also be used for other applications, such as problems involving cost/distance and time objectives.

\subsection{Contributions}
The major contributions of this paper are listed below. While the first two highlight methodological advances, the latter two summarize the practical value of this research.
\begin{itemize}
    \item \textbf{Exact integer-programming-based method:} A new Mixed Integer Program (MIP) formulation for BSTSPTW is presented. Specifically, our model improves the single commodity flow formulation for STSP \citep{letchford2013compact} to allow potential node and edge revisits and handle time-window constraints. Additionally, it incorporates multiple objectives and generates the efficiency frontier using scalarization. The MIP runtimes scale with the number of edges in the network. Consequently, this approach is useful for small networks.   

    \item \textbf{Local search heuristics:} Recognizing the computational challenges associated with large-scale, real-world instances, we propose a novel local search-based heuristic. Specifically, six operators are designed in the heuristic : \textsc{S3Opt}, \textsc{S3OptTW}, \textsc{RepairTW}, \textsc{FixedPerm}, \textsc{Quad} and \textsc{RandPermute}. Operators \textsc{S3Opt} and \textsc{S3OptTW} extend the conventional 3-opt move to the STSP.  The \textsc{RepairTW} operator makes a tour time-window feasible by destroying and repairing parts of the tour. \textsc{FixedPerm} improves the tours by optimizing the paths between customer locations. Lastly, \textsc{Quad} and \textsc{RandPermute} help escape local minima by introducing necessary diversification.

    \item \textbf{Real-world applications:} To showcase the effectiveness of the proposed approaches, an elaborate case study is presented using data from real-world networks, focusing on Amazon delivery routes in Austin, US \citep{merchan20222021}. Our findings reveal that the proposed approach runs within a reasonable computational budget of two hours and offers drivers multiple routes to choose from. 
    \item \textbf{Benchmarking:} For smaller Solomon-Potvin-Bengio datasets \citep{lopez2013travelling}, we match the results with MIP solutions. Given that the MIP fails to scale for real-world Amazon instances, we showcase the quality of the tours from our methods with those produced by adapting Lin–Kernighan-Helsgaun's heuristic and STSP solution techniques \citep{alvarez2019note}.
    \end{itemize}

The remainder of the paper is organized as follows. Section \ref{sec:litreview} summarizes the literature on related studies. The problem definition and the notation used in this work, and an example, are outlined in Section \ref{sec:preliminaries}. Section \ref{sec:exact} formulates a MIP for the BSTSPTW. We introduce local search heuristics in Section \ref{sec:heuristics} to address the computational challenges of real-world scenarios involving large graphs. Section \ref{sec:results} analyzes the performance of the proposed MIP and local search. Finally, in Section \ref{sec:conclusion}, we recap our findings and discuss future research directions.

\section{Literature Review}\label{sec:litreview}
This section reviews solution approaches for the TSP and its time-window variant, followed by literature on the multi-objective and Steiner versions. Finally, we examine specific studies that apply such methods in the context of EVs. 
The goal of the TSP is to find the shortest tour that visits all nodes in a graph and returns to the starting point. Typical objectives involve minimizing time, distance, or cost, and the graphs are assumed symmetric, with the cost of traveling from node $i$ to node $j$ being the same as that of $j$ to $i$. Over the years, several exact approaches such as cutting plane methods, branch-and-bound, and branch-and-cut have been proposed to solve Integer Programs (IPs) that model single-objective symmetric TSP. However, these methods are computationally expensive and do not scale well for large instances. Hence, it is common practice to employ heuristics such as neighborhood search, tour-improvement procedures such as the $k$-opt algorithm, tabu search, and insertion algorithms. 
Comprehensive reviews of exact and heuristic approaches to the TSP can be found in \cite{applegate2011traveling}. Other closely related problems to the TSP are the Chinese Postman Problem (CPP) and the Rural Postman Problem (RPP). CPP aims to find the optimal tour that visits every edge at least once \citep{aminu2006constraint}. The RPP, on the other hand, requires only a subset of edges to be serviced \citep{corberan2024theoretical}. Studies such as \cite{monroy2017adaptive} have also addressed variants of the RPP with time windows. The solutions to these problems may require edge revisits similar to our work and are handled using graph transformations. However, time window constraints are imposed only for the first visit to a required edge. 

Relaxing the cost symmetry on the edges yields the Asymmetric Traveling Salesperson Problem (ATSP). Heuristics designed for the symmetric TSP often fail or are ineffective for ATSP \citep{choi2003genetic}. Among ATSP heuristics, while 2-opt moves do not yield good results for the ATSP due to the inversion of tour segments, 3-opt-based methods can be effective. \cite{kanellakis1980local} proposed another popular ATSP heuristic that builds on the Lin-Kernighan (LK) heuristic for the symmetric TSP. It employs odd $k$-opt swaps and uses a depth-first search to find the best $k$ and double-bridge type exchanges for diversification. Helsgaun's heuristic \citep{helsgaun2000effective} is also a popular method for the ATSP, which transforms an ATSP instance into a symmetric TSP instance.

A practical variant of the TSP is the Traveling Salesperson Problem with Time Windows (TSPTW), where the objective is to minimize cost while visiting each node within its time window. \cite{savelsbergh1985local} showed that finding feasible solutions to the TSPTW is NP-complete and proposed a $k$-opt-based solution. Another common technique applies penalties to nodes when time windows are infeasible \citep{da2010general}. For example, \cite{ohlmann2007compressed} constructed a simulated annealing procedure with a variable penalty approach, and \cite{helsgaun2017extension} extended the LK heuristic for the symmetric TSP using separate penalties instead of adding them to the objective. An alternative objective in the TSPTW is to minimize the \textit{makespan} of the tour, i.e., total travel time \citep{kara2013new}. Modifications of the 2-opt move \citep{potvin1995exchange} and pre-processing steps \citep{dumas1995optimal} for the TSPTW have also been proposed. Recently, \cite{pralet2023iterated} employed a large neighborhood search based on the concept of \textit{insertion width}, which measures the amount of perturbation in the destruction-repair-based neighborhoods.

Objectives in real-world applications are usually multifold. However, the TSP is an NP-hard problem, and incorporating multiple objectives introduces additional challenges. As a result, heuristic methods are commonly employed to solve the Multi-Objective Traveling Salesperson Problem (MOTSP) based on local search or evolutionary algorithms. \cite{paquete2003two} presented a two-phase local search where, in the first phase, single-criteria versions of the problem were explored, and the second phase investigated the multi-objective version using weights on the objectives. \cite{angel2004dynasearch} introduced the ds-2-opt (dyna-search) neighborhood, specifically designed for exploring the non-dominated neighborhood of a tour. They used a dynamic program to calculate all independent 2-opt moves within a neighborhood. \cite{paquete2004pareto} proposed a Pareto local search method akin to \cite{angel2004dynasearch} with the distinction that their neighborhood was smaller and did not include any tour generated from a dominated solution. While this method yields satisfactory results, it is computationally expensive. To overcome this drawback, \cite{lust2010two} proposed a two-phase Pareto local search in which the first phase computes the convex subset of efficient solutions using weighted objectives, and the second phase generates the remaining efficient solutions using exchanges. Details on evolutionary algorithms for the MOTSP can be found in \cite{qamar2018comparative}.

In addition to multiple objectives, finding optimal tours that pass through a subset of nodes (as opposed to visiting all nodes) in the graph, also known as the Steiner Traveling Salesperson Problem (STSP), is particularly useful in the context of logistics \citep{letchford2013compact, zhang2015steiner}. One can transform a single-objective STSP instance into a TSP instance by generating a complete graph using the shortest paths between all terminal pairs \citep{alvarez2019note}. The nodes in the complete graph represent the terminals, and the edges are assigned the shortest path costs. Noting that generating a complete graph can be computationally expensive depending on the number of terminals, some studies directly solve the STSP using heuristics and MIP formulations. For instance, \cite{letchford2013compact} proposed a single commodity flow-based MIP formulation, and \cite{interian2017grasp} proposed a GRASP heuristic based on the 2-opt move for the symmetric STSP. However, for benchmark instances in the literature, \cite{alvarez2019note} found that the complete graph-based approach paired with the state-of-the-art TSP solvers such as Concorde \citep{applegate2003implementing} yielded faster results. 

In recent years, with the penetration of battery electric vehicles in last-mile logistics, the problem of TSP has become significantly more challenging due to factors such as charging locations, energy costs, hybrid vehicles, and battery capacity. For example, 
\cite{goeke2015routing} proposed an Adaptive Large Neighborhood Search (ALNS) to optimize the routing of a mixed fleet of EV and conventional combustion vehicles with an energy consumption model incorporating speed, gradient, and cargo load distribution. \cite{roberti2016electric} proposed a MILP formulation and a three-phase heuristic combined with dynamic programming to explore different EV recharging policies. \cite{keskin2016partial} introduced a more practical version of the problem that allows partial recharging. They formulated a MILP and developed an ALNS to solve it efficiently. \cite{lin2021electric} explored the last-mile delivery problem with time windows, where EVs can be charged or discharged at any station under time-variant electricity prices. The reader may refer to \cite{erdelic2019survey} for a comprehensive review of EV logistics. In contrast to these studies, we model operations where the fleet starts with a fully charged battery with enough capacity to complete the tour, with intermediate recharging through regenerative braking. We also integrate the effects of road gradients on energy consumption. While not widely recognized, road gradients are crucial role in last-mile routing and fuel consumption. Most papers estimate road grades using the altitude difference between the customer locations \citep{goeke2015routing, rao2016efficient}. However, \cite{brunner2021vehicle} observed that such simplification may induce errors in energy estimates when the altitude varies significantly along the route between two locations.   

To our knowledge, solution methods for the asymmetric BSTSPTW have not been explored, which, as discussed in Section \ref{sec:intro}, is an important problem in real-world EV freight routing. Reducing STSP to TSP is ineffective when time windows are involved (as will be elaborated in Section \ref{sec:heuristics}). Existing STSP formulations such as \cite{letchford2013compact} and \cite{ interian2017grasp} solve single-criterion problems without time windows. The only STSP heuristic we know of uses 2-opt moves and is limited to symmetric instances \citep{interian2017grasp}. Further, we apply these methods to a problem in EV logistics to demonstrate the practical advantages of this study.

\section{Preliminaries}
\label{sec:preliminaries}
\begin{figure}[H]
    \centering
    \begin{subfigure}{0.33\textwidth}
        \centering
        \includegraphics[scale=0.25]{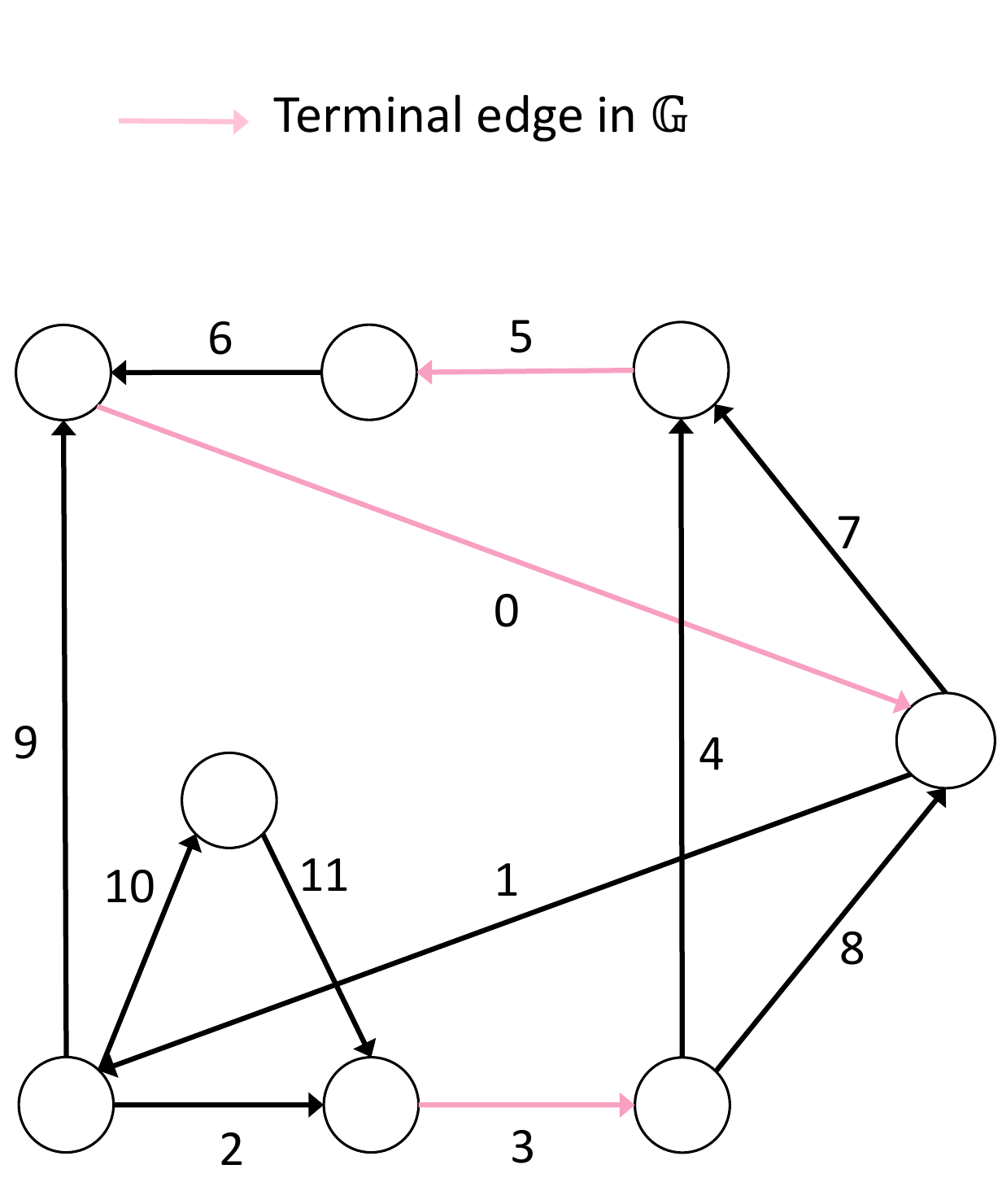}
        \caption{Toy road network $\primalgraph$}
        \label{fig:tourtosol}
    \end{subfigure}%
    \hfill
    \begin{subfigure}{0.36\textwidth}
        \centering
           \includegraphics[scale=0.25]{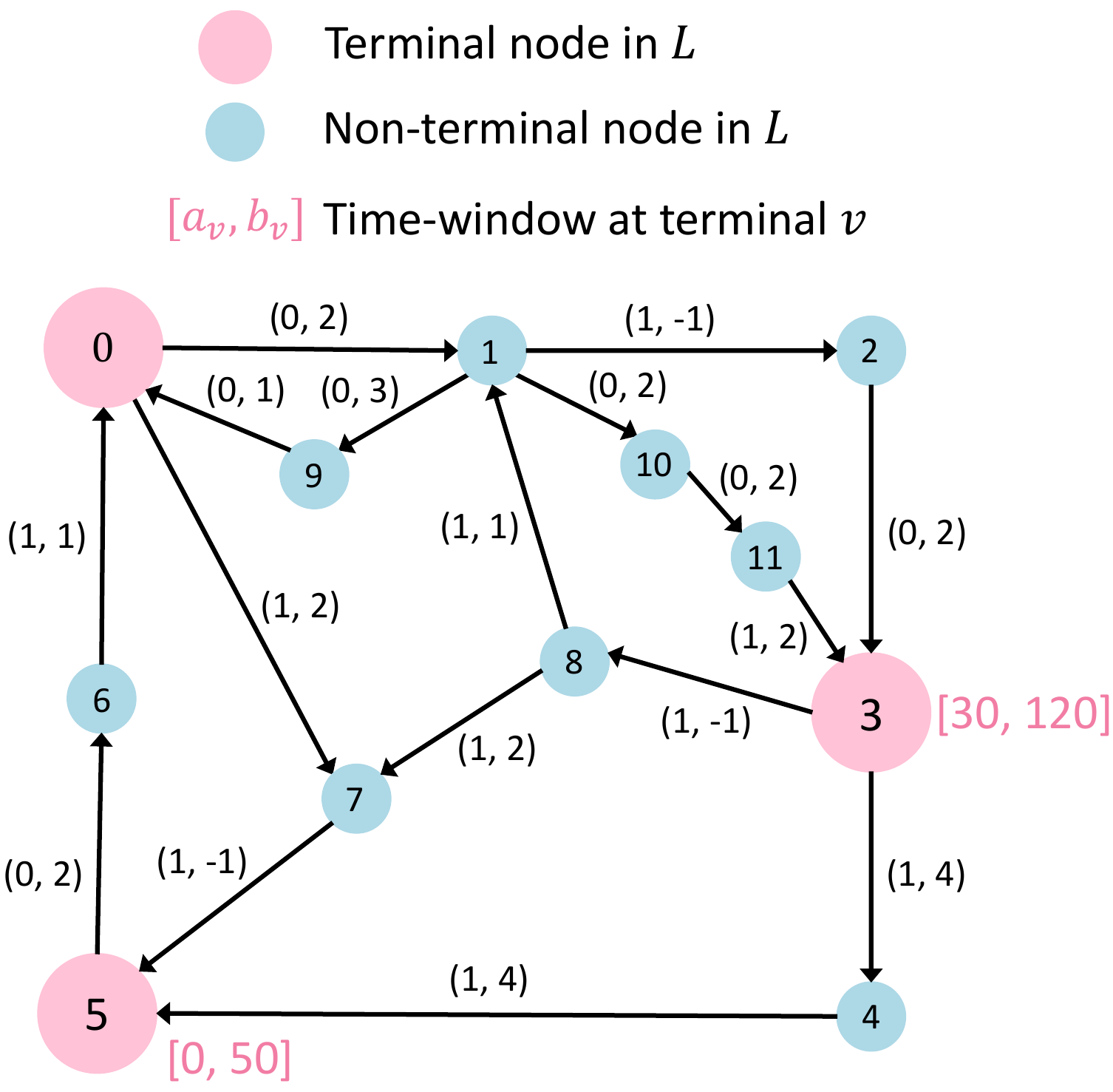}
        \caption{Line graph $\dualgraph$. Edge weights are left turns and energy}
        \label{fig:soltotour}
    \end{subfigure}%
    \begin{subfigure}{0.3\textwidth}
        \centering
           \includegraphics[scale=0.25]{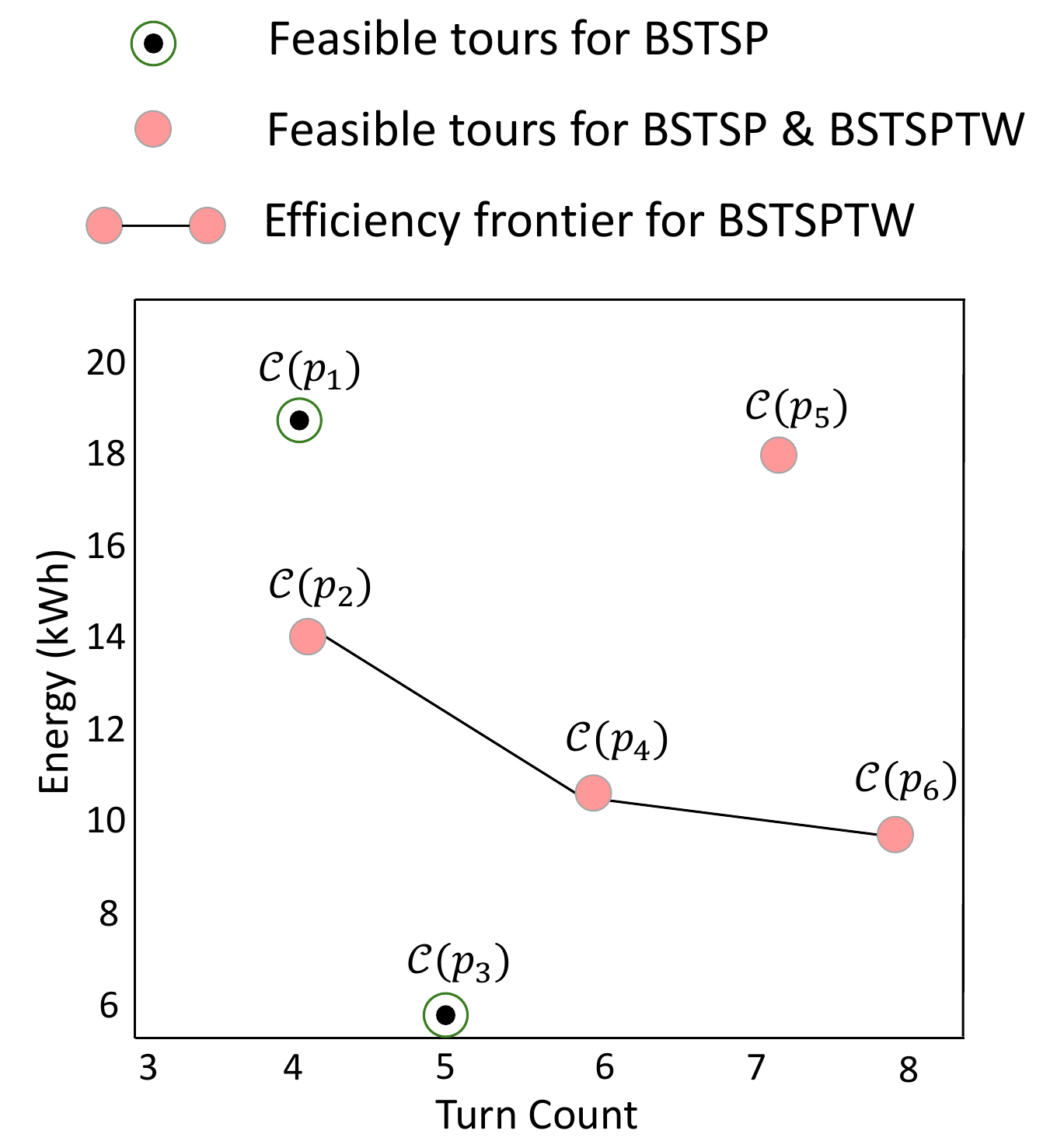}
        \caption{Pareto-optimal tours}
        \label{fig:soltotour}
    \end{subfigure}%
    \caption{Illustration of the BSTSP and BSTSPTW}
    \label{img:roadnet}
\end{figure}

This section formally introduces the BSTSPTW and the terminology used in this paper. The BSTSPTW aims to find the set of Pareto-optimal tours based on two objectives -- the number of left-hand turns and the energy consumption -- starting and ending at a depot while visiting customers within their time windows. Both objectives come with their unique set of challenges. First, the number of turns at a junction depends on three nodes (or two edges) and is not an edge attribute. Thus, one cannot directly use single- or bi-criteria shortest-path algorithms, such as Martin's algorithm \citep{martins1984multicriteria} to optimize turn movements. For example, see Figure \ref{img:roadnet}a, which shows a toy road network, where the labels indicate edge IDs and the pink edges are terminals that must be visited. A traveler moving along edge $0$ can choose to turn right to edge $1$ (the number of turns is zero) or make a left turn and move along edge $7$ (the number of turns is one). To keep the illustration simple, we ignore conflicts and label turns left and right based only on the angles between the edges. Second, many shortest-path algorithms assume non-negative edge weights, but since we model the energy consumed by vehicles with regenerative braking, the edge costs can be negative.

Let $\primalgraph = (\primalvertex, \primaledge, \primalterminal)$ be a directed graph, where $\primalvertex$ and $\primaledge$ denote the set of physical nodes and edges, respectively. Since the customer locations are often on the edges and not at the nodes, we let $\primalterminal \subseteq \primaledge$ denote the set of required edges that must be visited, also referred to as \textit{terminals}. The cardinality of $\primalterminal$ is indicated by $\sizedualterminal$. To address the turn issue, we reformulate the problem on a \textit{line graph}. Given a graph $\primalgraph$, the corresponding line graph $\dualgraph = (\dualvertex, \dualedge, \dualterminal)$ can be constructed by modeling the edges and turning movements in $\primalgraph$ as nodes $\dualvertex$ and edges $\dualedge$, respectively. The set of terminal edges $\primalterminal$ in $\primalgraph$ is transformed to a set of terminal nodes $\dualterminal$ in $\dualgraph$.\footnote{We could denote the line graph as $\dualgraph(\primalgraph)$ to highlight its dependence on the physical graph, but we refer to it as $\dualgraph$ for brevity.} Nodes in $\dualgraph$ are denoted using the symbols $u$ and $v$. Depot edge (node) in $\primalgraph$ ($\dualgraph$) is marked $0$. For example, consider Figure \ref{img:roadnet}b, which corresponds to the line graph of the network in Figure \ref{img:roadnet}a. Edge $1$ in the original graph is labeled as node $1$ in the line graph. The pink nodes in $\dualgraph$ represent the corresponding terminals. When traveling from $0$ to $7$ in $\primalgraph$, the turning direction depends on both edges $0$ and $7$. However, in the line graph $\dualgraph$, it becomes an edge property and depends on the edge connecting nodes $0$ and $7$.

The energy attributes in the line graph are defined assuming homogeneous edges in $\primalgraph$ (that is, the energy per unit length is the same along the entire edge) and that the customers are located at the midpoints of the edges (when customers are close to the endpoints, edges can be split into smaller segments to model energy usage more accurately). Mathematically, the energy consumption on a line graph's edge is set to the average energy required to traverse the corresponding edges in the original graph. For example, in Figure \ref{img:roadnet}b, the energy attribute on the edge connecting nodes $0$ and $7$ in $\dualgraph$ is set to $\frac{1}{2}(\energy{0} + \energy{7})$, where $\energy{0}$ and $\energy{7}$ denote the energy required to traverse the edges $0$ and $7$ in $\primalgraph$. Using a similar logic, One can derive the travel times on the edges in $\dualgraph$. For the remainder of this paper, unless explicitly stated, nodes and edges refer to those in the line graph. 
 
We now define additional notation for the line graph $L$. Let a \textit{path} be a set of nodes visited in a particular order. A path starting and ending at the same node is called a \textit{cycle}. A \textit{tour} is a cycle that visits each terminal and starts and ends at the depot. We represent the set of tours by $\tourset$. The cost attributes of an edge $(u,v) \in \dualedge$ are represented as an ordered pair $(\turn{u,v}, \energy{u,v})$, where $\turn{u,v}$ is the number of left turns and $\energy{u,v}$ denotes energy consumption. Likewise, for a path $p$, we define its \textit{cost} $\objectivevector{p}$ as a vector $(\turn{p}, \energy{p})$, where $\turn{p}$ and $\energy{p}$ denote the total number of left turns and the total energy required to traverse $p$, respectively.

\begin{definition}
    Given tours $p, q \in \tourset$, $p$ dominates $q$ (denoted by $\objectivevector{p} \preceq \objectivevector{q}$), iff $\turn{p} \leq \turn{q}$ and $\energy{p} \leq \energy{q}$.
\end{definition}

\begin{definition}
A set of tours $\paertooptimalset \subseteq \tourset$ is called \textit{Pareto-optimal} or \textit{non-dominated} iff for each tour $p\in \paertooptimalset$, there is no tour $q \in \tourset$, $q \neq p$ that dominates $p$. 
\end{definition}

The \textit{objective/criterion space} contains the cost vectors $\objectivevector{p}$ for all $p \in \tourset$. Let $\paertooptimalset$ denote the set of all Pareto-optimal tours for the BSTSPTW.  An \textit{efficiency frontier} $\efffrontier{\paertooptimalset} = \{\objectivevector{p}:p \in \paertooptimalset\}$ refers to the collection of points $\objectivevector{p}$ in the objective space for Pareto-optimal paths $p \in \paertooptimalset$. The methods proposed to solve the BSTSPTW maintain a non-dominated set $\localsearchndset$ that approximates $\paertooptimalset$. Figure \ref{img:roadnet}c shows the efficiency frontier for the line graph \ref{img:roadnet}b. Assuming constant travel times of $10$ minutes for each edge in the line graph, the connected pink points show the BSTSPTW efficiency frontier. The black points with green circles represent BSTSP tours (they violate time windows).

Table \ref{tab:extour} provides a detailed description of the tours. In general, node and edge revisits may be necessary in the line graph to discover optimal tours. To illustrate the need for revisits, a feature unique to our problem, consider terminal $5$. Assume that we start at the depot 0 at time $t=0$. Tour $p^1$ reaches terminal $5$ after 60 minutes, which does not satisfy its time window. Since terminal 5 must be reached in 50 minutes, the viable paths are $[\textbf{0}, 7, \textbf{5}]$ or $[\textbf{0}, 1, 2, \textbf{3}, 4, \textbf{5}]$. Due to the network structure, the subpath $[\textbf{0}, 7, \textbf{5}]$ forces a revisit to the depot node. Furthermore, in the Pareto-optimal tour $p^6$, visiting terminal $3$ requires revisiting both $5$ and $6$, and the edge $(5, 6)$. In fact, if the travel time on edge $(4, 5)$ is increased to 20 minutes while keeping the rest of the network the same, all BSTSPTW tours will require revisits.

\begin{table}[H]
\centering
\caption{\label{tab:extour} BSTSP and BSTSPTW tours for the network in Figure \ref{img:roadnet}. Terminals are highlighted in bold. By default, the depot is visited twice. A tour is considered to have revisits to the depot if it is visited at least three times.}
\begin{tabular}{lllllll}
\toprule
\multirow{2}{*}{\textbf{Tour}} & \multirow{2}{*}{\textbf{Node order}} & \multirow{2}{*}{\textbf{Turns}} & \multirow{2}{*}{\textbf{Energy}} &
\multirow{2}{*}{\textbf{TW feasible}} & \multicolumn{2}{l}{\textbf{Revisits}} \\
 &    &    &   &  & Node & Edge \\
 \midrule
$p^1$ & $[\textbf{0}, 1, {10}, {11}, \textbf{3}, 4, \textbf{5}, 6, \textbf{0}]$ & 4 & 19 & No & $\times$ & $\times$\\
$p^2$ & $[\textbf{0}, 1, 2, \textbf{3}, 4, \textbf{5}, 6, \textbf{0}]$ & 4 & 14 & Yes &$\times$ & $\times$ \\
$p^3$ & $[\textbf{0}, 1, 2, \textbf{3}, 8, 7, \textbf{5}, 6, \textbf{0}]$ & 5 & 6 & No & $\times$ & $\times$\\
$p^4$ & $[\textbf{0}, 7, \textbf{5}, 6, \textbf{0},  1, 2, \textbf{3}, 8, {1}, 9, \textbf{0}]$ & 6 & 11 & Yes & $\checkmark$ & $\times$\\
$p^5$ & $[\textbf{0}, 7, \textbf{5}, 6, \textbf{0}, 1, 2, \textbf{3}, 4, \textbf{5}, 6, \textbf{0}]$ & 7 & 18 & Yes & $\checkmark$ & $\checkmark$\\
$p^6$ & $[\textbf{0}, 7, \textbf{5}, 6, \textbf{0}, 1, 2, \textbf{3}, 8, 7, \textbf{5}, 6, \textbf{0}]$ & 8 & 10 & Yes & $\checkmark$ & $\checkmark$\\
\bottomrule
\end{tabular}
\end{table}

Formally, we define the BSTSP and BSTSPTW as follows. Given a graph $\dualgraph$, the objective of the BSTSP is to find a set of Pareto-optimal tours that start and return to the depot after visiting all the terminal nodes in $\dualterminal$ (at least once). The BSTSPTW additionally requires that each terminal $v \in \dualterminal$ is visited within its prescribed time window $[\earliesttime{v}, \latesttime{v}]$, where $\earliesttime{v}$ and $\latesttime{v}$ are measured with reference to the start time at depot 0. We ignore service times at terminals, but they can be trivially incorporated. Table \ref{tab:notn} summarizes the notation defined thus far.

\begin{table}[t]
\centering
\caption{Glossary of problem data}
\label{tab:notn}
\begin{tabular}{l l}
\toprule
\textbf{Notation} & \textbf{Description} \\ 
\hline
$\primalgraph = (\primalvertex, \primaledge, \primalterminal)$ & Road network with node set $\primalvertex$, edge set $\primaledge$, and terminal edge set $\primalterminal$ \\
$\dualgraph = (\dualvertex, \dualedge, \dualterminal)$ & Line graph with node set $\dualvertex$, edge set $\dualedge$, and terminal node set is $\dualterminal$ \\
$\sizedualterminal$ & Cardinality of $\dualterminal$\\ 
$v,p$ & Indices for nodes and paths/tours, respectively \\
$\objectivevector{p} = [\turn{p}, \energy{p}]$ & Cost vector, i.e., number of left turns and energy for path $p$ \\
$\energy{u,v}$ & Energy parameter of edge $(u,v)\in \dualedge$ \\
$\turn{u,v}$ & Turn count of edge $(u,v) \in \dualedge$ \\
$[\earliesttime{v}, \latesttime{v}]$ & Time window for terminal $v \in \dualterminal$ \\
$\paertooptimalset$ & Set of all Pareto-optimal tours for BSTSPTW \\
$\localsearchndset$ & Set of tours that approximate $\paertooptimalset$   \\
\hline
\end{tabular}
\end{table}

\section{Mixed Integer Programming Formulation}
\label{sec:exact}
This section presents a new MIP formulation for the BSTSPTW using decision variables that capture node and edge revisits. An optimal tour can revisit a node at most $\sizedualterminal$ times because revisiting a node is required only when there are terminals left to visit.\footnote{Imagine a Y-shaped or a hub-and-spoke network in which all terminals are the endpoints of cul-de-sacs. If the central node is the depot, it will be visited $\sizedualterminal$ times} When revisits occur, the time window constraint must be satisfied for at least one of the revisits. The proposed formulation addresses these constraints and allows for the identification of Pareto-optimal tours of the BSTSPTW. Our MIP builds on the work by \cite{letchford2013compact} in two ways: (1) time-window constraints (which are not the same as those of the classic TSPTW) are added while allowing revisits, and (2) the bi-objective problem is solved using scalarization \citep{geoffrion1968proper}. To formulate the MIP, we first create $\sizedualterminal$ copies of each node.\footnote{The number of node copies, decision variables, and constraints can be reduced through additional pre-processing. Since multiple revisits to a node $u$ can delay the time at which terminals are served, we can bound the maximum number of revisits by considering the shortest path time from the depot to $u$, the minimum travel time along a cycle that starts and ends at $u$ while visiting at least one terminal, and the latest times for serving other terminals, i.e., the $\latesttime{v}$ values.} It is assumed that a traveler starts at the depot with $\sizedualterminal - 1$ units of a commodity and drops exactly one unit at one of the copies of each terminal (other than the depot). The primary binary decision variable is $\edgevar{u}{i}{v}{j}$, where $u$ and $v$ represent nodes and $i$ and $j$ are their respective copy indices. Table \ref{tab:ipnotn} summarizes additional notation used in the MIP formulation.

\begin{table}[H]
\centering
\caption{MIP Terminology}
\label{tab:ipnotn}
\begin{tabular}{p{2.5cm}p{14cm}}
\hline
\textbf{Notation} & \textbf{Description} \\ 
\hline
$\edgevar{u}{i}{v}{j}$ & Binary variable which is 1 if the edge $(u,v)$ is part of the optimal tour, where node $u$'s $i^{th}$ copy is visited followed by node $v$'s $j^{th}$ copy, and is $0$ otherwise.\\
$\arrivaltimevar{v}{i}$ & Departure time variable at $i^{th}$ copy of node $v$\\
$\twnodevar{v}{i}$ & Binary variable which is $1$ if terminal $v$'s  $i^{th}$ copy is visited within its time windows, and is $0$ otherwise.\\
$\commodityflow{u}{i}{v}{j}$ & Number of commodities that pass through the edge $(u,v)$ when the $i^{th}$ copy of $u$ and the $j^{th}$ copy of $v$ are visited in succession. \\
$\alpha, \beta$ & Scalarization parameters for MIP \\
$\traveltime{u,v}$ & Travel time on edge $(u,v) \in \dualedge$ \\
$M_{u,v}, \underline{M}_{v}, \overline{M}_{v}$ & Big-M constants \\

\hline
\end{tabular}
\end{table}
The number of revisits to a node $v$ is given by $\sum_{(v, w) \in \dualedge} \sum_{i=1}^{\sizedualterminal} \sum_{j=1}^{\sizedualterminal} \edgevar{v}{i}{w}{j}$ and the number of revisits to edge $(u,v)$ is $\sum_{i=1}^{\sizedualterminal}\sum_{j=1}^{\sizedualterminal} \edgevar{u}{i}{v}{j}$. The objective function \eqref{eq:eq1} minimizes the weighted sum of two terms: the number of left turns and the energy. While the energy attribute can be negative on a subset of edges, we assume the problem is bounded and has no negative energy cycles. The weights of these parameters are represented by $\alpha$ and $\beta$, respectively. Equation \eqref{eq:eq3} ensures tour continuity every time a node copy is visited. Equations \eqref{eq:eq4} and \eqref{eq:eq5} guarantee that exactly one unit of the commodity is delivered at each terminal, and no deliveries are made at non-terminals across their copies (i.e., even if multiple copies of a terminal are visited, commodities are dropped off only once). 

\begin{align}
     \min  & \sum_{(u, v) \in \dualedge} \left(\alpha  \turn{u,v} + \beta  \energy{u,v} \right) \sum_{i=1}^{\sizedualterminal}\sum_{j=1}^{\sizedualterminal} \edgevar{u}{i}{v}{j} \label{eq:eq1}\\
      \text{s.t. } & \sum_{(u, v) \in \dualedge} \sum_{j=1}^{\sizedualterminal} \edgevar{u}{j}{v}{i} = \sum_{(v, w) \in \dualedge}\sum_{k=1}^{\sizedualterminal} \edgevar{v}{i}{w}{k}  && \forall~v \in \dualvertex, i \in \{1, \ldots, \sizedualterminal\}   \label{eq:eq3}\\
                 & \sum_{(u, v) \in \dualedge} \sum_{i=1}^{\sizedualterminal}\sum_{j=1}^{\sizedualterminal}  \commodityflow{u}{i}{v}{j}  - \sum_{(v, w) \in \dualedge} \sum_{j=1}^{\sizedualterminal}\sum_{k=1}^{\sizedualterminal}  \commodityflow{v}{j}{w}{k} = 1 &&   \forall~v \in \dualterminal \setminus \{\depot\}  \label{eq:eq4}\\
                 & \sum_{(u, v) \in \dualedge} \sum_{i=1}^{\sizedualterminal}\sum_{j=1}^{\sizedualterminal}  \commodityflow{u}{i}{v}{j}  - \sum_{(v, w) \in \dualedge} \sum_{j=1}^{\sizedualterminal}\sum_{k=1}^{\sizedualterminal}  \commodityflow{v}{j}{w}{k} = 0 &&  \forall~v \in \dualvertex \setminus \dualterminal \label{eq:eq5}\\
                 & 0 \leq \commodityflow{u}{i}{v}{j} \leq (\sizedualterminal - 1) \edgevar{u}{i}{v}{j}                             && \forall~(u, v) \in \dualedge,~i,j \in \{1, \ldots, \sizedualterminal\}   \label{eq:eq6}\\
                 & \arrivaltimevar{u}{i} + \traveltime{u,v} - M_{u,v}(1 - \edgevar{u}{i}{v}{j}) \leq \arrivaltimevar{v}{j} && \forall~(u,v) \, \in\dualedge,~i,j \in \{1, \ldots, \sizedualterminal\},~(v,j) \neq (0,1)  \label{eq:eq7} \\ 
                 & \sum_{(\depot,v) \in \dualedge}\sum_{i=1}^{\sizedualterminal}\edgevar{\depot}{1}{v}{i} = 1 &&   \label{eq:eq7.1} \\ 
                 & \twnodevar{v}{i} \leq \sum_{(v, w) \in \dualedge}\sum_{j=1}^{\sizedualterminal} x_{v,i,w,j} && \forall~v \in \dualterminal,~i \in \{1, \ldots, \sizedualterminal\}  \label{eq:eq8} \\
                 & \sum_{i=1}^{\sizedualterminal} \twnodevar{v}{i} \geq 1 && \forall~v \in \dualterminal \label{eq:eq9} \\
                 & \earliesttime{v} - \underline{M}_{v}(1 - 
                 \twnodevar{v}{i}) \leq \arrivaltimevar{v}{i} \leq  \latesttime{v} + \overline{M}_{v}(1 - \twnodevar{v}{i}) && \forall~v \in \dualterminal,~i \in \{1, \ldots, \sizedualterminal\}  \label{eq:eq10} \\
                 & \arrivaltimevar{v}{i} \in \mathbb{R_+},  \twnodevar{v}{i} \in \{0,1\}   && \forall~v \in \dualvertex,~i \in \{1, \ldots, \sizedualterminal\}  \label{eq:eq11}\\  
                & x_{u,i,v,j} \in \{0,1\} && \forall~(u, v) \in \dualedge,~i,j \in \{1, \ldots, \sizedualterminal\}  \label{eq:eq12}
\end{align}

 Constraint \eqref{eq:eq6} specifies that commodities can only pass through edges that are in the tour and bounds the decision variable $\commodityflow{u}{i}{v}{j}$. Equations \eqref{eq:eq4}--\eqref{eq:eq6} together eliminate subtours. Constraint \eqref{eq:eq7} ensures that departure times increase along any path within a tour and, consequently, avoids multiple visits to the same copy of the terminal. It is applied at all terminal copies except the first copy of the depot since the tour must terminate at it. Constraint \eqref{eq:eq7.1} extends the restriction on multiple visits to the first copy of the depot. The MIP implicitly determines the start time at node $\depot$. Inequalities \eqref{eq:eq8}--\eqref{eq:eq10} ensure that time windows are imposed for at least one of the copies of each terminal. In particular, according to \eqref{eq:eq8}, $\sum_{(v,w) \in \dualedge} \sum_{j=1}^{\sizedualterminal} \edgevar{v}{i}{w}{j}$ takes a value one when the $i^{th}$ copy of terminal $v$ is visited. In all other cases, the $\twnodevar{v}{i}$ values are forced to zero. Constraint \eqref{eq:eq9} guarantees that $\twnodevar{v}{i}$ is set to one for at least one of terminal $v$'s copies. Constraint \eqref{eq:eq8} and \eqref{eq:eq9} ensure that the traveler leaves each terminal at least once across all its copies. Whenever $\twnodevar{v}{i}$ is one, \eqref{eq:eq10} ensures that the time window at terminal $v$ is satisfied for its $i^{th}$ copy. For example, in tour $p^6=[\textbf{0}, 7, \textbf{5}, 6, \textbf{0}, 1, 2, \textbf{3}, 8, 7, \textbf{5}, 6, \textbf{0}]$ from Table \ref{tab:extour}, the first and second visits to node $5$ occur $20$ and $100$ minutes after departure from the depot node, i.e., $\arrivaltimevar{5}{1}=20$, $\arrivaltimevar{5}{2}=100$. Since the time-window conditions at node $5$ are satisfied only for the first visit, \eqref{eq:eq10} forces $\twnodevar{5}{1}$ to one and $\twnodevar{5}{2}$ and $\twnodevar{5}{3}$ to zeros. Waiting is modeled implicitly in Constraint \eqref{eq:eq10} since $\arrivaltimevar{v}{i}$ represents the departure time. Constraints \eqref{eq:eq11}--\eqref{eq:eq12} set non-negativity and binary restrictions on the decision variables. 

To set the big-$M$ values, we first calculate each node's earliest and the latest possible departure times. The earliest possible departure time at a node is defined using the shortest path duration from the depot to that node since this is the earliest we can depart from it. The latest possible departure time is the same across all nodes and is defined as the sum of $\max_{v\in\dualterminal}{\latesttime{v}}$ and the duration of the longest path among the Pareto-optimal paths from any terminal to the depot. From \eqref{eq:eq7}, note that $M_{u,v} \geq \arrivaltimevar{u}{i} + \traveltime{u,v} - \arrivaltimevar{v}{j}$. Hence, an upper bound for the value of $M_{u,v}$ is the sum of the latest possible departure time at node $u$, the travel time between the nodes $u$ and $v$, and the negative of the earliest possible departure time from node $v$. Likewise, we determine $\overline{M}_{v}$ by subtracting the latest time window at node $v$ (i.e., $\latesttime{v}$) from the latest possible departure time at node $v$; and $\underline{M}_{v}$ is set to the difference between $\earliesttime{v}$ and the earliest possible departure time from node $v$.

Although node copies in the formulation capture revisits, an initial visit to node $v$ may have a higher copy index than a subsequent visit. In other words, should a node be visited on two separate occasions, the index associated with its copies could be any two integers in $\{1, \ldots, \sizedualterminal\}$ and not necessarily copy 1 and copy 2. For example, consider Figure \ref{fig:soltotour}, which shows the set of $\edgevar{u}{i}{v}{j}$ variables that satisfy the MIP constraints and are equal to one for the tour $p^6=[\textbf{0}, 7, \textbf{5}, 6, \textbf{0}, 1, 2, \textbf{3}, 8, 7, \textbf{5}, 6, \textbf{0}]$ from Table \ref{tab:extour}. The numbers next to the nodes represent the arrival times. Here, the second visit to terminal $5$ does not happen at the second copy but is instead connected to the third copy.

\begin{figure}[H]
    \centering
        \begin{subfigure}{0.46\textwidth}
        \centering
           \includegraphics[width=\linewidth,page=2]{Journal_Images/mipsol.pdf}
        \caption{Constructing tours from a general solution}
        \label{fig:soltotour}
    \end{subfigure}%
    \hfill
    \begin{subfigure}{0.45\textwidth}
        \centering
        \includegraphics[width=\linewidth,page=1]{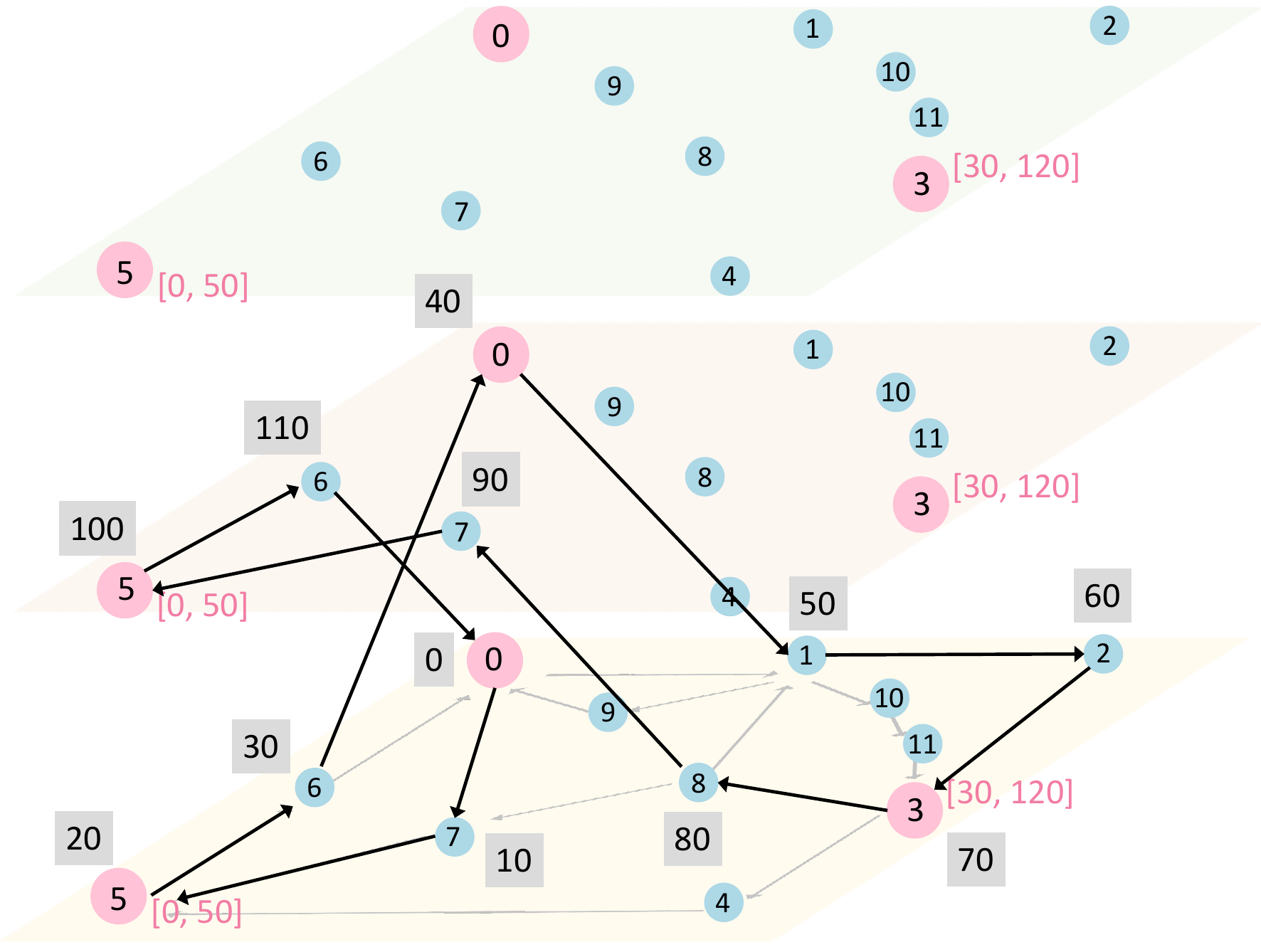}
        \caption{Tours where copies correspond to revisit order}
        \label{fig:tourtosol}
    \end{subfigure}%
    \caption{Interpreting solutions from the MIP formulation. The original edges are shown in the lowest layer. Values in grey boxes denote the arrival time.}
    \label{img:mipsol}
\end{figure}
 
To deal with this issue, a precedence order can be imposed by adding constraints of the form $\sum_{v:(u,v) \in \dualedge} \sum_{j = 1}^{\sizedualterminal} \edgevar{u}{i}{v}{j} \geq \sum_{v:(u,v) \in \dualedge} \sum_{j = 1}^{\sizedualterminal} \edgevar{u}{i+1}{v}{j}$ for all $ u \in \dualvertex \setminus \{0\}, i \in \{1, \ldots, \sizedualterminal - 1\}$. These inequalities require a node in a particular layer to be visited only if its copies in the lower layers are visited. While the inequalities are valid for the above formulation and may be added as lazy cuts, they can increase the number of constraints. They would, however, result in solutions similar to the one shown in Figure \ref{fig:tourtosol}. The node copies in the third layer remain unused since no node in the network is visited three times. In our experiments, we noticed that adding these extra constraints slows the MIP code. Instead, using a post-processing step, any solution to the MIP can be used to derive a feasible BSTSPTW tour simply by following the $x$ variables starting from the first copy of the depot.

\begin{figure}[H]
    \centering
    \includegraphics[scale=.45]{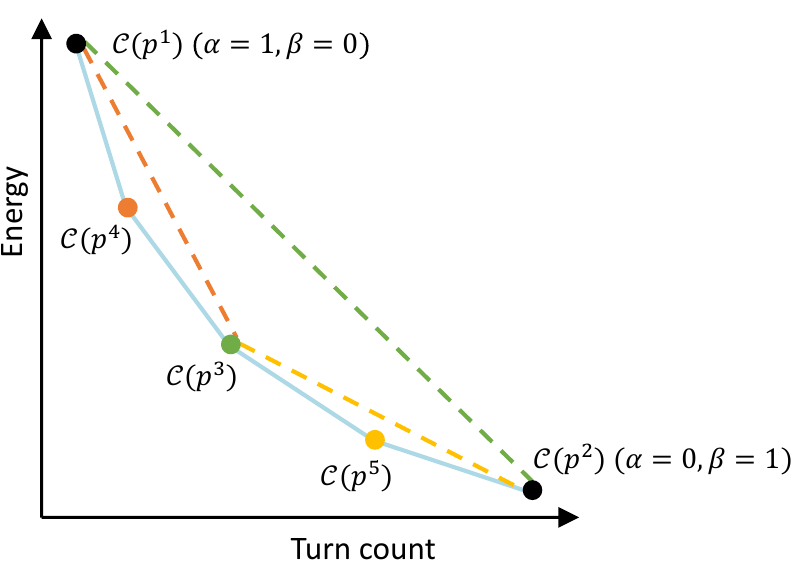}
    \caption{Scalarization technique}     
    \label{fig:scl}
\end{figure}

We now describe \textit{scalarization}, a standard technique to enumerate points on the efficiency frontier of bi-criteria optimization problems. Let the optimal tour for fixed $\alpha$ and $\beta$ be given by \textsc{MixIntProg}($\alpha, \beta$). We first start by finding the extreme points, which correspond to the cases where $(\alpha, \beta)$ is $(1,0)$ and $(0,1)$. These two cases are single-objective STSPs, where only the number of turns or energy is optimized. For subsequent iterations, we use a recursive algorithm, which sets the ratio of the weights (i.e., $\alpha / \beta$) to the absolute value of the slope of the line connecting the endpoints of unexplored regions of the efficiency frontier. For example, in Figure \ref{fig:scl}, suppose the solutions to the single-criterion versions obtained by setting $(\alpha, \beta)$ to (1,0) and (0,1) are $p^1$ and $p^2$, respectively. We invoke \textsc{MixIntProg} using the absolute value of the slope of the line joining $\objectivevector{p^1}$ and $\objectivevector{p^2}$, which results in a new tour $p^3$.  Similarly, in the next iteration, the recursive algorithm uses the slope between $\objectivevector{p^1}$ and $\objectivevector{p^3}$ to discover tour $p^4$ and the slope between $\objectivevector{p^3}$ and $\objectivevector{p^2}$ to find tour $p^5$. Recursion continues until no new path is found. 

Algorithm \ref{alg:IP_sol} (\textsc{Scalarization}) summarizes the pseudocode for this process. Note that scalarization always yields a convex-shaped efficiency frontier as it uses a convex combination of objectives \citep{ marler2010weighted}. However, the efficiency frontier for multi-criteria problems is not necessarily convex \citep{climaco2016approach}. Hence, the algorithm may return a subset $\localsearchndset$ of the complete Pareto-optimal tour set $\paertooptimalset$. Lines 1--2 call \textsc{MixIntProg} to get the extreme points. Line 6 adds the newly discovered tours to $\localsearchndset$. Line 7 calls the recursive function \textsc{NewTour}, which starts by setting the absolute value of the slope of the segment connecting the two tours. It then solves the MIP with normalized weights $\alpha = \omega/(1+\omega)$ and $\beta = 1/ 1+\omega$. Next, it checks if the termination criteria are satisfied, i.e., if the newly generated point coincides with the points from which the slope was initially derived or if the maximum time is exceeded. If not, we add the new tour to $\localsearchndset$, and the recursion continues.

\begin{algorithm}[H]
\caption{\textsc{Scalarization}}\label{alg:IP_sol}
\begin{algorithmic}[1]
\State $p \gets \textsc{MixIntProg}(\alpha = 1, \beta = 0)$ \Comment Optimize number of turns

\State $q \gets \textsc{MixIntProg}(\alpha = 0, \beta = 1)$ \Comment Optimize energy usage

\If {$\objectivevector{p} = \objectivevector{q}$} 
\State \textbf{return}
\Else
\State $\textit{Z} \gets \{p, q$\}
\State \textsc{NewTour}$(p, q, \localsearchndset)$ \Comment{Find new tours on the efficiency frontier}\\
\Return $\localsearchndset$
\EndIf

\vspace{2mm}
\setcounter{ALG@line}{0}
\Procedure{\textsc{NewTour}}{$p, q, \localsearchndset$}
\State $\omega \gets \left | (\energy{p} - \energy{q}) / (\turn{p} - \turn{q}) \right| $
\State $r \gets \textsc{MixIntProg} \left(\alpha = \omega / (1+\omega), \beta = 1 / (1+\omega) \right)$ \Comment{Solve the MIP}

\If {$\objectivevector{r} = \objectivevector{p}$ \textbf{or} $\objectivevector{r} = \objectivevector{q}$ \textbf{or} time limit exceeded} 
\State \textbf{return}
\EndIf
\State $\localsearchndset \gets \localsearchndset \cup \{r\}$ 
\State $\textsc{NewTour}(p, r, \localsearchndset)$
\State $\textsc{NewTour}(r, q, \localsearchndset)$
\EndProcedure

\end{algorithmic}
\end{algorithm}

\section{Local Search}
\label{sec:heuristics}
While Algorithm \ref{alg:IP_sol}, i.e., Scalarization combined with MIP, can discover exact BSTSPTW solutions, it does not scale well for large-scale real-world networks. This necessitates using heuristics to find near-optimal solutions in a reasonable time since the problem at hand is operational. 

\begin{figure}[H]
    \centering
\includegraphics[scale=0.35]{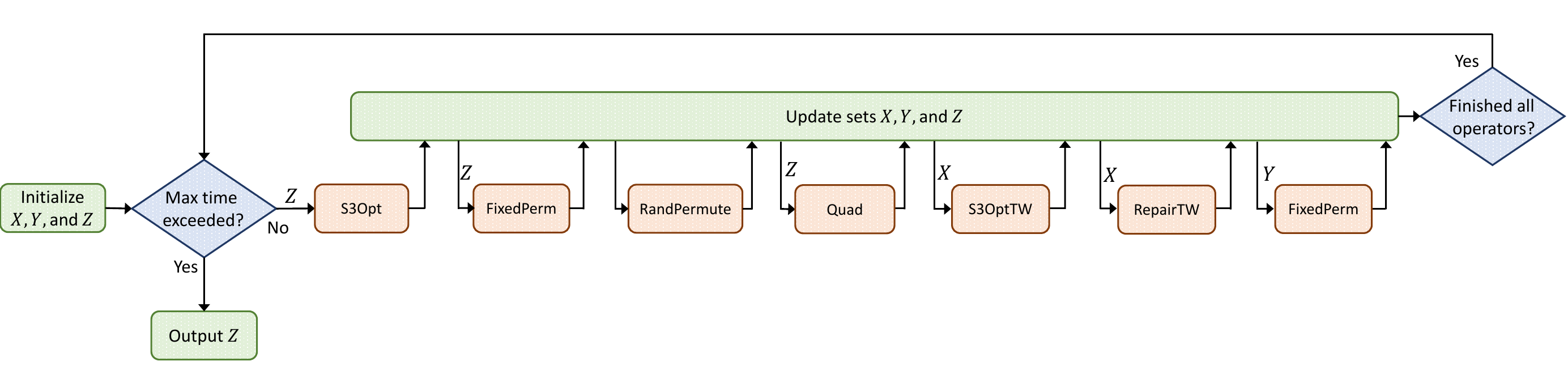}
    \caption{Framework of the proposed local search. Operators are shown in orange.}
    \label{fig:flowchart}
\end{figure}

In this section, we propose a new local search-based heuristic that starts with three sets: time-window feasible Pareto-optimal tours ($\localsearchndset$), time-window feasible non-Pareto-optimal tours ($\setY$), and tours that potentially violate the time windows ($\explorationset$). six operators are applied sequentially in each iteration (\textsc{FixedPerm} is applied twice), as shown in Figure \ref{fig:flowchart}. These operators select tours from the aforementioned sets and generate new tours, subsequently updating sets $\explorationset$, $\setY$, and $\localsearchndset$. To keep the figure compact, we represent the steps that update these sets in a single block. We repeat this procedure until a stopping criterion is satisfied. Section \ref{ls:trans} describes a procedure for generating initial tours to populate sets  $\explorationset$, $\setY$, and $\localsearchndset$. The operators are explained in Section \ref{ls:nbd}. Finally, Section \ref{ls:full} summarizes the complete local search procedure.

\subsection{Initial Solution Generation}\label{ls:trans}
\label{ls:init}
\begin{figure}[h]
    \centering
    \includegraphics[scale=0.47]{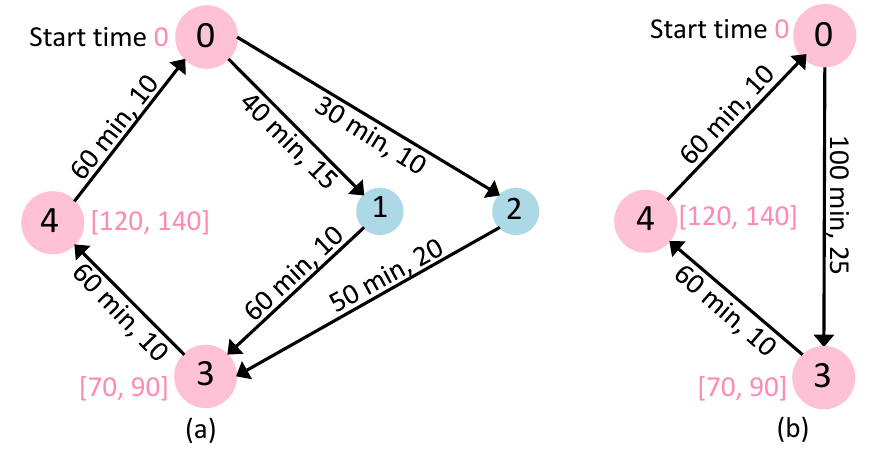}
    \caption{Issues with STSP to TSP transformations. (a) STSPTW instance (b) Corresponding TSPTW instance}
    \label{img:tranformation}
\end{figure}

As mentioned earlier, we can transform the single-objective STSP to a TSP by generating a complete graph where the terminals represent nodes, and edges are assigned costs based on the shortest paths between the terminals. When edge costs can take negative values, the shortest paths can be found using Johnson's algorithm \citep{johnson1977efficient}. Johnson's algorithm works in two stages. The first stage involves applying the Bellman-Ford algorithm to create a graph with positive edge costs. Subsequently, in the second stage, Dijkstra's algorithm is used to find the shortest path with modified costs. In the case of the BSTSP, we can combine the scalarization technique discussed in Section \ref{sec:exact} with single-objective STSPs and use the complete graph method described above to generate Pareto-optimal tours. However, the above technique fails for problems with time windows. When we convert an STSPTW to a TSPTW instance, the shortest-cost paths between terminals may lead to time-window infeasibilities since reducing makespan is not our objective. To see why, consider the example in Figure \ref{img:tranformation}a, which shows an STSPTW instance where the terminal nodes are pink. The edge attributes indicate the travel times and costs between nodes. Time windows are shown in square brackets next to each node. For simplicity, consider the single-objective case that minimizes cost. Figure \ref{img:tranformation}b represents the corresponding TSPTW instance. Observe that the least-cost path from nodes $0$ to $3$, i.e., $[\textbf{0}, 1, \textbf{3}]$, violates the time-window constraint at node $3$. However, in the STSPTW instance (Figure \ref{img:tranformation}a), there exists a time-window feasible tour $[\textbf{0}, 2, \textbf{3}, \textbf{4}, \textbf{0}]$. To address this issue, one could create a multigraph from the STSP instance with all possible paths connecting two terminals as edges, but this would result in a graph with potentially an exponential number of edges, making the TSP intractable.

To overcome these challenges, we propose to generate initial solutions for the local search using a modified version of \textsc{Scalarization} in Algorithm \ref{alg:IP_sol}. For a given set of edge-attribute weights ($\alpha$ and $\beta$), we first generate the corresponding single-criterion STSP instance. Next, we convert the STSP to a TSP instance by finding the shortest path between every terminal pair using Johnson's algorithm. If the resulting TSP instance contains negative edge weights, we increase all edge weights by the absolute value of the smallest edge weight. This is done because, to the best of our knowledge, state-of-the-art TSP solvers such as Lin-Kernighan-Helsgaun (LKH) do not support negative edge weights. Also, this transformation maintains correctness as every solution to a standard TSP (without revisits) has the same number of edges. Finally, we solve the TSP instance using the LKH solver. This process is repeated for a fixed time budget. The resulting tours are grouped into three sets $\localsearchndset$, $\setY$, and $\explorationset$ based on earlier definitions.

While updating these sets, we assess the time-window feasibility of a tour $p$ using a penalty function $\penaltyfunction{p}$, which provides an upper bound on the sum of the time-window penalties across all terminals in tour $p$, i.e., $\penaltyfunction{p} =  \sum_{v \in \dualterminal} \penaltyfunction{p, v}$.
Here, $\penaltyfunction{p, v}$ is the minimum penalty of terminal $v$ (considering all its revisits, indexed by $k$) in tour $p$, i.e., $\penaltyfunction{p, v} = \min_{k}\penaltyfunction{p, v, k}$, and $\penaltyfunction{p, v, k} = \max\big\{0, \patharrivaltime{p}{v,k} - \latesttime{v}\big\}$. The functions $\penaltyfunction{p, v, k}$ and $\patharrivaltime{p}{v,k}$ are the penalty and arrival time at terminal $v$ at its $k^{th}$ revisit, respectively, assuming that each terminal is served during its first visit. Since we allow waiting at terminal nodes, $\penaltyfunction{p, v, k}=0$ when $\patharrivaltime{p}{v, k} < \latesttime{v}$. Note that $\penaltyfunction{p}$ is an upper bound on the actual penalty of the tour. To understand why, consider a tour $p=[\textbf{0},\textbf{1},\textbf{2},\textbf{1},\textbf{0}]$, where the set of terminals is $\dualterminal=\{1,2,3,4,5\}$, and the time windows are defined as $[\earliesttime{1}, \latesttime{1}]=[30, 60]$ and $[\earliesttime{2}, \latesttime{2}]= [10, 30]$, with a travel time of 10 minutes between every pair of nodes. Suppose we assume that terminal 1 is served during its first visit. In that case, we have $\patharrivaltime{p}{1, 1}=10$ and $\patharrivaltime{p}{2, 1}=40$ (since there is a 20-minute wait at terminal 1), which implies that the time-window feasibility penalty $\penaltyfunction{p}=10$. However, if it is served during its second revisit, the tour is feasible for BSTSPTW. Calculating the actual tour penalty is computationally intensive as it requires multiple scans of the tour to check for the feasibility of terminals depending on when we serve them. As we will see shortly, $\penaltyfunction{p}$ is essential to every operator and must be invoked for every new path. Hence, a quick-and-dirty upper bound is better than the actual penalty that takes longer to compute. Table \ref{tab:LSnotn} summarizes some additional notation the local search operators use.

\begin{table}[H]
\centering
\caption{Local search terminology}
\label{tab:LSnotn}
\begin{tabular}{p{1.5cm}p{14cm}}
\hline
\textbf{Notation} & \textbf{Description} \\ 
    \hline
    \textbf{Variables/Functions} \\
    \hline
$p^{k}, q^{k}$ & General indices for paths and cycles\\
$\tourindex{p}{i}$ & $i^{th}$ node in path $p$\\
$\subpath{p}{i}{j} $ & Subpath from nodes $\tourindex{p}{i}$  to $\tourindex{p}{j}$ in path $p$ (including endpoints) \\ 
$\length{p},\length{\subpath{p}{i}{j}}$ & Number of edges in the path $p$, $\subpath{p}{i}{j}$, respectively \\
$\penaltyfunction{p, v, k}$ & Time-window penalty at the $k^{th}$ visit of terminal $v$ in tour $p$ \\
$\penaltyfunction{p, v}$ & Minimum (across all revisits) time-window penalty at terminal $v$ in tour $p$ \\
$\penaltyfunction{p}$ &  Time-window penalty of tour $p$ \\
$\patharrivaltime{p}{v,k}$ & Arrival time at the $k^{th}$ visit of terminal $v$ along tour $p$ \\
$\gain{\localsearchndset}{\subpathcollection}$ & A gain function which measures the efficiency of a set of subpaths $\subpathcollection$ w.r.t. to tours in $\localsearchndset$\\ 
$p \concatenate q$ & Path formed by concatenating two paths $p$ and $q$\\
$\meanenergy{\dualedge}$ & Average energy consumption across the edges in $\dualedge$\\ 
$\meanenergy{\localsearchndset}$ & Average energy consumption across the edges in tours belonging to $\localsearchndset$\\ 
$\meanenergy{\subpath{p}{u}{v}}$ & Average energy consumption across the edges in path $\subpath{p}{u}{v}$ \\
$\meanturns{\dualedge}$ & Average number of left turns across edges in $\dualedge$\\ 
$\meanturns{\localsearchndset}$ & Average number of left turns across the edges in tours belonging to $\localsearchndset$\\ 
$\meanturns{\subpath{p}{u}{v}}$ & Average number of left turns across the edges in path $\subpath{p}{u}{v}$\\
$\stdenergy$ & Standard deviation of energy consumption the edges in $\dualedge$\\ 
$\stdturns$ & Standard deviation of the number of left turns across the edges in $\dualedge$\\ 
$\hyperParameter{1}, \ldots, \hyperParameter{6}$& Hyperparameters used in the local search operators \\
    \hline
    \textbf{Sets} \\
    \hline
$\candidates$ & Set of candidates used in \textsc{S3Opt} and \textsc{S3OptTw} operators\\
$\explorationset$ & Set of BSTSP feasible tours that violate time windows for at least one terminal \\
$\setY$ & Set of BSTSPTW feasible, non-Pareto-optimal tours\\
$\adjterminalvector{p}$ & Set of tuples $(i,j)$ such that $\tourindex{p}{i}$ and $\tourindex{p}{j}$ are adjacent terminals in path $p$\\
$\violatedterminalset{p}$ & Set of terminals that do not satisfy time windows in path $p$ $ 
(\violatedterminalset{p}\subseteq\dualterminal\setminus\{\depot\})$\\
$\pathbank{u}{v}$ & Set of Pareto-optimal paths from terminal $u$ to $v$ that are reused from previous iterations \\
$\pathbanksubset{u}{v}$ & Set of paths (not necessarily Pareto-optimal) from terminal $u$ to $v$ \\
\hline
\end{tabular}
\end{table}

\subsection{Search Operators}
\label{ls:nbd}
This section introduces six operators: \textsc{S3Opt}, \textsc{S3OptTW}, \textsc{RepairTW},  \textsc{FixedPerm}, \textsc{Quad}, and \textsc{RandPermute}. Together, they help intensify and diversify solutions. Table \ref{tab:nbd} summarizes the characteristics of these operators.

\begin{table}[H]
\centering
\small
\caption{Summary of operators. Column \textit{Type} indicates if the operator helps in Intensification (Int) or Diversification (Div).}
\label{tab:nbd}
\begin{tabular}{p{2.2cm} p{0.8cm} p{0.8cm} p{10cm}}
\hline
\textbf{Name}& \textbf{Type} & \textbf{Input} & \textbf{Purpose}  \\
\hline
\textsc{S3Opt} & Int & $\localsearchndset$ & Find the best 3-opt move that improves the objectives   \\
\textsc{S3OptTW}& Int & $\explorationset$    & Find the best 3-opt move while fixing time-window infeasibilities \\
\textsc{RepairTW}& Int & $\explorationset$  & Repair time-window infeasibilities by removing and inserting terminals
  \\
\textsc{FixedPerm} & Int  & $\setY/\localsearchndset$   & Find new tours by optimizing the paths between the terminals   \\
\textsc{Quad} & Div & $\localsearchndset$ & Achieve diversification by changing the order of four terminal pairs \\
\textsc{RandPermute}& Div & -     & Explore a new random permutation of terminals \\
\hline
\end{tabular}

\end{table}

\subsubsection{\textsc{S3Opt} and \textsc{S3OptTW}}
\label{sssec:iopt}
\begin{figure}[H]
    \centering
    \includegraphics[scale=0.38]{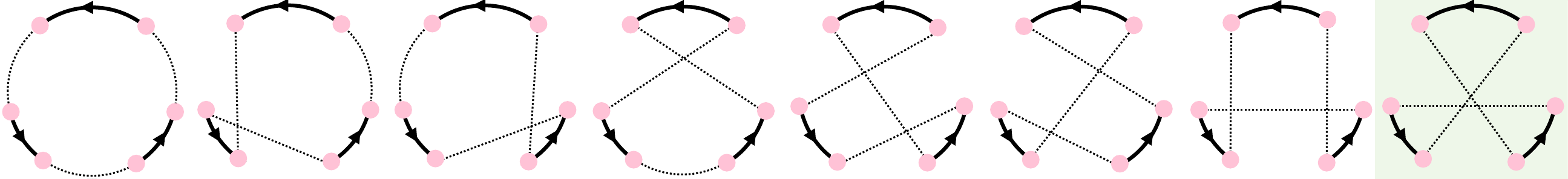}
    \caption{All possible swaps in 3-opt, only one of which preserves directionality}
    \label{fig:3opts}
\end{figure}

A conventional 3-opt move involves deleting three edges in a TSP tour and adding three edges to find a better tour. In the current situation, this approach is not effective for two reasons. First, in the STSP environment, the modified tour can skip terminal nodes. Second, due to the asymmetry of edges in the network, only one of the seven possible 3-opt moves preserves the direction of the edges (see Figure \ref{fig:3opts}). In all other cases, since the path segments are reversed, the number of left turns is significantly affected. To overcome this issue, we propose two improved Steiner versions of 3-opt -- the \textsc{S3Opt} and \textsc{S3OptTW} operators.

\begin{figure}[H]
    \centering
    \includegraphics[scale=0.47]{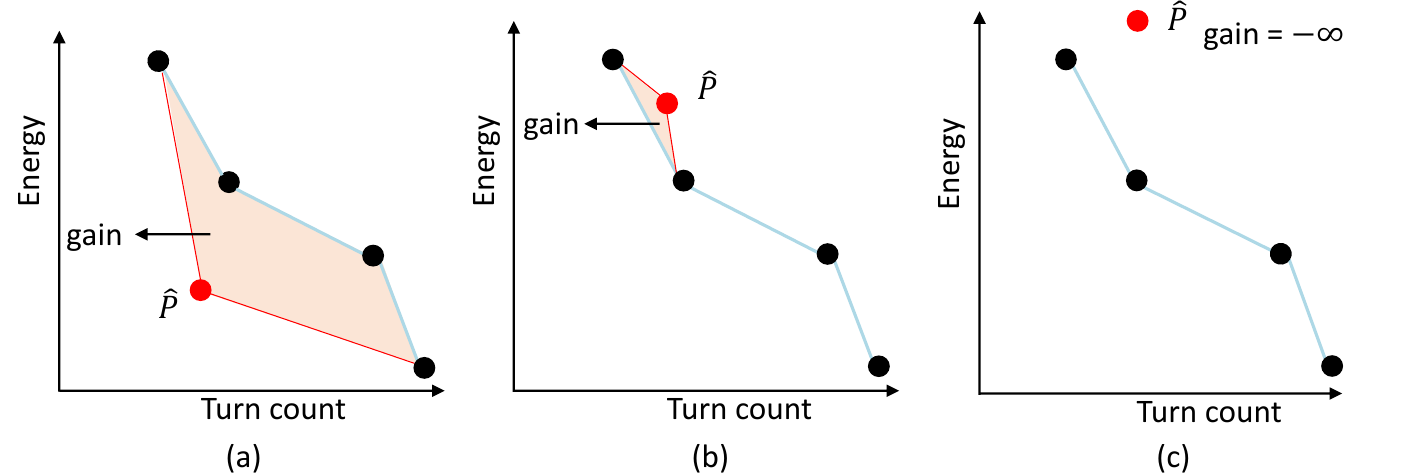}
    \caption{Illustration of the gain function}
    \label{fig:gainFigure}
\end{figure}

Our proposed method relies on a gain function that measures the improvements from replacing an existing edge(s) in a tour with a new edge(s) for the bi-criterion case. Note that replacing existing edge(s) with new edge(s) may or may not result in a connected path but the output can be viewed as a collection of subpaths (with potentially degenerate elements). Suppose this intermediate collection of subpaths is denoted as $\subpathcollection$, and let the associated cost vector, $\objectivevector{\subpathcollection}$, represent the sum of edge attributes of the subpaths. The number of possible $\subpathcollection$s grows exponentially with problem size. Hence, we use the gain function to prioritize edge swaps that may potentially yield better tours. Mathematically, the gain function,  $\gain{\localsearchndset}{\subpathcollection}$, is defined as the area between the original efficiency frontier $\efffrontier{\localsearchndset}$ and the efficiency frontier obtained after adding the point $\objectivevector{\subpathcollection}$. See Figure \ref{fig:gainFigure} for reference. Black points represent the objective values of the non-dominated tours in $\localsearchndset$, and the red point marks the sum of cost attributes of the edges in $\subpathcollection$. The shaded area represents the gain. If the cost attributes of $\subpathcollection$ are \say{non-dominated}, the gain is positive (Figures \ref{fig:gainFigure}a and \ref{fig:gainFigure}b); if it is dominated, the gain is set to $-\infty$ (Figure \ref{fig:gainFigure}c). Both \textsc{S3Opt} and \textsc{S3OptTW} moves involve two main steps with a few key differences.

\begin{figure}[H]
    \centering
    \includegraphics[scale=0.43]{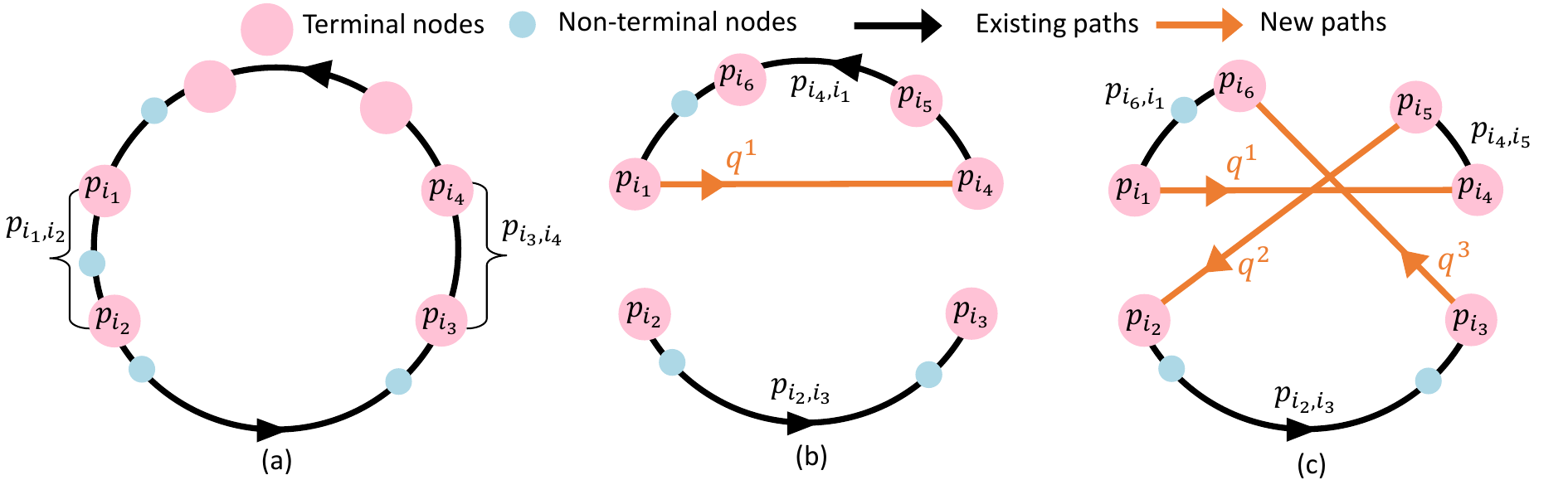}
    \caption{\textsc{S3Opt} and \textsc{S3OptTW} operators illustration}
    \label{img:iopt}
\end{figure}

\noindent\textbf{Step 1:} \textsc{S3Opt} and \textsc{S3OptTW} moves begin by picking a tour $p$ from the sets $\localsearchndset$ and $\explorationset$, respectively (details on the selection procedure will be provided later). We skip these operators if tour $p$ contains less than 6 terminals. Next, we find the set of tuples $\candidates$, referred to as \textit{candidate set}, whose elements are of the form ($\tourindex{p}{i_1},\tourindex{p}{i_2}, \tourindex{p}{i_3}, \tourindex{p}{i_4}, q^1$), and meet the conditions mentioned below. Here, $i_1, i_2, i_3$, and $i_4$ are indices of nodes in $p$. See Figure \ref{img:iopt}b for reference.
\begin{enumerate}[(a)]
    \item $(\tourindex{p}{i_1},\tourindex{p}{i_2})$ and $(\tourindex{p}{i_3}, \tourindex{p}{i_4})$ should be adjacent terminal pairs in tour $p$, i.e., there are no other terminals between $\tourindex{p}{i_1}$ and $\tourindex{p}{i_2}$, and between $\tourindex{p}{i_3}$ and $\tourindex{p}{i_4}$ (non-terminal nodes are allowed). The adjacency criteria ensure that terminals are not skipped when generating a new tour. The path $q^1 \in \pathbank{\tourindex{p}{i_1}}{\tourindex{p}{i_4}}$, where  $\pathbank{u}{v}$ denotes the set of Pareto-optimal paths from terminal $u$ to $v$.
    \item There should be at least two terminals between $\tourindex{p}{i_4}$ and $\tourindex{p}{i_1}$ (not counting $\tourindex{p}{i_4}$ and $\tourindex{p}{i_1}$) in the subpath $\subpath{p}{i_4}{i_1}$, where $\subpath{p}{i}{j}$ denotes the subpath connecting nodes $\tourindex{p}{i}$ to $\tourindex{p}{j}$ in path $p$ (including endpoints).    
    \item For the \textsc{S3OptTW} operator, in addition to the above two conditions, $(\tourindex{p}{i_1}, \tourindex{p}{i_2})$ is selected such that the time-window penalty at terminal $\tourindex{p}{i_2}$ is non-zero, i.e., $\penaltyfunction{p,\tourindex{p}{i_2}}\neq0$. 
\end{enumerate}     

The candidate set $\candidates$ cardinality is generally large for \textsc{S3Opt} and can slow the local search procedure. To address this challenge, an additional filtering step is applied. Let $\meanenergy{\localsearchndset}$ denote the mean energy consumption across the edges in tours belonging to set $\localsearchndset$. If the energy consumption in the segment $\subpath{p}{i_1}{i_2}$, i.e., $\energy{\subpath{p}{i_1}{i_2}}$ is less than $ \length{\subpath{p}{i_1}{i_2}}\meanenergy{\localsearchndset}$, where $\length{\subpath{p}{i}{j}}$ is the number of edges in the path $\subpath{p}{i}{j}$, then all tuples ($\tourindex{p}{i_1},\tourindex{p}{i_2}, \tourindex{p}{i_3}, \tourindex{p}{i_4}, q^1$) that meet this criterion are removed. For the candidate tuples retained after the pruning step, the energy consumption on path $\subpath{p}{i_1}{i_2}$ is higher than the average, and hence replacing $\subpath{p}{i_1}{i_2}$ with $q^1$ will likely improve the tour. We apply the above filtration step on $\subpath{p}{i_1}{i_2}$ and $\subpath{p}{i_3}{i_4}$ using both the cost attributes, i.e., energy and number of turns. 
For each remaining candidate, $(\tourindex{p}{i_1},\tourindex{p}{i_2}, \tourindex{p}{i_3}, \tourindex{p}{i_4}, q^1) \in D$, we add the segment $q^1$ and remove subpaths $\subpath{p}{i_1}{i_2}$ and $\subpath{p}{i_3}{i_4}$ from tour $p$. This results in a cycle $q^1 \concatenate \subpath{p}{i_4}{i_1}$ and a subpath $\subpath{p}{i_2}{i_3}$ as shown in Figure \ref{img:iopt}b. Given two paths, $p$ and $q$, such that the last node of $p$ is the same as the first node of $q$,
the expression $p\concatenate q$ denotes the path formed by concatenating the elements of $q$, starting from its second node to path $p$. For example, given $p=[1,2]$, $q=[2,3]$, we have $p\concatenate q=[1,2,3]$. Let $\subpathcollection$ denote this collection of edges. Next, we evaluate the gain (i.e., $\gain{\localsearchndset}{\subpathcollection}$) for each such $\subpathcollection$, and the candidate with the best gain is passed to Step 2. 

For \textsc{S3OptTW}, the additional condition (c) keeps the candidates set $\candidates$ computationally tractable; hence, filtration w.r.t energy/turns and gain is not required. All elements of set $\candidates$ are passed to Step 2.

\noindent\textbf{Step 2:} Let $\adjterminalvector{p}$ denote the set of index tuples $(i,j)$ such that $\tourindex{p}{i}$ and $\tourindex{p}{j}$ are adjacent terminals in path $p$. For instance, in $p^6=[\textbf{0}, 7, \textbf{5}, 6, \textbf{0}, 1, 2, \textbf{3}, 8, 7, \textbf{5}, 6, \textbf{0}]$ from Table \ref{tab:extour}, $\adjterminalvector{p^6} = [(1, 3), (3, 5), (5, 8), (8, 11), (11, 13)]$ (node indices start from 1 onward). For each pair $(i_5$, $i_6)\in\adjterminalvector{\subpath{p}{i_4}{i_1}}$ (condition (b) in Step 1 ensures the existence of at least one such pair in subpath $\subpath{p}{i_4}{i_1}$), we first initialize $\pathbank{\tourindex{p}{i_5}}{\tourindex{p}{i_2}}$ and $\pathbank{\tourindex{p}{i_3}}{\tourindex{p}{i_6}}$ with respective bi-criteria Pareto-optimal paths, if empty. Next, we connect the terminals $\tourindex{p}{i_5}$ to $\tourindex{p}{i_2}$ and $\tourindex{p}{i_3}$ to $\tourindex{p}{i_6}$ using a randomly selected path $q^2$ from $\pathbank{\tourindex{p}{i_5}}{\tourindex{p}{i_2}}$ and $q^3$ from $\pathbank{\tourindex{p}{i_3}}{\tourindex{p}{i_6}}$. This results in a new tour $q \equiv \subpath{p}{i_2}{i_3} \concatenate q^3 \concatenate \subpath{p}{i_6}{i_1} \concatenate q^1 \concatenate \subpath{p}{i_4}{i_5} \concatenate q^2$ as shown in Figure \ref{img:iopt}c. Iterating across all cycles from Step 1 generates a set of new feasible tours. Recall that in \textsc{S3Opt}, only the candidate with the best gain is passed to Step 2. However, in \textsc{S3OptTW}, all elements of $D$ from Step 1 are passed to Step 2. Thus, the number of new tours can explode. To keep the computations tractable, the new tours are filtered based on Pareto-optimality.

\begin{figure}[H]
    \centering
    \begin{subfigure}[b]{0.27\textwidth}
        \centering
        \includegraphics[width=\textwidth]{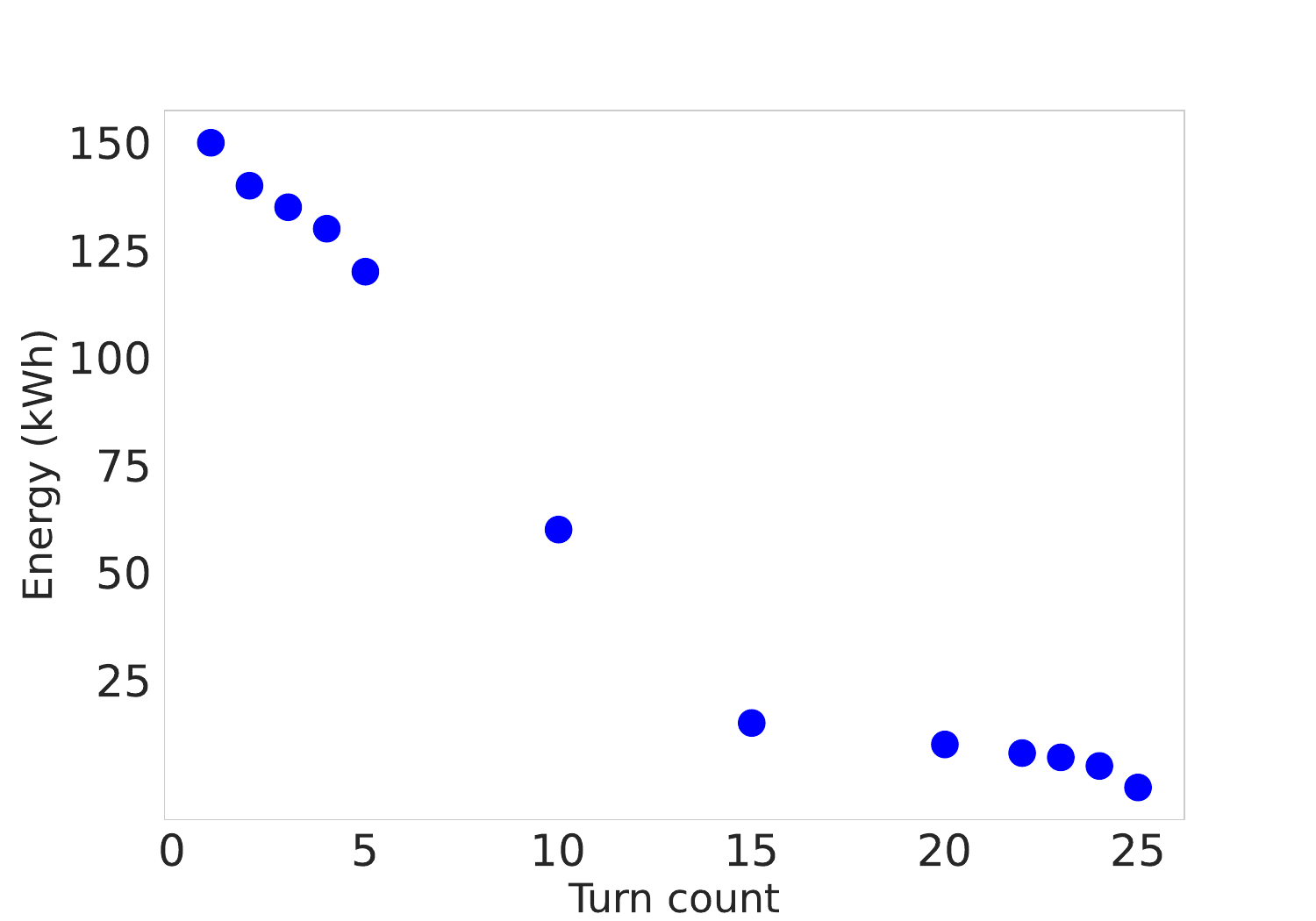}
        \caption{Current $\localsearchndset$}
        \label{fig:subfig-a}
    \end{subfigure}
    \begin{subfigure}[b]{0.27\textwidth}
        \centering
        \includegraphics[width=\textwidth]{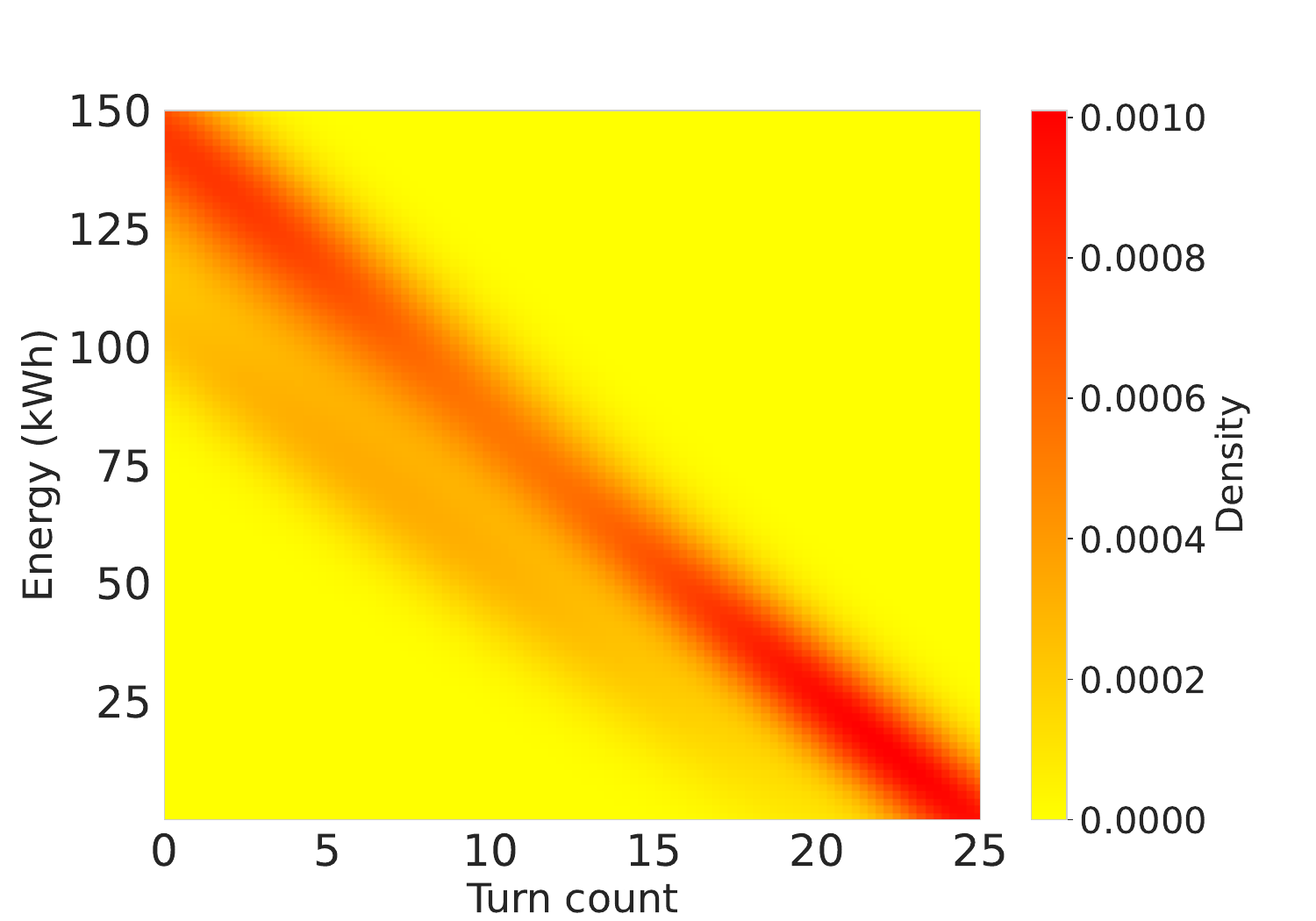}
        \caption{Fitted density function}
        \label{fig:subfig-b}
    \end{subfigure}
    \begin{subfigure}[b]{0.27\textwidth}
        \centering
        \includegraphics[width=\textwidth]{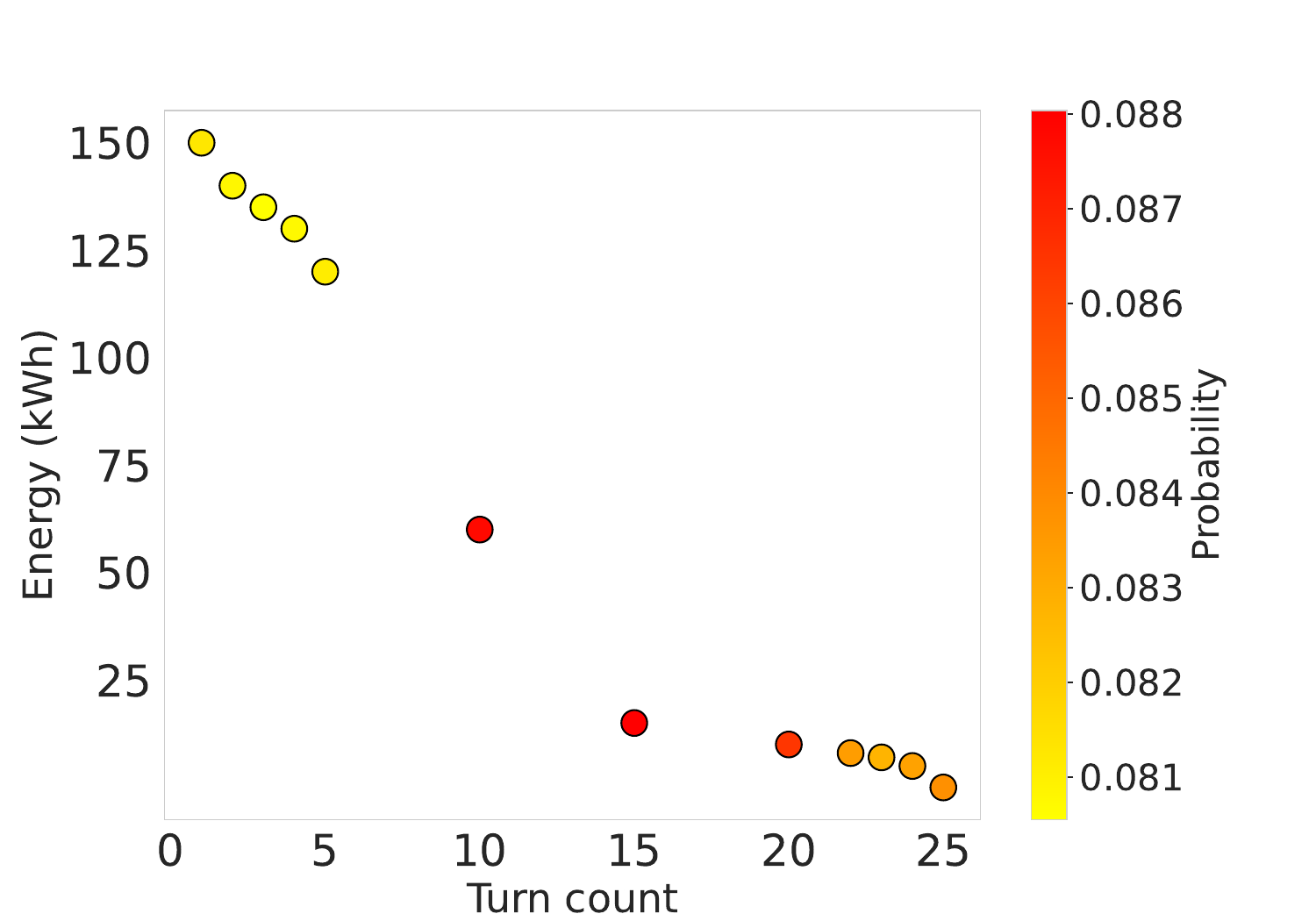}
        \caption{Probability of selection}
        \label{fig:subfig-c}
    \end{subfigure}
    \caption{Tour selection in the \textsc{S3Opt} operator}
    \label{fig:kdillustration}
\end{figure}

In Step 1, instead of randomly selecting a tour in \textsc{S3Opt}, we diversify the search by choosing the tour from a relatively less explored portion of the objective space. To do this, we first estimate the Gaussian kernel density function (\texttt{GKD}) based on $\efffrontier{\localsearchndset}$. Next, for a tour $p \in \localsearchndset$, we define a relative Gaussian kernel density function for tour \texttt{rGKD}($p$) as $\texttt{GKD}\big(\objectivevector{p}\big) / \sum_{q\in\localsearchndset}\texttt{GKD}\big(\objectivevector{q}\big)$. Subsequently, for each tour $p \in \localsearchndset$, we determine the probability of selection using  \eqref{Eq:kde}. The idea is to assign higher probabilities to low-density points by subtracting the densities from 1, thereby prioritizing less explored regions within the objective space. Figure \ref{fig:kdillustration}a shows the efficiency frontier, while Figure \ref{fig:kdillustration}b showcases the corresponding Gaussian kernel density function. 
Figure \ref{fig:kdillustration}c shows the probability of selection of each tour. For \textsc{S3OptTW}, tour $p$ is selected randomly from $\explorationset$. 

\begin{equation}
    \label{Eq:kde}
\texttt{pmf}(p)  = \frac{1-\texttt{rGKD}(p)}{\sum_{q\in \localsearchndset} \big(1-\texttt{rGKD}(q)\big)}
\end{equation}

Algorithm \ref{algo:S3Opt} summarizes the \textsc{S3Opt} operator. Lines 2--3 select a tour based on probabilities from  \eqref{Eq:kde}. Next, Lines 4--8 initialize the candidate set $\candidates$. Lines 9 and 10 find the element with the maximum gain in set $\candidates$. Finally, Lines 11--14 create new tours by connecting the resulting cycle and subpath from Step 1.

\begin{algorithm}[H]
\caption{\textsc{S3Opt}($\localsearchndset$)}\label{algo:S3Opt}
\begin{algorithmic}[1]
\State $\texttt{tour\_set}\gets \emptyset$
\State $\texttt{pmf}(p)\gets $ \eqref{Eq:kde} $\forall p \in \localsearchndset$
\State Select a tour $p\in\localsearchndset$ based on the probabilities $\texttt{pmf}(p)$
\State $\candidates \gets \emptyset$ 
\For{$(\tourindex{p}{i_1},\tourindex{p}{i_2}, \tourindex{p}{i_3}, \tourindex{p}{i_4}, q^1$) in $p$ that satisfy conditions (a) and (b) in Step 1} 
\If{$\energy{\subpath{p}{i_1}{i_2}} \geq  \length{\subpath{p}{i_1}{i_2}}\meanenergy{\localsearchndset}  \textbf{ and } \turn{\subpath{p}{i_1}{i_2}} \geq  \length{\subpath{p}{i_1}{i_2}}\meanturns{\localsearchndset}$}
\If{$\energy{\subpath{p}{i_3}{i_4}} \geq   \length{\subpath{p}{i_3}{i_4}}\meanenergy{\localsearchndset} \textbf{ and } \turn{\subpath{p}{i_3}{i_4}} \geq   \length{\subpath{p}{i_3}{i_4}}\meanturns{\localsearchndset}$ }
\State $\candidates \gets \candidates \cup \{(\tourindex{p}{i_1},\tourindex{p}{i_2}, \tourindex{p}{i_3}, \tourindex{p}{i_4}, q^1)\}$
\EndIf
\EndIf
\EndFor
\State $P^\prime \gets \{(q^1 \concatenate \subpath{p}{i_4} {i_1} , \subpath{p}{i_2}{i_3}): (\tourindex{p}{i_1},\tourindex{p}{i_2}, \tourindex{p}{i_3}, \tourindex{p}{i_4}, q^1)\in \candidates\} $ \Comment{Add segment $q^1$ and remove $\subpath{p}{i_1}{i_2}$, $\subpath{p}{i_3}{i_4}$ from $p$}
\State $(q^1 \concatenate \subpath{p}{i_4} {i_1} , \subpath{p}{i_2}{i_3}) \gets  \argmax_{\subpathcollection\in P^\prime}\gain{\localsearchndset}{\subpathcollection}$
\For {$(i_5, i_6) \in \adjterminalvector{ \subpath{p}{i_4}{i_1}}$} 
\State $q^3\gets$ Select a random path from $\pathbank{\tourindex{p}{i_3}}{\tourindex{p}{i_6}}$
\State $q^2 \gets $ Select a random path from $\pathbank{\tourindex{p}{i_5}}{\tourindex{p}{i_2}}$
\State $\texttt{tour\_set}\leftarrow\texttt{tour\_set}\cup \{\subpath{p}{i_2}{i_3} \concatenate q^3 \concatenate \subpath{p}{i_6}{i_1} \concatenate q^1 \concatenate \subpath{p}{i_4}{i_5} \concatenate q^2\}$ 
\EndFor
\\\Return \texttt{tour\_set}
\end{algorithmic}
\end{algorithm}

Algorithm \ref{algo:S3OptTW} describes the \textsc{S3OptTW} operator. It starts by randomly selecting a tour from set $\explorationset$. Lines 4--6 initialize the candidate set $\candidates$ with the additional condition that the terminal $\tourindex{p}{i_2}$ violates its time window. Lines 7--11 create new tours using Step 2. Finally, Lines 12 and 13 filters the resulting tours to only contain Pareto-optimal subsets using a function $\textsc{ParetoUpdate}$. $\textsc{ParetoUpdate}$ simply adds new non-dominated paths to $\texttt{tour\_set}$. Additionally, any existing tours dominated by the newly added tour are removed, ensuring that $\texttt{tour\_set}$ is always a non-dominated set. 

\begin{algorithm}[H]
\caption{\textsc{S3OptTW}($\explorationset$)}\label{algo:S3OptTW}
\begin{algorithmic}[1]
\State $\texttt{new\_tours}, \texttt{tour\_set}\gets \emptyset, \emptyset$
\State Select a random tour $p$ from $\explorationset$ 
\State $\candidates \gets \emptyset$ 
\For{$(\tourindex{p}{i_1},\tourindex{p}{i_2}, \tourindex{p}{i_3}, \tourindex{p}{i_4}, q^1$) in $p$ that satisfy conditions (a), (b), and (c) in Step 1} 
\State $\candidates \gets \candidates \cup \{(\tourindex{p}{i_1},\tourindex{p}{i_2}, \tourindex{p}{i_3}, \tourindex{p}{i_4}, q^1)\}$
\EndFor
\For{$(\tourindex{p}{i_1},\tourindex{p}{i_2}, \tourindex{p}{i_3}, \tourindex{p}{i_4}, q^1) \in \candidates$} 
\For {$(i_5, i_6) \in \adjterminalvector{ \subpath{p}{i_4}{i_1}}$} 
\State $q^3\gets$ Select a random path from  $\pathbank{\tourindex{p}{i_3}}{\tourindex{p}{i_6}}$
\State $q^2 \gets $ Select a random path from $\pathbank{\tourindex{p}{i_5}}{\tourindex{p}{i_2}}$
\State $\texttt{new\_tours}\leftarrow\texttt{new\_tours}\cup \{\subpath{p}{i_2}{i_3} \concatenate q^3 \concatenate \subpath{p}{i_6}{i_1} \concatenate q^1 \concatenate \subpath{p}{i_4}{i_5} \concatenate q^2\}$ 
\EndFor
\EndFor
\For{tour $p\in \texttt{new\_tours}$}
\State $\texttt{tour\_set}\gets \textsc{ParetoUpdate}(p, \texttt{tour\_set})$ 
\EndFor
\\\Return \texttt{tour\_set}
\end{algorithmic}
\end{algorithm}

\subsubsection{\textsc{RepairTW}}\label{gr}
The main aim of the \textsc{RepairTW} heuristic is to tackle the time-window infeasibility of tours in $\explorationset$. In a nutshell, this operator removes terminals that have a positive time-window penalty and attempts to find a position between two terminals where a smaller subset of the removed terminals can be reinserted. It operates in two steps: Destroy and Repair.

\begin{figure}[H]
    \centering
    \includegraphics[scale=0.4]{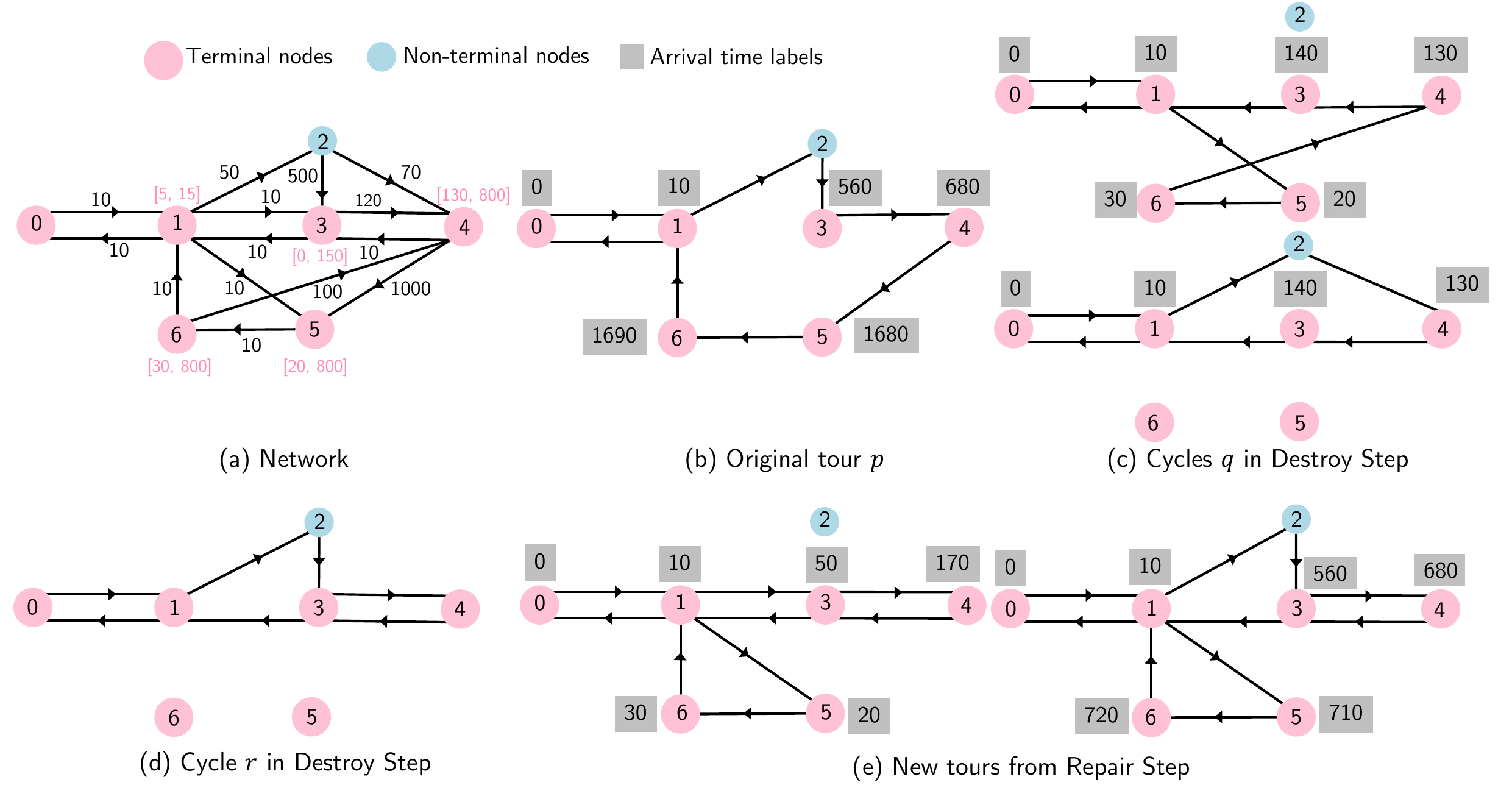}
    \caption{\textsc{RepairTW} operator illustration}
    \label{fig:GRexample}
\end{figure}

\textbf{Step 1 (Destroy):} For a randomly selected tour $p\in \explorationset$, let $\violatedterminalset{p} \subseteq \dualterminal \setminus \{\depot\}$ denote the subset of terminals with strictly positive penalty, i.e., $\penaltyfunction{p,v}> 0$, where $v \in \dualterminal$. If $\violatedterminalset{p} = \dualterminal \setminus \{\depot\}$, i.e., all terminals violate their time windows in tour $p$, then we skip the \textsc{RepairTW} operator. Otherwise, we delete the subpaths connecting every occurrence of terminal $v \in \violatedterminalset{p}$ in tour $p$ to its adjacent terminals. A sequence of time-window feasible terminals is then obtained. For illustration, consider the network shown in Figure \ref{fig:GRexample}a. The edge weights represent travel times. Figure \ref{fig:GRexample}b shows a tour $p = [\textbf{0}, \textbf{1}, 2, \textbf{3}, \textbf{4}, \textbf{5}, \textbf{6}, \textbf{1}, \textbf{0}]$, where $\adjterminalvector{p} = [(1, 2), (2, 4), (4, 5), (5, 6), (6, 7), (7, 8), (8, 9)]$. The grey node labels denote the arrival time at the first visit of the terminal. Thus, we have  $\violatedterminalset{p} = \{3, 5,6 \}$. After deleting the sub-paths connected to terminals in $\violatedterminalset{p}$, we are left with a sequence $[\textbf{0}, \textbf{1}, \textbf{4}, \textbf{1}, \textbf{0}]$, in which all terminals satisfy their time windows. 

The removed portions between terminals with feasible time windows are then reconnected using the shortest-time paths, resulting in a cycle $q$ (which may or may not be a tour since a tour is a cycle that passes through all the terminals). If the shortest paths between the terminals happen to pass through the remaining terminals, and if the resulting tour satisfies $\penaltyfunction{q}=0$, the procedure stops. In the earlier example, the shortest paths from $1$ to $4$ are $[\textbf{1},\textbf{5},\textbf{6},\textbf{4}]$ and $[\textbf{1},2,\textbf{4}]$, and from $4$ to $1$ is $[\textbf{4},\textbf{3},\textbf{1}]$. Two cases arise as seen in  Figure \ref{fig:GRexample}c: (1) $q=[\textbf{0}, \textbf{1}, \textbf{5}, \textbf{6}, \textbf{4}, \textbf{3}, \textbf{1}, \textbf{0}]$, in which case the algorithm stops since $q$ covers all the terminals and satisfies $\penaltyfunction{q}=0$. (2) $q=[\textbf{0}, \textbf{1}, 2, \textbf{4}, \textbf{3}, \textbf{1}, \textbf{0}]$, which is not feasible since terminals $5$ and $6$ are not part of the cycle. 

Given a cycle $q$ which does not contain all terminals, we find a subset of terminals $\feasibleinsertionset{q}{i}{j} \subseteq \violatedterminalset{q}~\forall~(i, j)\in\adjterminalvector{q}$ that can be inserted between $i^{th}$ and $j^{th}$ position in cycle $q$, one at a time using shortest time paths, while keeping the cycle $q$ time-window feasible. We then find a pair $(i^*, j^*)$ that allows us to insert the maximum number of terminals, i.e., $(i^*, j^*) \in \arg\max_{(i,j)\in \adjterminalvector{q}}|\feasibleinsertionset{q}{i}{j}|$ (ties are broken randomly). Continuing with the previous example, we have $\adjterminalvector{q} = [(1, 2), (2, 4), (4, 5), (5, 6), (6, 7)]$, $(i^*, j^*)=(2, 4)$, and $\feasibleinsertionset{q}{i^*}{j^*}=\{5,6\}$.

Next, we go back to the original tour $p$ and remove the terminals in $\feasibleinsertionset{q}{i^*}{j^*}$ and the subpaths connected to them. The repair step will attempt to insert only these terminals back to create a new tour. We do this because the number of elements in $\feasibleinsertionset{q}{i^*}{j^*}$ is less than that in $\violatedterminalset{q}$, increasing the chances of satisfying time windows after they are reinserted. As before, we connect the deleted portions using the shortest time paths, and if the resulting cycle (say $r$) passes through all terminals and is time-window feasible, the procedure terminates. If not, we proceed to the repair step. 

\textbf{Step 2 (Repair):} Let $\violatedterminalset{r}$ denote the set of terminals in $\dualterminal$ that are missing from cycle $r$. Removing terminals $5$ and $6$ from tour $p$ yields cycle $r = [\textbf{0}, \textbf{1}, 2, \textbf{3}, \textbf{4}, \textbf{3}, \textbf{1}, \textbf{0}]$ as illustrated in Figure \ref{fig:GRexample}d. In this step, we build new tours that insert these missing terminals in cycle $r$. To this end, we start by finding the shortest time path, $p^1$, that starts from terminal $\tourindex{q}{i^*}$ and visits all the terminals in $\violatedterminalset{r}$ in a greedy nearest-neighbor manner. In the earlier example, $|\feasibleinsertionset{q}{i^*}{j^*}| = \{5, 6\}$ and $\tourindex{q}{i^*} = 1$. Thus, $p^1=[1,5,6]$. Next, we insert path $p^1$ in cycle $r$. Specifically, for all pairs $(i, j)\in\adjterminalvector{r}$ such that $\tourindex{r}{i}=\tourindex{q}{i^*}$, we build new tours by replacing the subpath $\subpath{r}{i}{j}$ with $p^1\concatenate p^2$, where $p^2$ is the shortest path from the last node in $p^1$ to the terminal $\tourindex{r}{j}$.

In the earlier example, $  \adjterminalvector{r} = [(1,2), (2,4), (4,5), (5,6), (6,7), (7,8)]$ and $\tourindex{r}{2}=\tourindex{r}{7}=\tourindex{q}{i^*}=1$, and there are two possible places to insert the missing terminals -- between $(i,j)= (2,4)$ and $(i,j)= (7,8)$. In the first case, we have $\subpath{r}{2}{4} = [1,2,3]$ and  $p^2= [6,1,3]$. The resulting new tour from replacing the subpath $\subpath{r}{2}{4}$ with $p^1\concatenate p^2 = [1,5, 6,1,3]$, is $[\textbf{0},\textbf{1},\textbf{5},\textbf{6},\textbf{1},\textbf{3},\textbf{4},\textbf{3},\textbf{1},\textbf{0}]$, as shown in Figure \ref{fig:GRexample}e (left). Similarly, in the second case, we get $\subpath{r}{7}{8} = [1,0]$, $p^2= [6,1,0]$,  $p^1\concatenate p^2 = [1,5, 6,1,0]$, and the new tour is $[\textbf{0},\textbf{1},2,\textbf{3},\textbf{4},\textbf{3},\textbf{1},\textbf{5},\textbf{6},\textbf{1},\textbf{0}]$, as shown in Figure \ref{fig:GRexample}e (right). The tour on the left in Figure \ref{fig:GRexample}e adheres to the time windows, but the tour in the right panel of Figure \ref{fig:GRexample}e violates the time window at terminal 3. Yet, the repair process reduces the number of terminals violating time windows from three to one. Additional tours can be generated by using precomputed Pareto-optimal paths instead of exclusively relying on the shortest paths,

The pseudocode for \textsc{RepairTW} operator is presented in Algorithm \ref{algo:gf}. Lines 1--4 randomly select a tour $p\in\explorationset$ and initialize cycle $q$ by removing the terminals in set  $\violatedterminalset{q}$ from tour $p$. Lines 7--10 create cycle $r$. For each revisit of terminal $\tourindex{q}{i^*}$ in tour $p$, Lines 13--16 create a new tour by inserting the missing terminals in cycle $r$. Note that since $|r|$ represents the number of edges in the path $r$, $|r|+1$ denotes the number of nodes in $r$.

\begin{algorithm}
\caption{\textsc{RepairTW}($\explorationset$)}\label{algo:gf}
\begin{algorithmic}[1]
\State $\texttt{tour\_set}\gets \emptyset$ 
\State Select a random tour $p$ from $\explorationset$
\State $\violatedterminalset{p}= \{v:v \in \dualterminal, \penaltyfunction{p,v} \neq 0\}$ 
\State Create cycle $q$ by removing terminals $v\in \violatedterminalset{p}$  \Comment{Step 1 (Destroy)}
\If{cycle $q$ has all terminals and $\penaltyfunction{q} = 0$}
\State \textbf{return} $\{q\}$
\EndIf
\For {$(i, j) \in\adjterminalvector{q}$}
\State $\feasibleinsertionset{q}{i}{j}\gets\{k:p_k\in \violatedterminalset{p}, p_k$ can be inserted between $\tourindex{q}{i}$ and $\tourindex{q}{j}$ s.t. $\penaltyfunction{q, k}=0$ and $\penaltyfunction{q, k}=0 \}$
\EndFor
\State $(i^*, j^*) \gets \arg\max_{(i,j)\in \adjterminalvector{q}}|\feasibleinsertionset{q}{i}{j}|$
\State  $r \gets$  remove terminals in $\feasibleinsertionset{q}{i^*}{j^*}$  from $p$ and fill the missing portions using shortest-time paths 
\If{cycle $r$ has all terminals and $\penaltyfunction{r} = 0$}
\State \textbf{return} $\{r\}$
\EndIf
\State $ \violatedterminalset{r}\gets  \{v:v\in \dualterminal, v\not\in r\}$
\For{$(i,j)\in\adjterminalvector{r}$} 
\If{$\tourindex{r}{i}=\tourindex{q}{i^*}$}
\State Initialize $p^1$ and $p^2$ \Comment{Step 2 (Repair)} 
\State  $\texttt{new\_tour}\gets\subpath{r}{i}{j}\concatenate p^1 \concatenate p^2 \concatenate \subpath{r}{j}{|r|+1}$   
\State $\texttt{tour\_set}\gets \texttt{tour\_set} \cup \{\texttt{new\_tour}\}$ 
\EndIf
\EndFor
\\\Return \texttt{tour\_set}

\end{algorithmic}
\end{algorithm}

\subsubsection{\textsc{FixedPerm}}
\label{subsub:fpm}
The \textsc{RepairTW} operator changes the permutation of the terminals and repairs time-window infeasibilities. However, it makes minimal effort in optimizing the path between terminals. We handle such inter-terminal path optimization using the \textsc{FixedPerm} operator. Given a tour $p$, we find new tours such that the order of terminals visited is maintained. A naive way to do this is to improve the paths between all the adjacent terminals using a bi-criteria shortest path algorithm. Since multiple Paerto-optimal paths exist between each pair of terminals, we get a set of tours as output. However, computing bi-objective shortest paths for all adjacent terminals is time-consuming. In addition, combining these new paths can result in an exponential number of tours. Thus, we must efficiently select a subset of adjacent terminal pairs. To this end, we assign priorities to all adjacent terminal pairs based on the difference between the energy consumption (number of turns) of the path connecting the terminals and the mean energy (number of turns) consumption of the edges in the graph. 

Specifically, given a tour $p$, we first randomly select an objective, say energy. Next, for all adjacent terminal pairs $(\tourindex{p}{i}, \tourindex{p}{j})$, where $(i, j)\in \adjterminalvector{p}$, we calculate the energy consumption of the path $\subpath{p}{i}{j}$, i.e., $\energy{\subpath{p}{i}{j}}$. 
For each $(i, j)\in \adjterminalvector{p}$, we assign a priority, $\temrinalPriority{i,j}$ using \eqref{eq:prob}, which help identify segments with cost attributes that are higher than the average. Thus, replacing segments having higher priority with other subpaths will likely improve the tour.  
\begin{equation}
\temrinalPriority{i, j} = 
\dfrac{\max\left\{0, \energy{\subpath{p}{i}{j}} - \length{\subpath{p}{i}{j}}\meanenergy{\dualedge}\right\}}{\sum_{(i', j') \in \adjterminalvector{p}} \max\{0,\ \energy{\subpath{p}{i'}{j'}} - \length{\subpath{p}{i'}{j'}} \meanenergy{\dualedge}\}}  \label{eq:prob}
\end{equation}
For $(i, j) \in \adjterminalvector{p}$, let $\pathbanksubset{\tourindex{p}{i}}{\tourindex{p}{j}}$ denote a set of paths from terminal $\tourindex{p}{i}$ to $\tourindex{p}{j}$ (not necessarily Pareto-optimal). $\pathbanksubset{\tourindex{p}{i}}{\tourindex{p}{j}}$ is  initialized with the existing subpath connecting terminals $\tourindex{p}{i}$ and $\tourindex{p}{j}$ (i.e., $\subpath{p}{i}{j}$) in the current tour $p$. Let $\pathbank{\tourindex{p}{i}}{\tourindex{p}{j}}$ be the set of all Pareto-optimal paths between $\tourindex{p}{i}$ and $\tourindex{p}{j}$. We scan the top $\hyperParameter{1}$ elements of $\adjterminalvector{p}$ based on the $\temrinalPriority{i, j}$ values, where $\hyperParameter{1}$ is a hyperparameter and add paths from $\pathbank{\tourindex{p}{i}}{\tourindex{p}{j}}$ to $\pathbanksubset{\tourindex{p}{i}}{\tourindex{p}{j}}$. The algorithms populates new tours by concatenating paths in $\pathbanksubset{\tourindex{p}{i}}{\tourindex{p}{j}}~\forall~(i, j) \in \adjterminalvector{p}$. We save and reuse the path set $\pathbank{\tourindex{p}{i}}{\tourindex{p}{j}}$ to avoid solving additional shortest path problems in other operators and future iterations.

The initial selection of a tour $p$ for \textsc{FixedPerm} is made from a set of paths that are \textit{skewed}. A tour $p \in \localsearchndset$ is assumed to be \textit{skewed} if either of its cost attributes ($\energy{p}$ or $\turn{p}$) are $\hyperParameter{2}$ standard deviations away from their respective means over the edges of the network, i.e., if either $(\energy{p} - \length{p}\meanenergy{\dualedge}) \notin [- \hyperParameter{2} \stdenergy, \hyperParameter{2}  \stdenergy]$ or 
$(\turn{p} - \length{p} \meanturns{\dualedge})\notin [- \hyperParameter{2} \stdturns, \hyperParameter{2}  \stdturns]$, where $\hyperParameter{2}$ is another hyperparameter. This step increases the chances of investigating unexplored regions of the decision space.

Algorithm \ref{algo:FixedPermNbd} illustrates the pseudocode for \textsc{FixedPerm}. Lines 1--5 select a random tour $p$ from the set of skewed tours, and initialize the variables $\temrinalPriority{i,j}$ and $\pathbanksubset{\tourindex{p}{i}}{\tourindex{p}{j}}$ for $(i,j)\in \adjterminalvector{p}$. The for-loop in Lines 6 and 7 populates first $\hyperParameter{1}$  elements of $\pathbanksubset{\tourindex{p}{i}}{\tourindex{p}{j}}$ with Pareto-optimal paths. Finally, Line 8 generates a set of new tours. The concatenation operator applied to a set of vectors concatenates every possible combination, similar to the cartesian product.

\begin{algorithm}[H]
\caption{\textsc{FixedPerm}($\localsearchndset$ or $\setY, \hyperParameter{1}, \hyperParameter{2} $)}\label{algo:FixedPermNbd}
\begin{algorithmic}[1]
\State $\texttt{skewed\_tours} \gets \emptyset$
\State Add tour $p \in \localsearchndset$ to $\texttt{skewed\_tours}$ if $(\energy{p} - \length{p}\meanenergy{\dualedge})\notin [- \hyperParameter{2} \stdenergy, \hyperParameter{2}  \stdenergy]$ and 
$(\turn{p} - \length{p} \meanturns{\dualedge})\notin [- \hyperParameter{2} \stdturns, \hyperParameter{2}  \stdturns]$
\State Select a random tour $p$ from \texttt{skewed\_tours} 
\State Update $\temrinalPriority{i,j}$ using \eqref{eq:prob} $\forall \ (i, j) \in \ \adjterminalvector{p}$
\State $\pathbanksubset{\tourindex{p}{i}}{\tourindex{p}{j}}  \gets \{\subpath{p}{i}{j}\}\ \forall\  (i, j) \in \adjterminalvector{p} $
\For{top $\hyperParameter{1}$ elements $(i,j)\in \adjterminalvector{p}$ based on $\temrinalPriority{i,j}$}
\State $\pathbanksubset{\tourindex{p}{i}}{\tourindex{p}{j}} \gets \pathbanksubset{\tourindex{p}{i}}{\tourindex{p}{j}}  \cup \pathbank{\tourindex{p}{i}}{\tourindex{p}{j}}$ 
\EndFor
\State $\texttt{tour\_set} \gets \concatenate_{(i,j) \in\ \adjterminalvector{p}}  \pathbanksubset{\tourindex{p}{i}}{\tourindex{p}{j}}$
\\\Return \texttt{tour\_set}
\end{algorithmic}
\end{algorithm}

\subsubsection{\textsc{Other Operators}}
\textbf{Quad}:
The \textsc{Quad} operator is a 4-opt-like move that introduces randomness in the local search process by generating multiple tours and diversifying the search. For a randomly selected tour $p\in\localsearchndset$, it selects four adjacent terminal pairs $(\tourindex{p}{i_8}, \tourindex{p}{i_1}), (\tourindex{p}{i_2}, \tourindex{p}{i_3}), (\tourindex{p}{i_4}, \tourindex{p}{i_5}),$ and $ (\tourindex{p}{i_6}, \tourindex{p}{i_7})$ and replaces the paths between them to create new feasible tours for BSTSP (i.e., tours that may violate time windows). The selected pairs should not have terminals in common. We complete the tour by joining terminals $\tourindex{p}{i_8}$ to $\tourindex{p}{i_5}$, $\tourindex{p}{i_6}$ to $\tourindex{p}{i_3}$, $\tourindex{p}{i_4}$ to $\tourindex{p}{i_1}$, and $\tourindex{p}{i_2}$ to $\tourindex{p}{i_7}$, which creates a double-bridge connection that cannot be reached using sequential 3-opt moves \citep{kanellakis1980local}. To join these terminal pairs, we use paths in $\pathbank{\tourindex{p}{i_8}}{\tourindex{p}{i_5}}$, $\pathbank{\tourindex{p}{i_6}}{\tourindex{p}{i_3}}$, $\pathbank{\tourindex{p}{i_4}}{\tourindex{p}{i_1}}$, and $\pathbank{\tourindex{p}{i_2}}{\tourindex{p}{i_7}}$, respectively.

\textbf{RandPermute}:
\textsc{RandPermute} involves creating a random permutation of terminals and connecting them using precomputed bi-objective shortest paths. This operator aims to introduce diversification in the search process.

\subsection{Local Search Procedure}\label{ls:full}
This section describes the overall local search procedure. The set $\explorationset$ is maintained to track only the top $\hyperParameter{3}$ BSTSP tours with the least time-windows penalties. This helps limit our attention to promising tours that are likely to satisfy time window restrictions after a few local search operations. To account for the differences in runtimes of the operators, we introduce three additional hyperparameters: $\hyperParameter{4}$, $\hyperParameter{5}$, and $\hyperParameter{6}$, representing the number of times \textsc{FixedPerm}, \textsc{Quad}, and \textsc{RepairTW} are repeated in each iteration, respectively. 

Algorithm \ref{algo:ls} presents the pseudocode for the local search procedure. Lines 1--2 initialize the six hyperparameters $(\hyperParameter{1}, \hyperParameter{2},\hyperParameter{3}, \hyperParameter{4}, \hyperParameter{5}, \hyperParameter{6})$, and the set to store Pareto-optimal paths $\pathbank{u}{v}$ between terminal pairs $u$ and $v$. Line 3 initializes the sets $\explorationset$, $\setY$, and $\localsearchndset$ as outlined in Section \ref{ls:init}. Lines 4--25 contain the main \textit{while} loop. Each iteration starts by updating the mean energy consumption and turns corresponding to edges belonging to tour in $\localsearchndset$. Following this we apply the \textsc{S3Opt} operator on set $\localsearchndset$ and update the sets $\explorationset$, $\setY$, and $\localsearchndset$ using the \textsc{UpdateSets} function (described later). Next, we apply \textsc{FixedPerm} $\hyperParameter{4}$ times. Since both \textsc{S3Opt} and \textsc{FixedPerm} aim at intensification, we then apply the diversification moves \textsc{RandPermute} and \textsc{Quad}. Afterward, we apply \textsc{S3OptTW} and \textsc{RepairTW} to improve time-window infeasibilities. Finally, we call \textsc{FixedPerm} on set $\setY$. The proposed sequence aims first to achieve intensification, followed by applying diversification moves. This enriches the three sets $\explorationset, \setY$, and $\localsearchndset$ and the set of shortest paths $\pathbank{u}{v}$, which saves computational effort in the \textsc{S3OptTW} and \textsc{RepairTW} moves. Applying these operators earlier within the iteration was less effective since their performance depends on the tours in set $\explorationset$. Initially, $\explorationset$ contains tours with very high time-window penalties, making them less likely to become time-window feasible. However, as the iteration progresses, the penalties of tours in $\explorationset$ decrease due to the max penalty clause in \textsc{UpdateSets}, described next. 

\begin{algorithm}[H]
\caption{\textsc{LocalSearch}}\label{algo:ls}
\begin{algorithmic}[1]
\State Initialize global variables $\hyperParameter{1}$, $\hyperParameter{2}$, $\hyperParameter{3}$, $\hyperParameter{4}$, $\hyperParameter{5}$, $\hyperParameter{6}$ \Comment{Hyperparameters}
\State Initialize sets $\explorationset, \setY, \localsearchndset$ \Comment{Refer Section \ref{ls:init}} 
\While{\textit{maximum time not exceeded}}  
\State Update $ \meanenergy{\localsearchndset}$ and $ \meanturns{\localsearchndset}$
\State $\texttt{tour\_set} \gets \textsc{S3Opt}(\localsearchndset)$
\State $\explorationset, \setY, \localsearchndset  \gets \textsc{UpdateSets}(\texttt{tour\_set}, \explorationset, \setY, \localsearchndset $)

\For{$\hyperParameter{4}$ times}
\State $\texttt{tour\_set} \gets \textsc{FixedPerm}(\localsearchndset)$
\State $\explorationset, \setY, \localsearchndset  \gets \textsc{UpdateSets}(\texttt{tour\_set}, \explorationset, \setY, \localsearchndset $)
\EndFor

\State $\texttt{tour\_set} \gets \textsc{RandPermute}()$
\State $\explorationset, \setY, \localsearchndset  \gets \textsc{UpdateSets}(\texttt{tour\_set}, \explorationset, \setY, \localsearchndset $)

\For{$\hyperParameter{5}$ times}
\State $\texttt{tour\_set} \gets \textsc{Quad}(\localsearchndset)$
\State $\explorationset, \setY, \localsearchndset  \gets \textsc{UpdateSets}(\texttt{tour\_set}, \explorationset, \setY, \localsearchndset $)
\EndFor

\State $\texttt{tour\_set} \gets \textsc{S3OptTW}(\explorationset)$
\State $\explorationset, \setY, \localsearchndset  \gets \textsc{UpdateSets}(\texttt{tour\_set}, \explorationset, \setY, \localsearchndset $)

\For{$\hyperParameter{6}$ times}
\State $\texttt{tour\_set} \gets \textsc{RepairTW}(\explorationset)$
\State $\explorationset, \setY, \localsearchndset  \gets \textsc{UpdateSets}(\texttt{tour\_set}, \explorationset, \setY, \localsearchndset $)
\EndFor

\For{$\hyperParameter{4}$ times}
\State $\texttt{tour\_set} \gets \textsc{FixedPerm}(\setY)$
\State $\explorationset, \setY, \localsearchndset  \gets \textsc{UpdateSets}(\texttt{tour\_set}, \explorationset, \setY, \localsearchndset $)
\EndFor


\EndWhile\\
\Return $\localsearchndset$

\vspace{2mm}
\setcounter{ALG@line}{0}
\Procedure{\textsc{UpdateSets}}{$\texttt{tour\_set}, \explorationset, \setY, \localsearchndset $}

\State $\texttt{max\_penalty} \gets \max_{q\in \explorationset}\penaltyfunction{q}$  \Comment{If $\explorationset=\emptyset$, $\texttt{max\_penalty} = \infty$}
\For{tour $p \in \texttt{tour\_set}$} 
\If{$\penaltyfunction{p} = 0$}
    \State $\localsearchndset \gets \textsc{ParetoUpdate}(p, \localsearchndset)$
    \If{$p\not\in \localsearchndset$ }
    \State $\setY \gets \setY \cup \{p\}$
    \EndIf
\ElsIf{$\penaltyfunction{p} < \texttt{max\_penalty}$ \textbf{or} $|\explorationset| < \hyperParameter{3}$}
    \State $\texttt{max\_penalty} \gets \max_{q\in \explorationset}\penaltyfunction{q}$
    \State{$\explorationset \gets \explorationset \cup \{p\}$}
\EndIf
\EndFor
\State Sort the tours $p\in\setY$ based on $\energy{p}$ or $\turn{p}$ (selected randomly)
\State Sort the tours $p\in\explorationset$ based on $\penaltyfunction{p}$ + tour wait time
\State Update sets $\explorationset$ and $\setY$ to keep first $\hyperParameter{3}$ sorted tours
\\\Return $\explorationset, \setY, \localsearchndset$
\EndProcedure
\end{algorithmic}
\end{algorithm}

The function \textsc{UpdateSets} iterates over all new tours. Each time-window feasible tour is either added to $\setY$ or $\localsearchndset$ depending on its Pareto-optimality. We add time-window infeasible tours to $\explorationset$ if they improve the maximum penalty of tours in $\explorationset$ or if the cardinality of $\explorationset$ is less than a threshold. To keep the computations tractable, we limit the size of sets $\explorationset$ and $\setY$. For the set $\setY$, based on experimentation, its size is limited to the top $\hyperParameter{3}$ tours sorted by one of the cost attributes, chosen randomly in each iteration.

The set $\explorationset$ stores only the top $\hyperParameter{3}$ tours based on the sum of tour penalty and \textit{tour wait time}. Wait time for a tour $p$ is calculated as the total wait time across all terminals in $p$. For a terminal $v \in \dualterminal$ that meets its time window requirements (i.e., $\penaltyfunction{p, v}=0$), the wait time is $\max\{0, \patharrivaltime{p}{v,1} - \earliesttime{v}\}$. The wait time is set to zero for terminals that do not meet their time windows. Sorting tours by the sum of penalty and wait time, rather than penalty alone, is beneficial because a tour with a low penalty can have a high wait time. This is because all terminals are assumed to be serviced during their first visit while calculating the penalties and wait times. Moreover, minimizing tour wait time also helps address parking issues that arise from early arrival.

\section{Results}\label{sec:results}
All algorithms discussed in previous sections were implemented in Python 3.9 and run on an Intel(R) Xeon(R) Gold 6154 CPU clocked at 3.00GHz with 512 GB RAM. The MIP model was solved using IBM's CPLEX 22.1.1. Section \ref{ip} compares the solution quality between the MIP and the local search method for small benchmark instances. Section \ref{ls} analyzes the effectiveness of the local search in real-world scenarios. The source codes are available at \href{https://github.com/transnetlab/last-mile-ev-logistics.git}{github.com/transnetlab/last-mile-ev-logistics}.

\subsection{Benchmark: Local Search vs. MIP Formulation}\label{ip}
We benchmarked the local search using the Solomon-Potvin-Bengio dataset \citep{lopez2013travelling}, comprising 30 asymmetric instances for the TSPTW. We selected the top six instances with the highest number of terminals and adapted them to the Steiner case by creating two scenarios where 10\% and 20\% of the nodes were designated as terminals, resulting in 12 problems in total. Terminal nodes retain their original time windows, while those not designated as terminals have their time windows relaxed. Turn costs (0 or 1) were introduced randomly, and energy consumption values were drawn from a uniform distribution between 0 and 10. We set the maximum permissible number of terminal revisits to either the number of terminals or four, whichever is lesser. Due to the small size of these instances, hyperparameters required for the local search procedure were set to high values, ensuring thorough exploration of the solution space. Each test instance was allowed a maximum duration of four hours in MIP and two hours in local search. To keep the MIP tractable, we cap the number of copies of each vertex to $\min\{\sizedualterminal, 4\}$, thus limiting the number of revisits to 4 for instances with more terminals.

\begin{figure}[H]
    \centering
    \begin{subfigure}{0.24\textwidth}
        \centering
        \includegraphics[width=\linewidth]{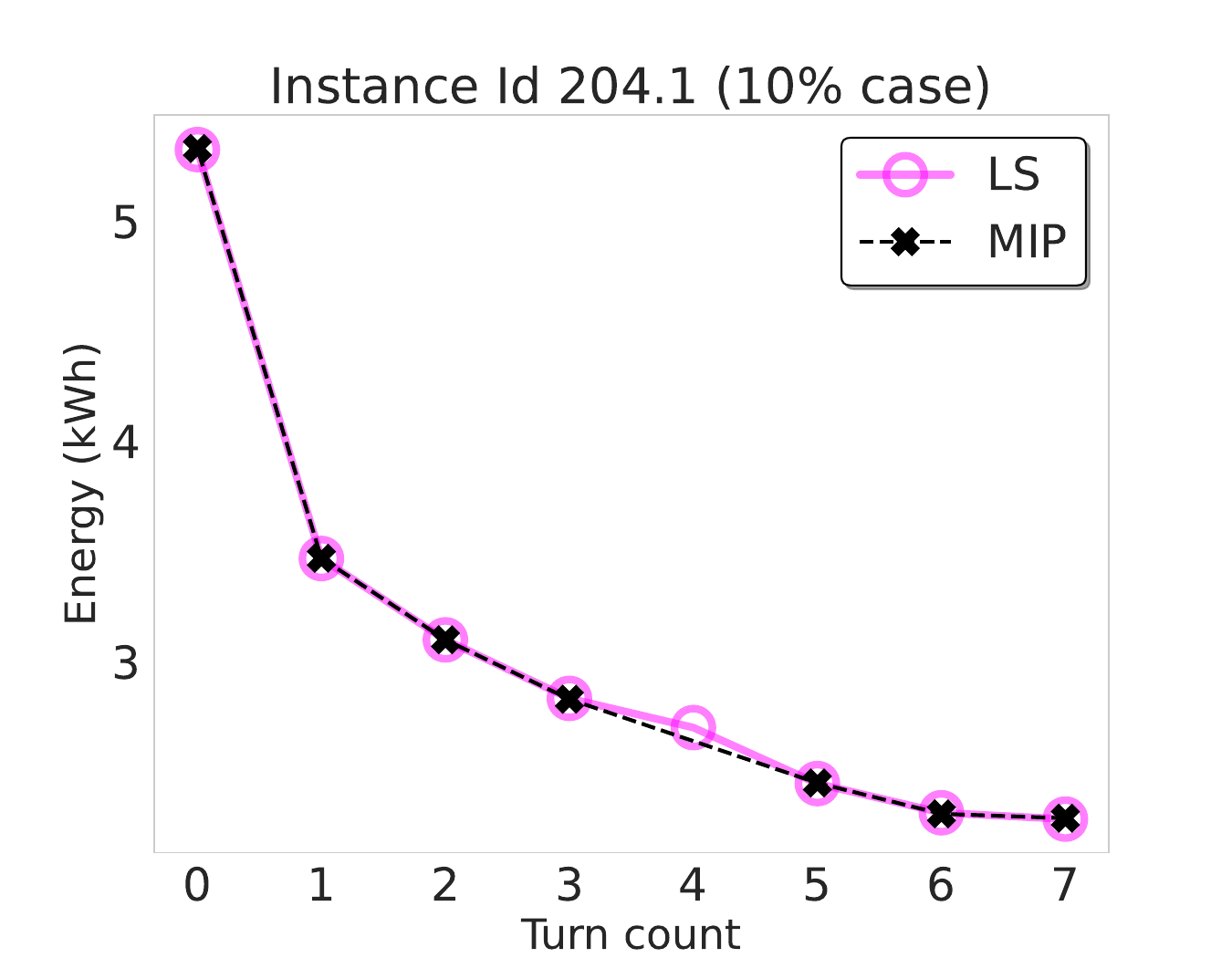}
    \end{subfigure}%
    \begin{subfigure}{0.24\textwidth}
        \centering
            \includegraphics[width=\linewidth]{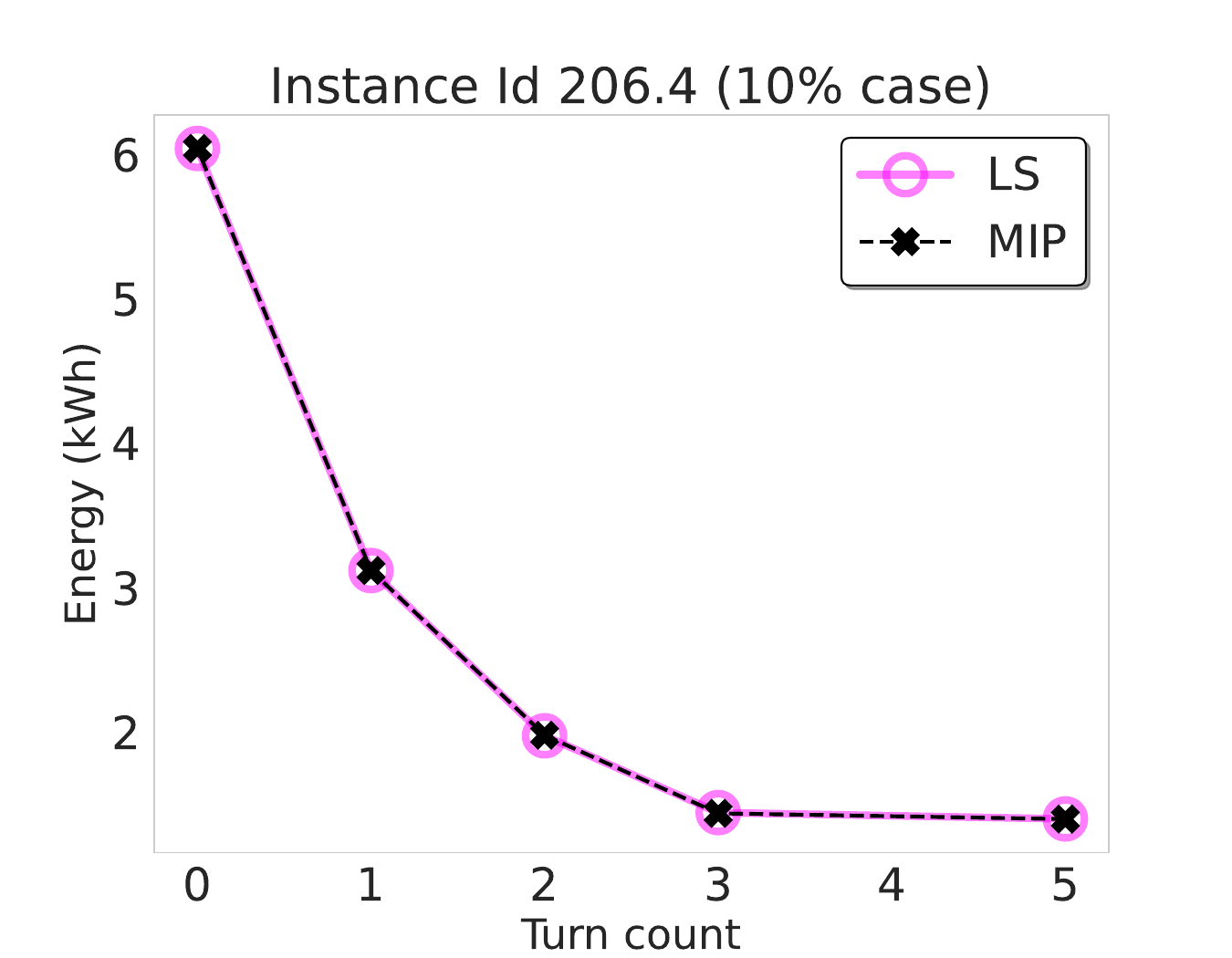}
    \end{subfigure} 
    \begin{subfigure}{0.24\textwidth}
        \centering
        \includegraphics[width=\linewidth]{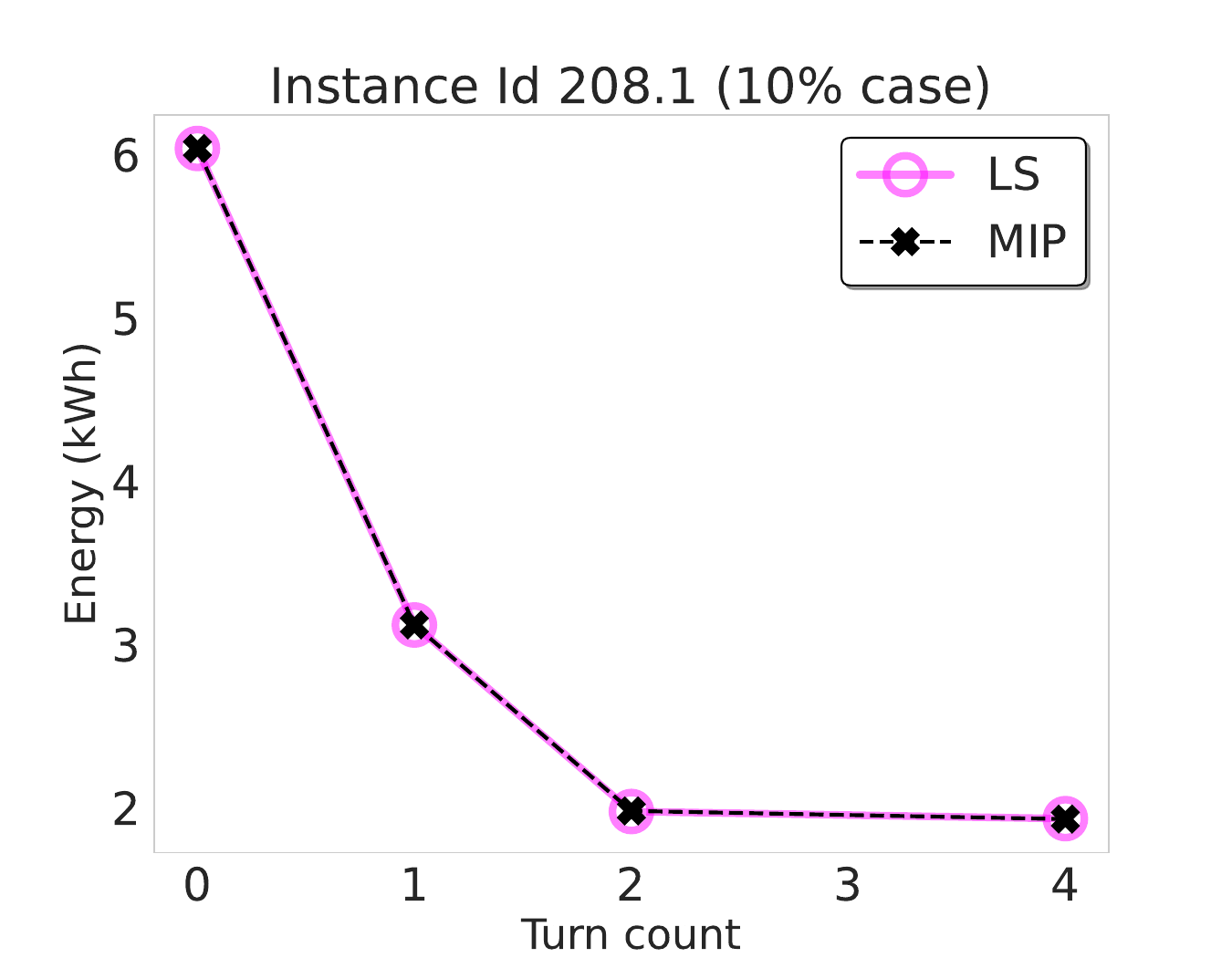}
    \end{subfigure}
     \begin{subfigure}{0.24\textwidth}
        \centering
        \includegraphics[width=\linewidth]{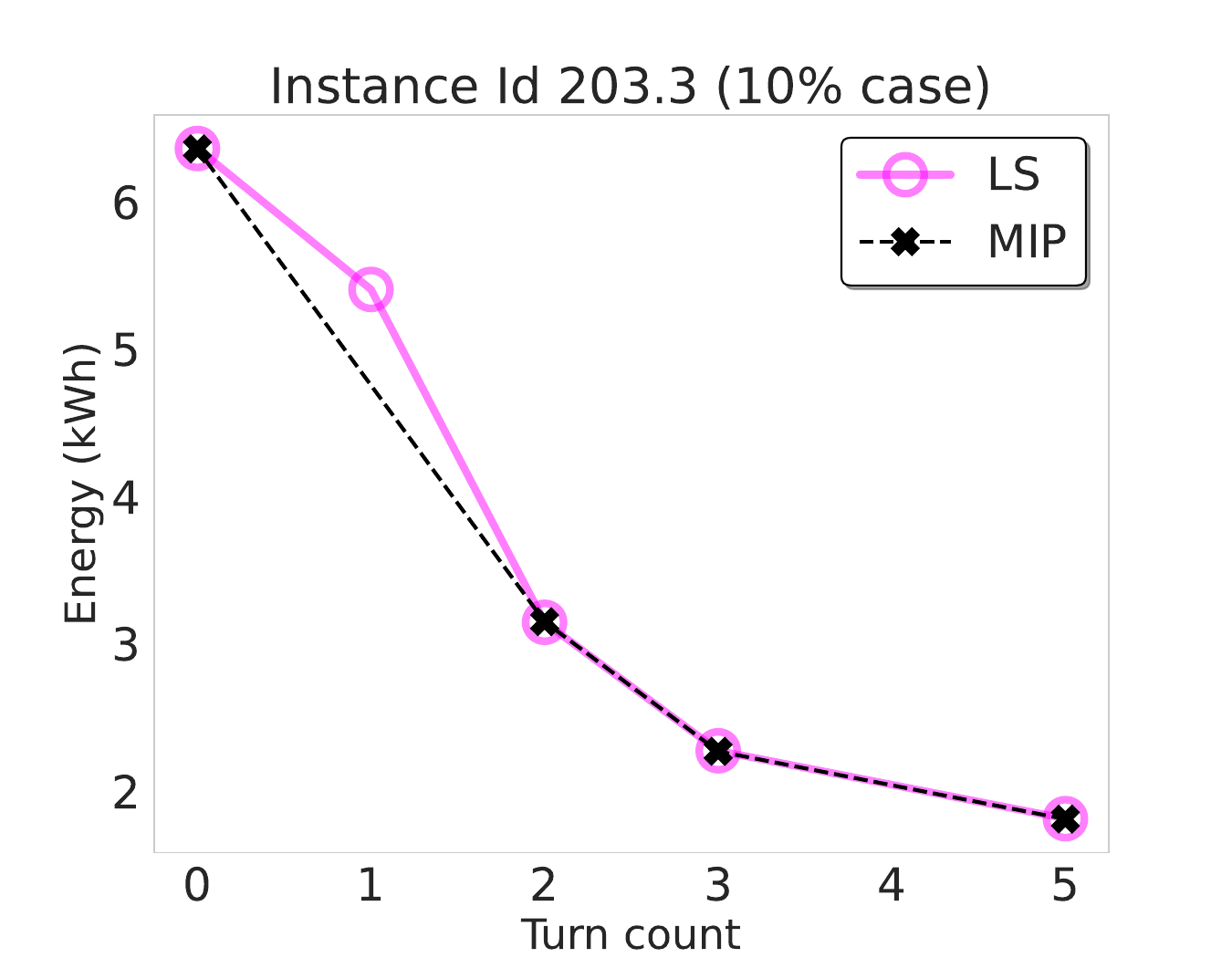}
    \end{subfigure}\\
    \begin{subfigure}{0.24\textwidth}
        \centering
        \includegraphics[width=\linewidth]{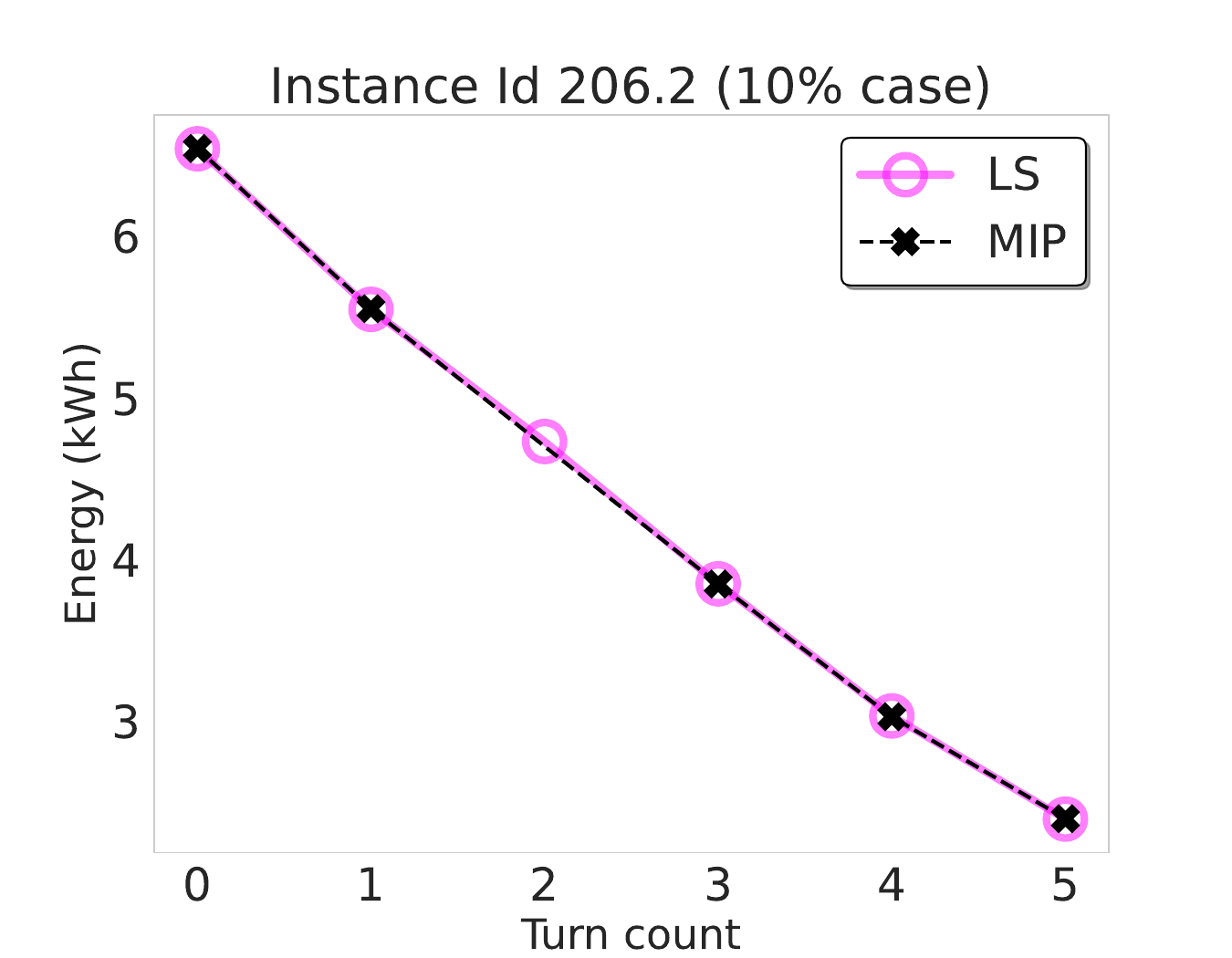}
    \end{subfigure}%
    \begin{subfigure}{0.24\textwidth}
        \centering
        \includegraphics[width=\linewidth]{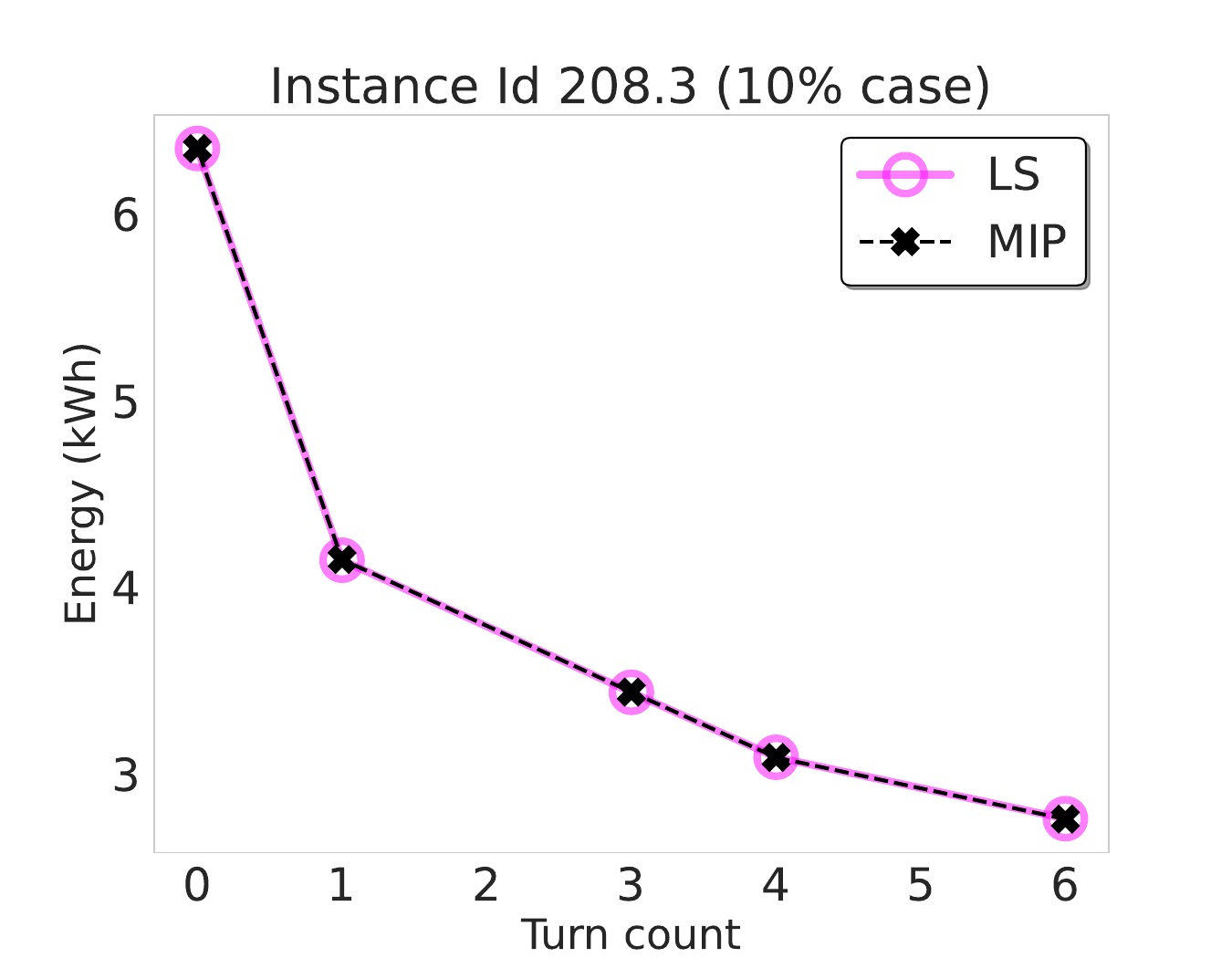}
    \end{subfigure} 
    \begin{subfigure}{0.24\textwidth}
        \centering
        \includegraphics[width=\linewidth]{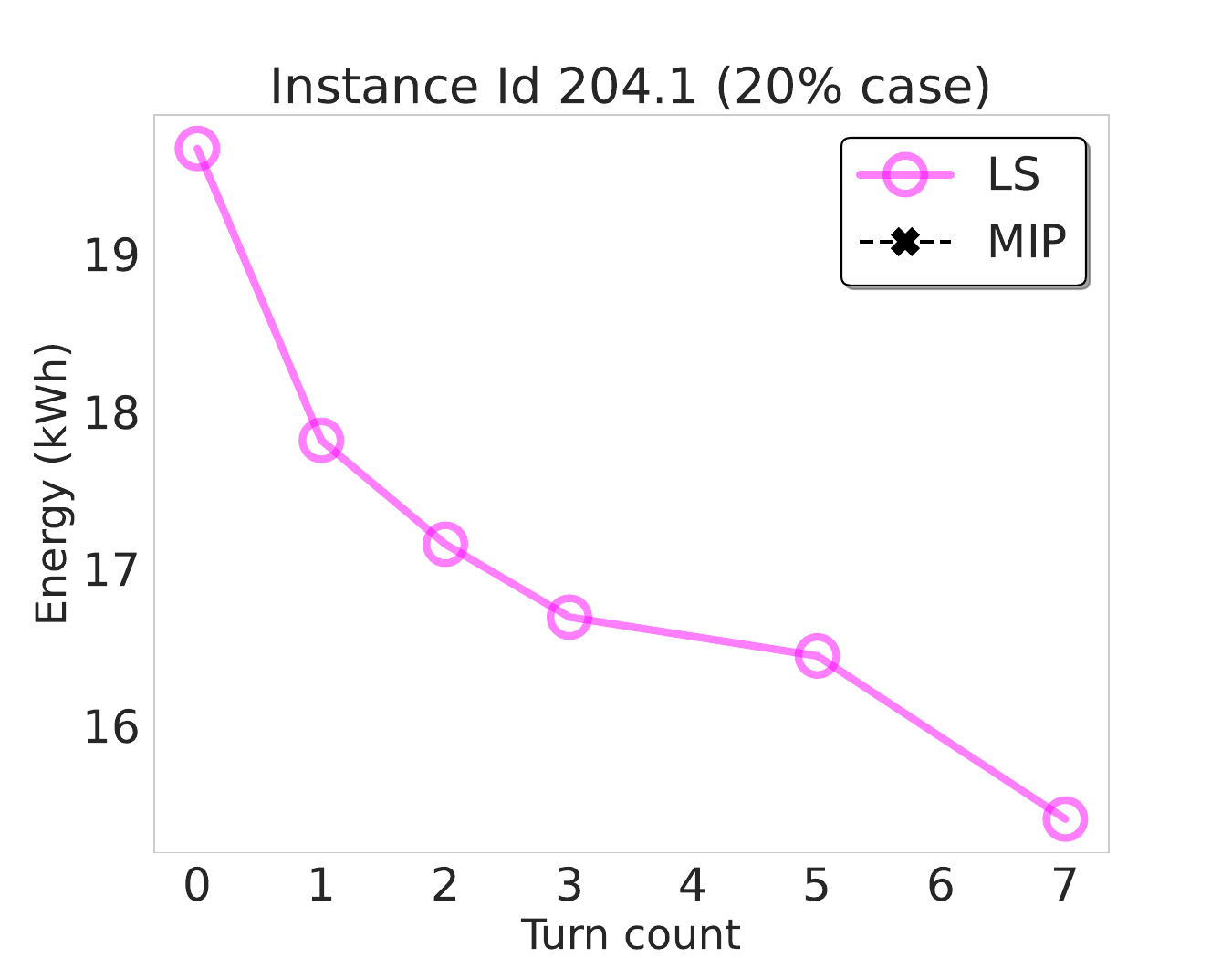}
    \end{subfigure}
     \begin{subfigure}{0.24\textwidth}
        \centering
        \includegraphics[width=\linewidth]{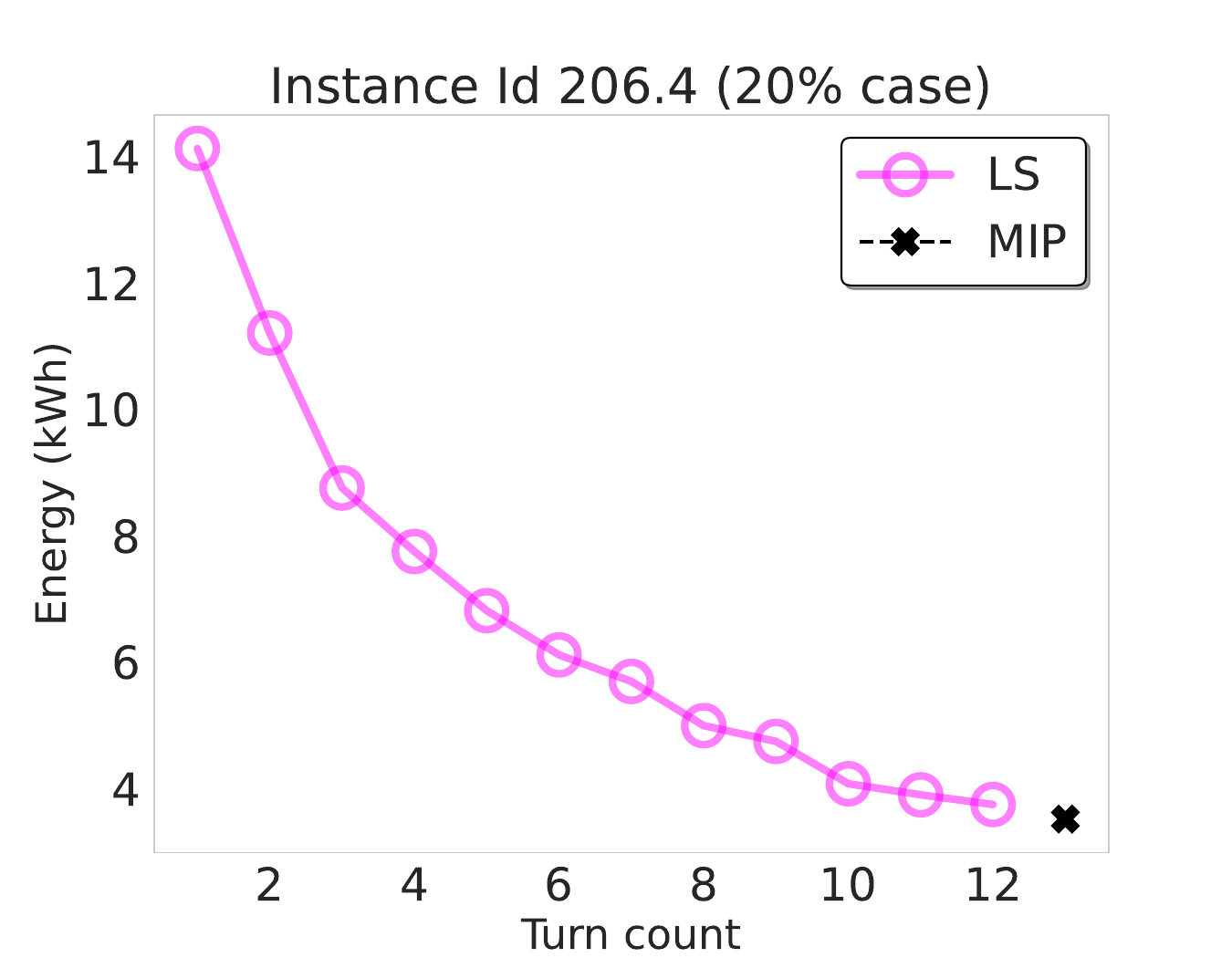}
    \end{subfigure}
    \begin{subfigure}{0.24\textwidth}
        \centering
        \includegraphics[width=\linewidth]{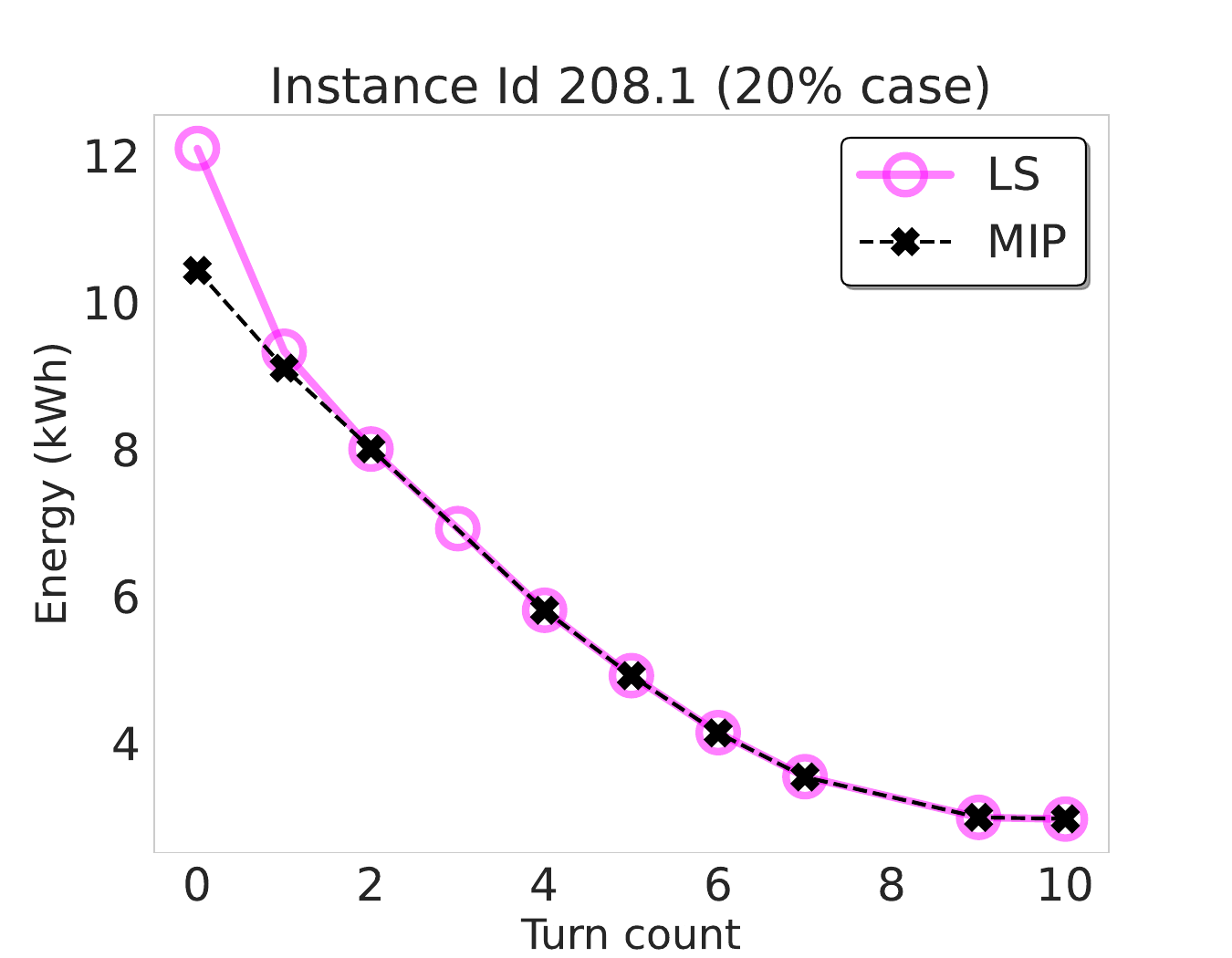}
    \end{subfigure}
        \begin{subfigure}{0.24\textwidth}
        \centering
        \includegraphics[width=\linewidth]{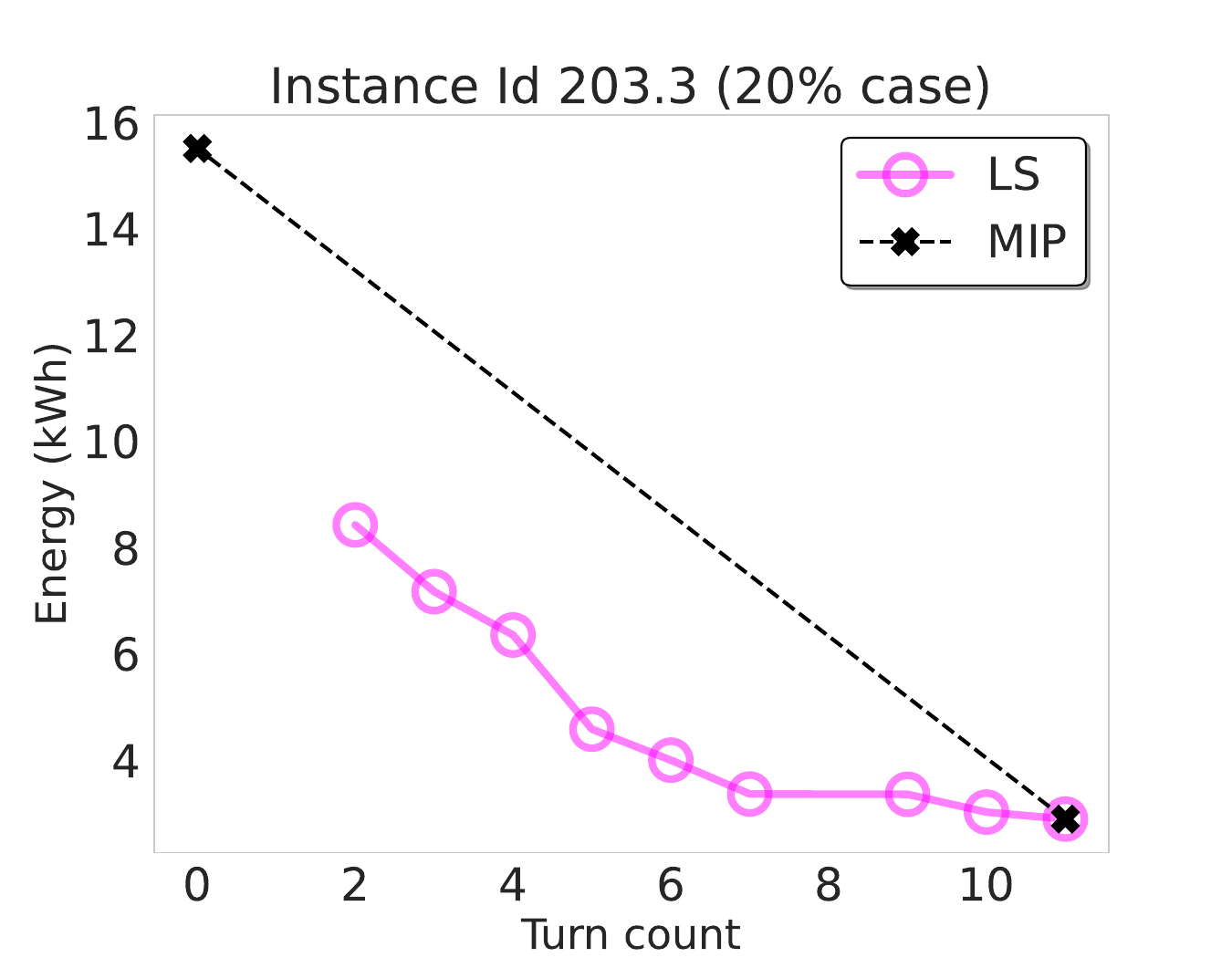}
    \end{subfigure}%
    \begin{subfigure}{0.24\textwidth}
        \centering
        \includegraphics[width=\linewidth]{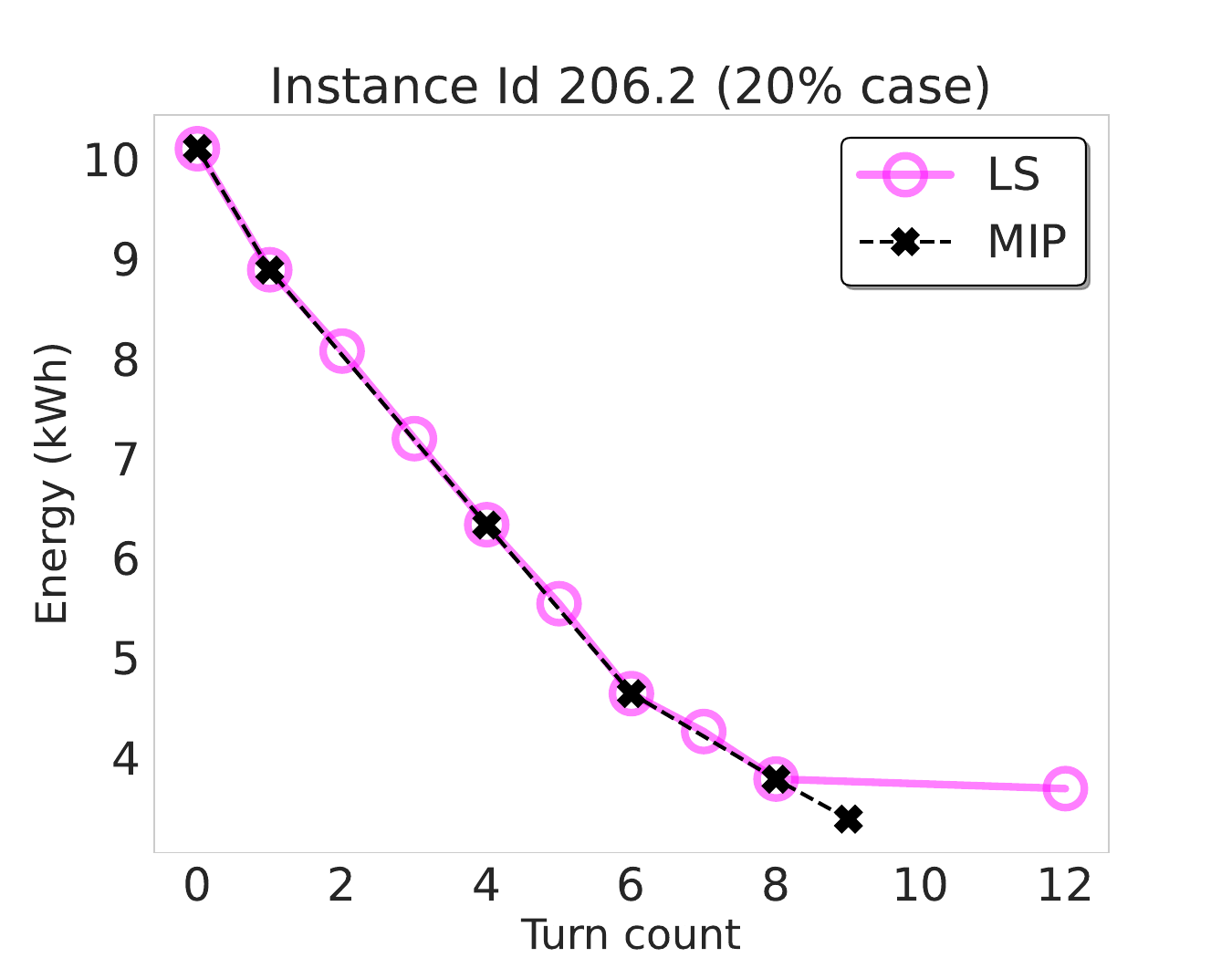}
    \end{subfigure}
        \begin{subfigure}{0.24\textwidth}
        \centering
        \includegraphics[width=\linewidth]{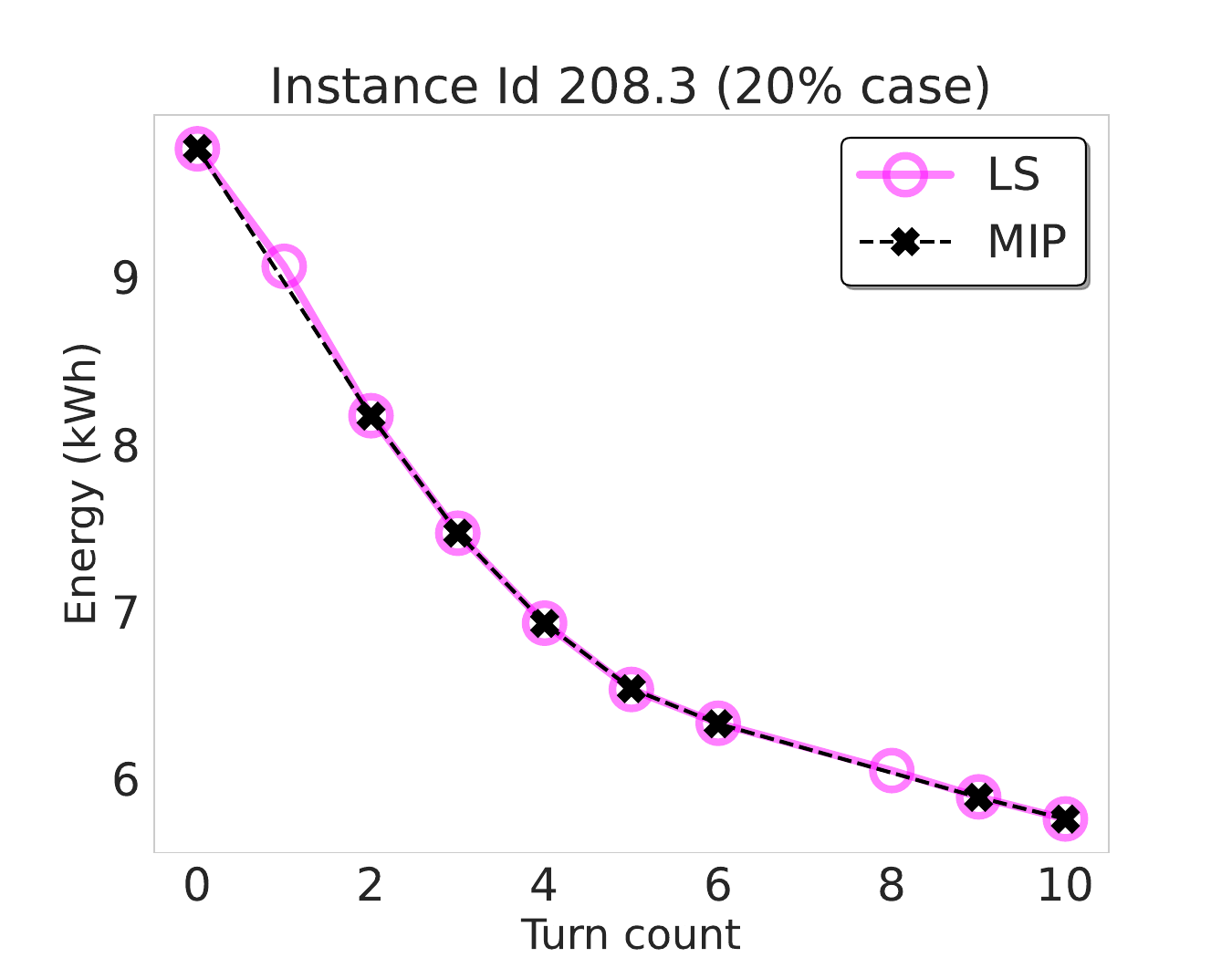}
    \end{subfigure}\\
    \caption{Efficiency frontiers from MIP and local search on benchmark instances}
    \label{fig:30mins}
\end{figure}

\begin{table}[H]
    \centering
    \caption{Benchmark results for local search vs. MIP (\textit{Scenarios}: \% of nodes acting as terminals,  \textit{Instance}: Instance Id, \textit{Nodes}: Number of nodes in the network, \textit{Edges}: Number of edges in the network, \textit{$\sizedualterminal$}: Number of terminals, \textit{$|\localsearchndset_{MIP}|$}: Pareto-optimal BSTSPTW tours count from MIP, \textit{$|\localsearchndset_{LS}|$}: Pareto-optimal BSTSPTW tours count from local search) }
    \label{tab:MIPresults}
    \begin{tabular}{ccccc|c|c}
        \hline
        \multirow{1}{*}{\textbf{Scenarios}} & \multirow{1}{*}{\textbf{Id}} & \multirow{1}{*}{\textbf{Nodes}} & \multirow{1}{*}{\textbf{Edges}} & \multirow{1}{*}{\boldsymbol{$\sizedualterminal$}} & \multicolumn{1}{c|}{\boldsymbol{$|\localsearchndset_{MIP}|$}}& \multicolumn{1}{c}{\boldsymbol{$|\localsearchndset_{LS}|$}}   \\
        \hline
        \multirow{6}{*}{10\%} 
        & 204.1 &  46 & 2070 & 4 & 7   & 8    \\
        & 206.4 &  38 & 1406 & 3 & 5   & 5    \\%
        & 208.1 &  38 & 1406 & 3 & 4   & 4   \\%
        & 203.3 &  37 & 1332 & 3 & 4   & 5   \\%
        & 206.2 &  37 & 1332 & 3 & 5   & 6   \\%
        & 208.3 &  36 & 1260 & 3 & 5   & 5   \\%
        \hline
        \multirow{6}{*}{20\%}
        & 204.1 &  46 & 2070 & 9 & 0   & 6   \\
        & 206.4 &  38 & 1406 & 7 & 1   & 12    \\%
        & 208.1 &  38 & 1406 & 7 & 9   & 10   \\%
        & 203.3 &  37 & 1332 & 6 & 2  & 9   \\%
        & 206.2 &  37 & 1332 & 6 & 6   & 10   \\%
        & 208.3 &  36 & 1260 & 6 & 8   & 10   \\%
        \hline
    \end{tabular}
\end{table}
Table \ref{tab:MIPresults} summarizes the results in the two scenarios. Column \textit{Id} corresponds to the instance name in the original dataset. Columns \textit{Nodes}, \textit{Edges}, and $\sizedualterminal$ denote the number of nodes, edges, and terminals, respectively. Values in column $|\localsearchndset_{MIP}|$ and $|\localsearchndset_{LS}|$ indicate the final number of Pareto-optimal BSTSPTW tours from Algorithm \ref{alg:IP_sol} (Scalarization combined with MIP) and local search, respectively. Figure \ref{fig:30mins} shows the efficiency frontier corresponding to the columns \textit{$|\localsearchndset_{MIP}|$} and \textit{$|\localsearchndset_{LS}|$}.

Results indicate that our local search outperformed MIP in almost all cases. Although MIP is an exact method, local search gave more points in a few cases because MIP can only discover a convex subset of the Pareto-optimal BSTSPTW solutions (refer Section \ref{sec:exact}). For test case 204.1 (20\%), the MIP failed to generate a single tour in four hours. The performance of the MIP degraded rapidly when either the maximum permissible number of terminal revisits or the proportion of nodes designated as terminals is increased.

\subsection{Real-world Application}\label{ls}
To demonstrate the performance of our model in real-world settings, we used the Amazon last-mile routing research challenge dataset \citep{merchan20222021} comprising historical route information of Amazon's last-mile deliveries. Each route is characterized by driver-operated sequence of delivery stops, package dimensions, and delivery time windows. Specifically, we focused on data from Austin, USA, which has 214 routes. Out of these routes, only 87 that meet the following criteria were analyzed: (a) All terminals should be reachable from each other, (b) There should be no negative edge cycle either in the primal or line graph, and (c) At least one terminal should have a time window. The number of total delivery locations ranged from 80 to 209 (an average of 146.1), and the ones with non-trivial time windows were between 1 to 31 (an average of 9.6). We obtained road network attributes such as length and speed limits from \href{https://www.openstreetmap.org/}{www.openstreetmap.org} (OSM). The direction of turns (left, right, or straight) between two road links was determined using a bearing angle threshold of $45^{\circ}$. We did not have accurate information on stop signs and signals and hence a left turn at an intersection was assumed to be conflicting if there were at least two incoming or outgoing arcs. Our analysis revealed that the original time windows in the Amazon dataset were quite lenient, extending up to 8 hours, making nearly all BSTSP tours time-window feasible. To increase the difficulty of the instances, we reduced the upper limit of all the time windows by 60\%.

For estimating energy consumption, we adopt a physics-based model derived from \cite{travesset2015transport} as outlined in  \eqref{eq:energy}. This model quantifies the energy required to traverse a distance $\delta$, considering several parameters including vehicle mass ($\Delta$), fictive mass for rolling inertia ($\Delta_f$), coefficient of rolling resistance ($\nu$), road gradient angle ($\phi$), air density ($\rho$), drag coefficient ($C_x$), equivalent vehicle cross-section ($A$), vehicle speed (\textit{VS}), acceleration ($\eta$), and opposing wind speed (\textit{WS}). 
\begin{equation}
\text{energy} =\left[\Delta g(\nu\cos\phi+\sin\phi)+0.5\left(\rho C_xA(\text{\textit{VS}}+\text{\textit{WS}})^2\right)+(\Delta+\Delta_f)\eta\right]\delta
\label{eq:energy}
\end{equation}
In the present paper, we consider Lion-8 straight electric trucks (\href{https://thelionelectric.com/documents/en/Lion8_all_applications.pdf}{thelionelectric.com}) and use the following parameter values based on literature: $\Delta=27,216$ kg, $g=9.8$ m/s$^2$, $\nu =0.0058, \rho = 1.1$ kg/m$^3$ $, C_x=0.6, A=5.4$m$^2$ $, \text{\textit{WS}} = 0$. For each link, the vehicle speed (\textit{VS}) was set to its OSM speed limit and was assumed to be constant throughout the link ($\eta =0$). Regenerative braking is assumed for links where the energy from  \eqref{eq:energy} is negative, indicating energy gain, with a regenerative efficiency factor of 0.7 \citep{travesset2015transport}. To determine road gradients, we leverage the United States Geological Survey Elevation point query service (\href{https://apps.nationalmap.gov/epqs/}{nationalmap.gov/epqs}).

\begin{figure}[H]
    \centering
    \begin{subfigure}{0.45\textwidth}
        \centering
        \includegraphics[width=.87\linewidth]{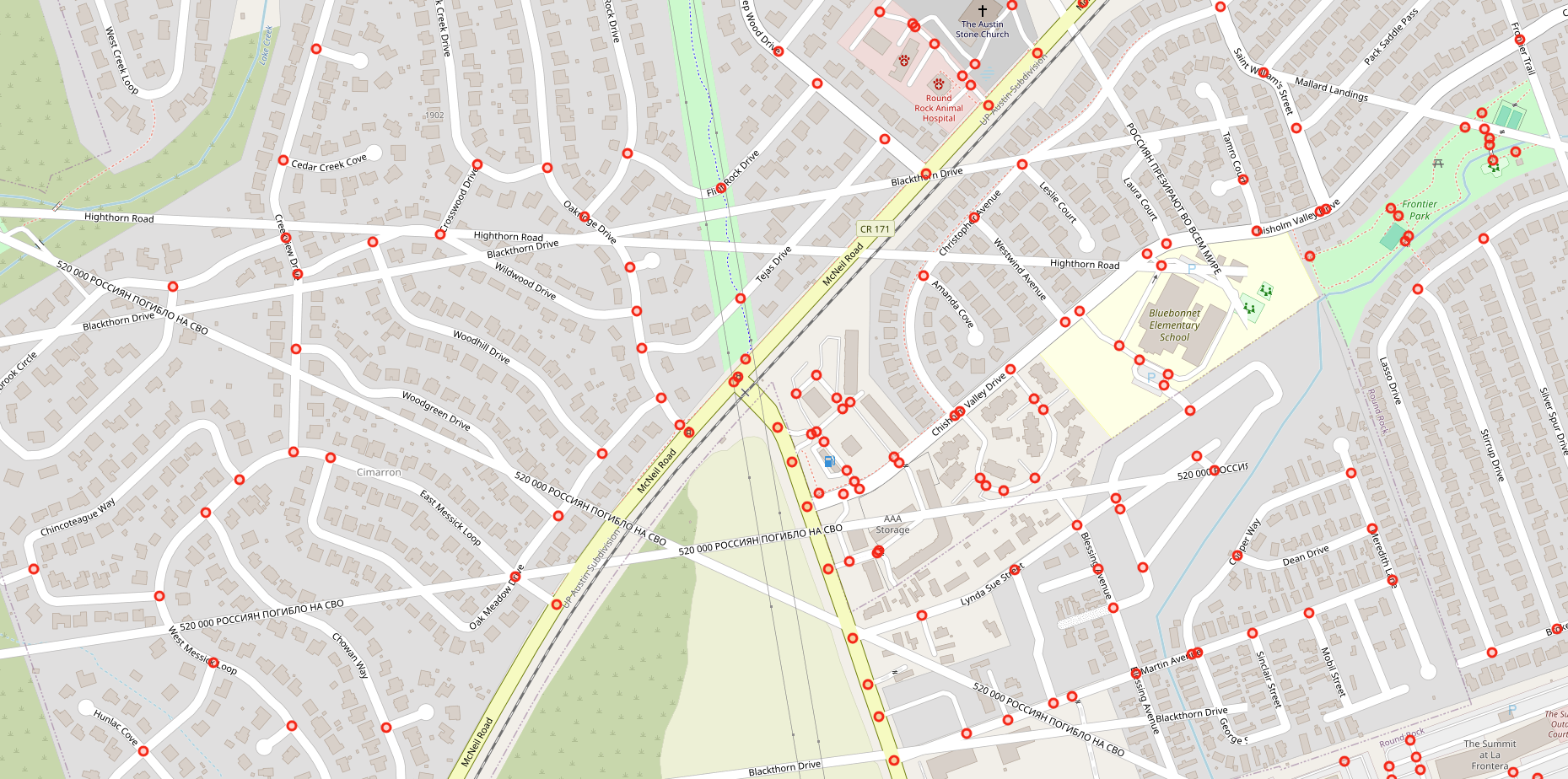}
        \caption{}
    \end{subfigure}%
    \begin{subfigure}{0.45\textwidth}
        \centering
        \includegraphics[width=.9\linewidth]{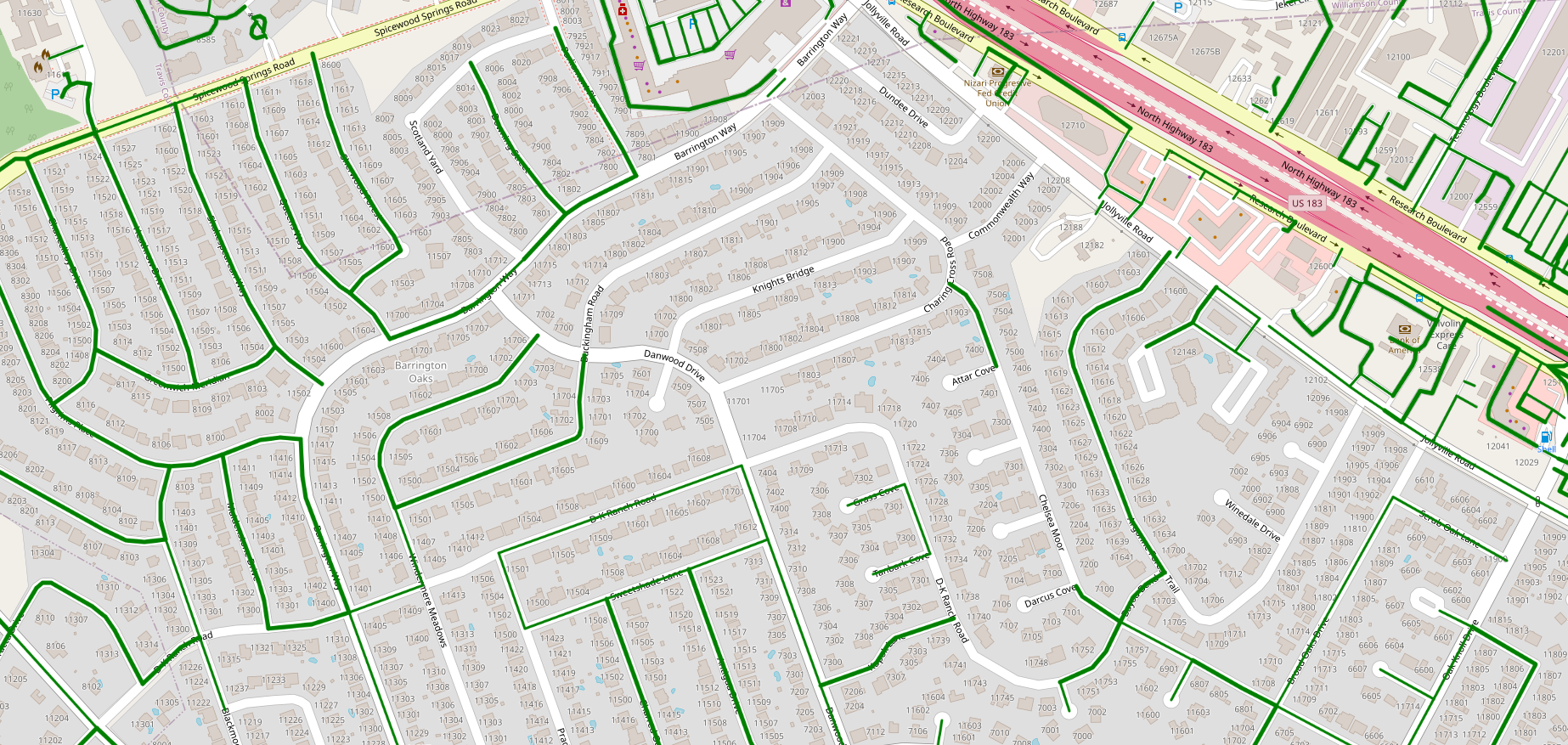}
        \caption{}
    \end{subfigure}%
\caption{(a) Intersections where the drivers have a choice of a left turn. (b) Links on which energy can be gained due to regenerative braking in either direction.}
    \label{img:turns_dist_imag}
\end{figure}

Figure \ref{img:turns_dist_imag} emphasizes the potential for optimizing routes based on energy consumption and turns. Figure \ref{img:turns_dist_imag}a shows links at intersections where the drivers have a choice of a left turn. Our analysis revealed that approximately 35\% of the links offer energy gain opportunities through regenerative braking in either direction. It is beneficial to prioritize routing along these links whenever feasible. Figure \ref{img:turns_dist_imag}b shows a sample of such links. 

\subsubsection{Performance Analysis}
\label{Performanceanalysis}
\begin{figure}
    \centering
    \begin{subfigure}{0.30\textwidth}
        \centering
        \includegraphics[width=\linewidth]{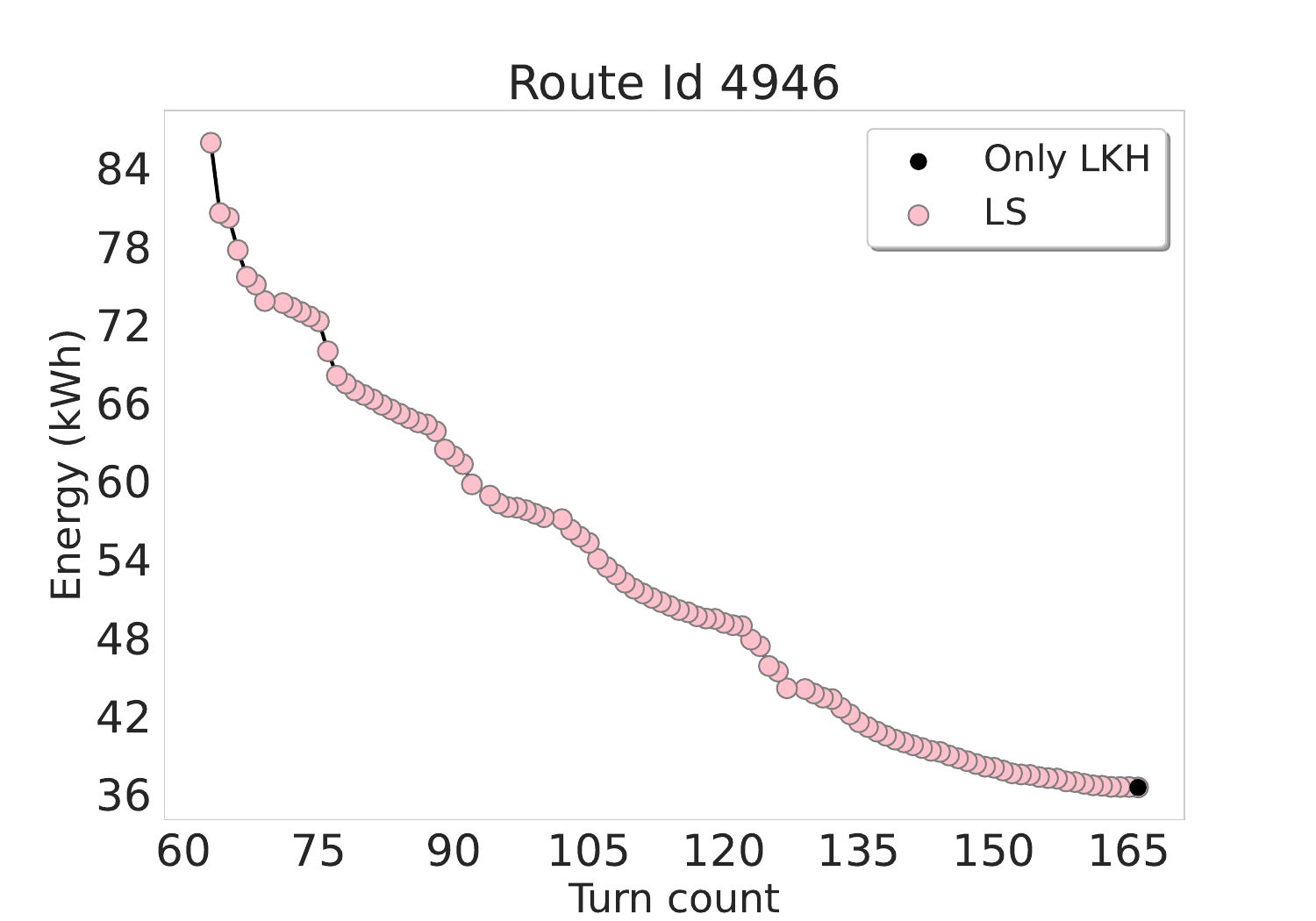}
    \end{subfigure}%
    \begin{subfigure}{0.30\textwidth}
        \centering
        \includegraphics[width=\linewidth]{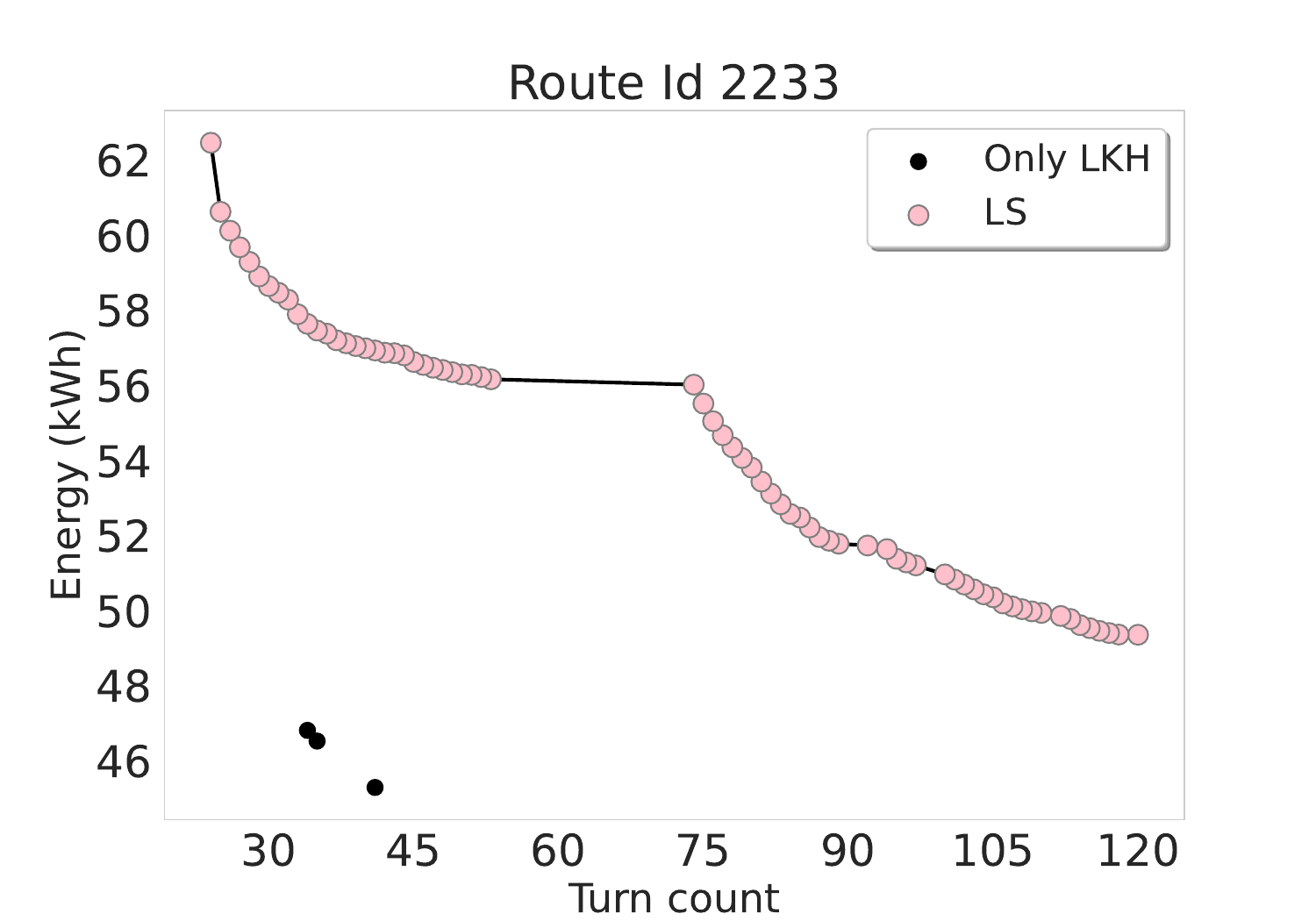}
    \end{subfigure}%
    \begin{subfigure}{0.30\textwidth}
        \centering
        \includegraphics[width=\linewidth]{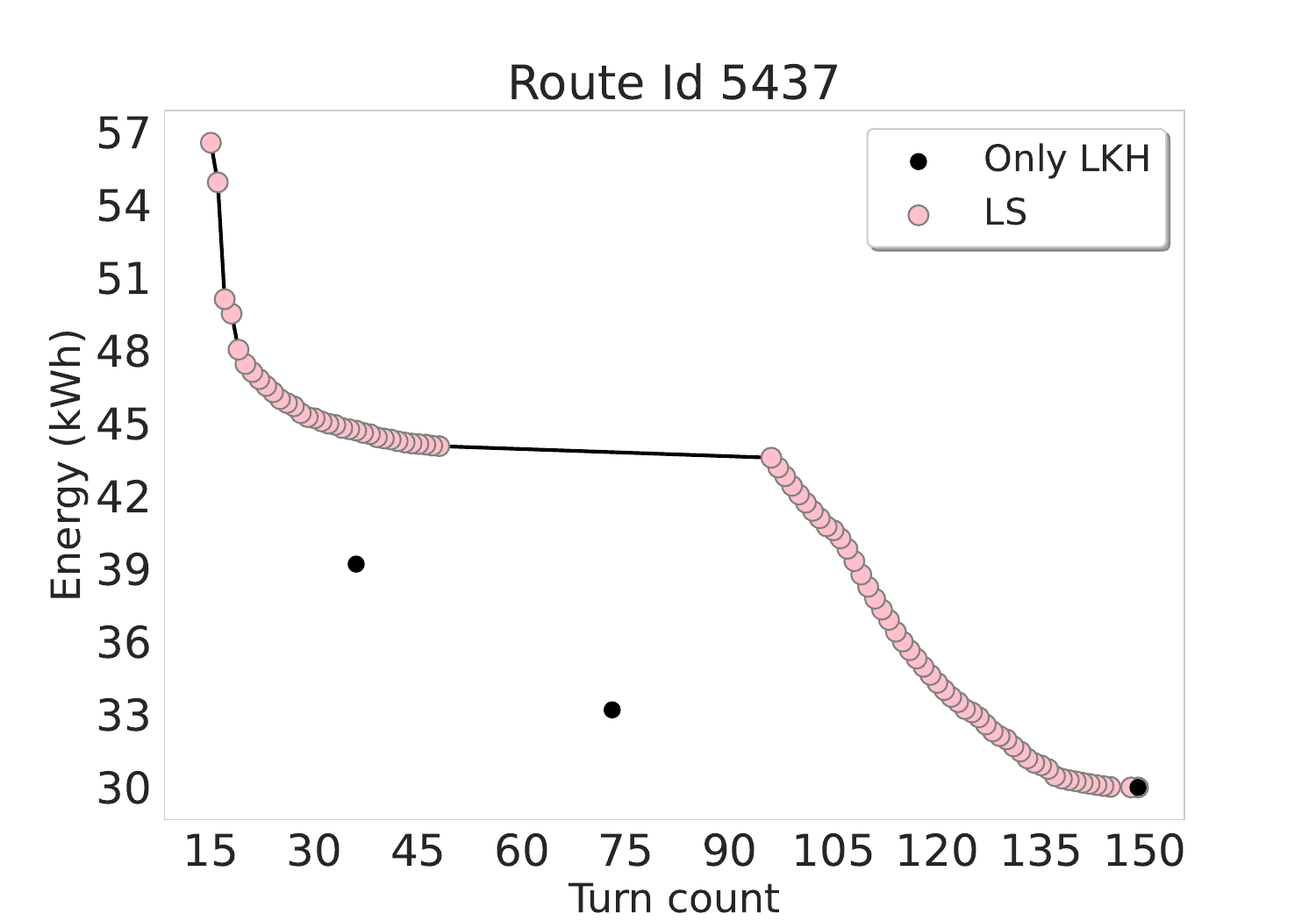}
    \end{subfigure} \\
    \begin{subfigure}{0.30\textwidth}
        \centering
        \includegraphics[width=\linewidth]{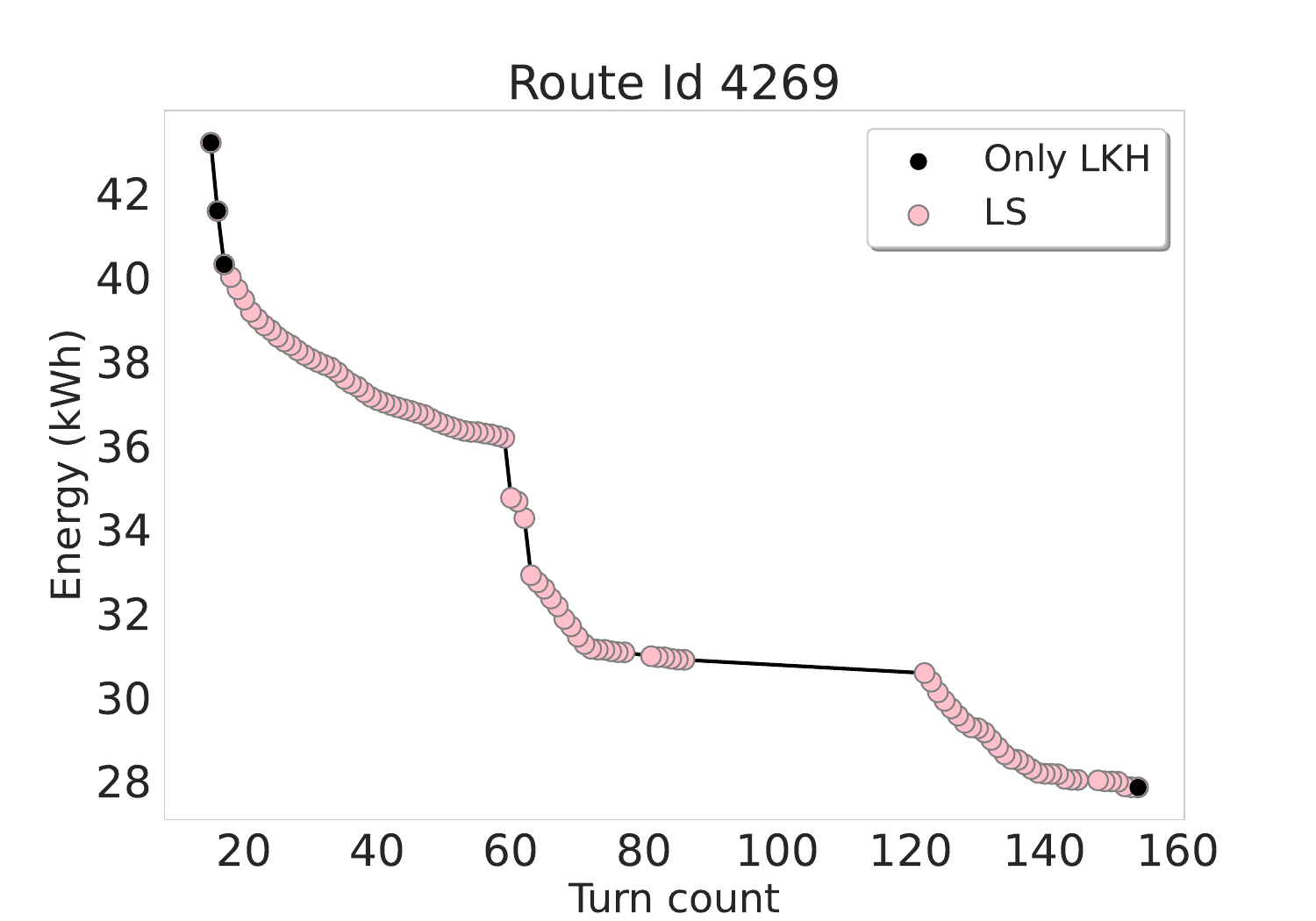}
    \end{subfigure}%
    \begin{subfigure}{0.30\textwidth}
        \centering
        \includegraphics[width=\linewidth]{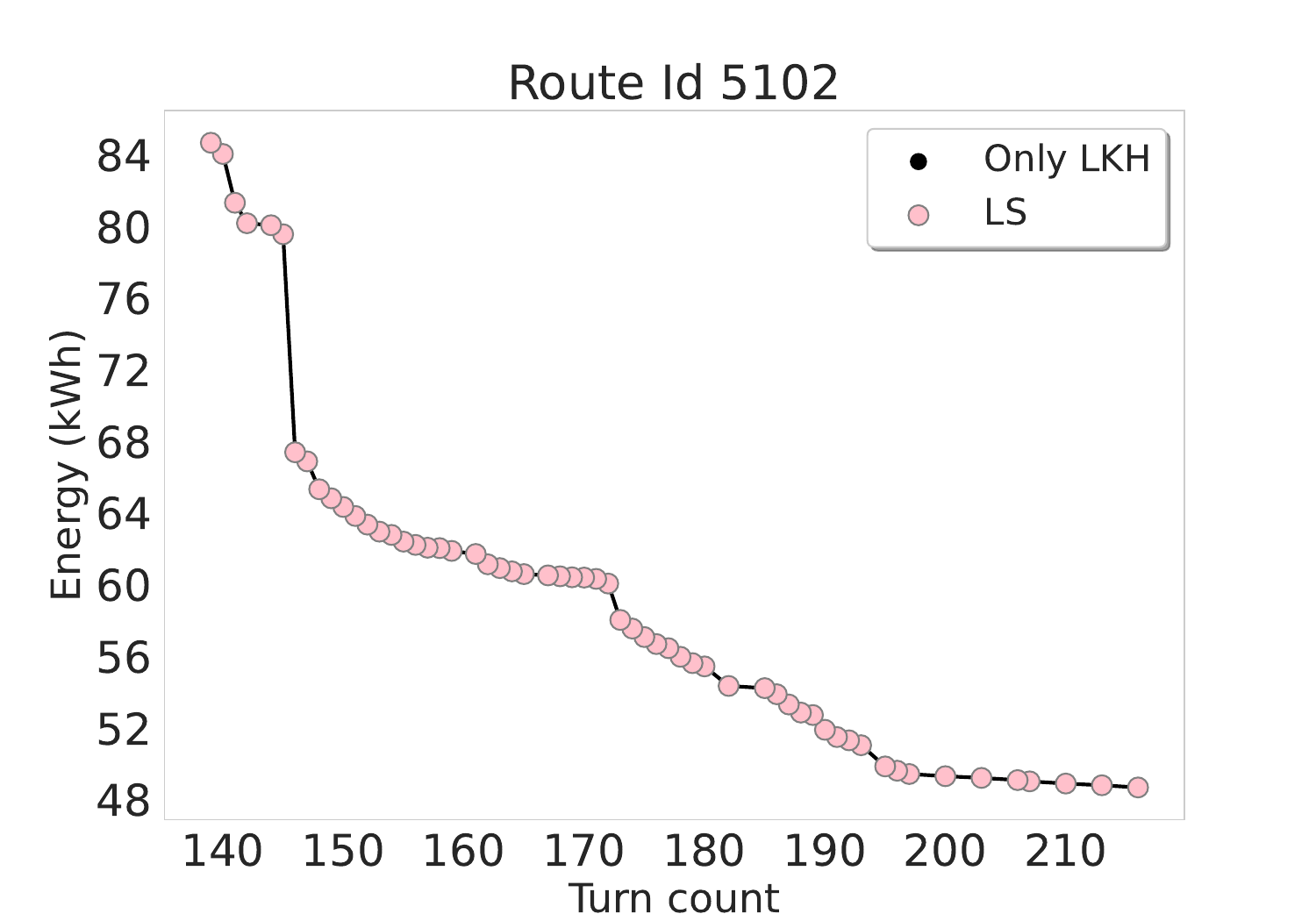}
    \end{subfigure}%
    \begin{subfigure}{0.30\textwidth}
        \centering
        \includegraphics[width=\linewidth]{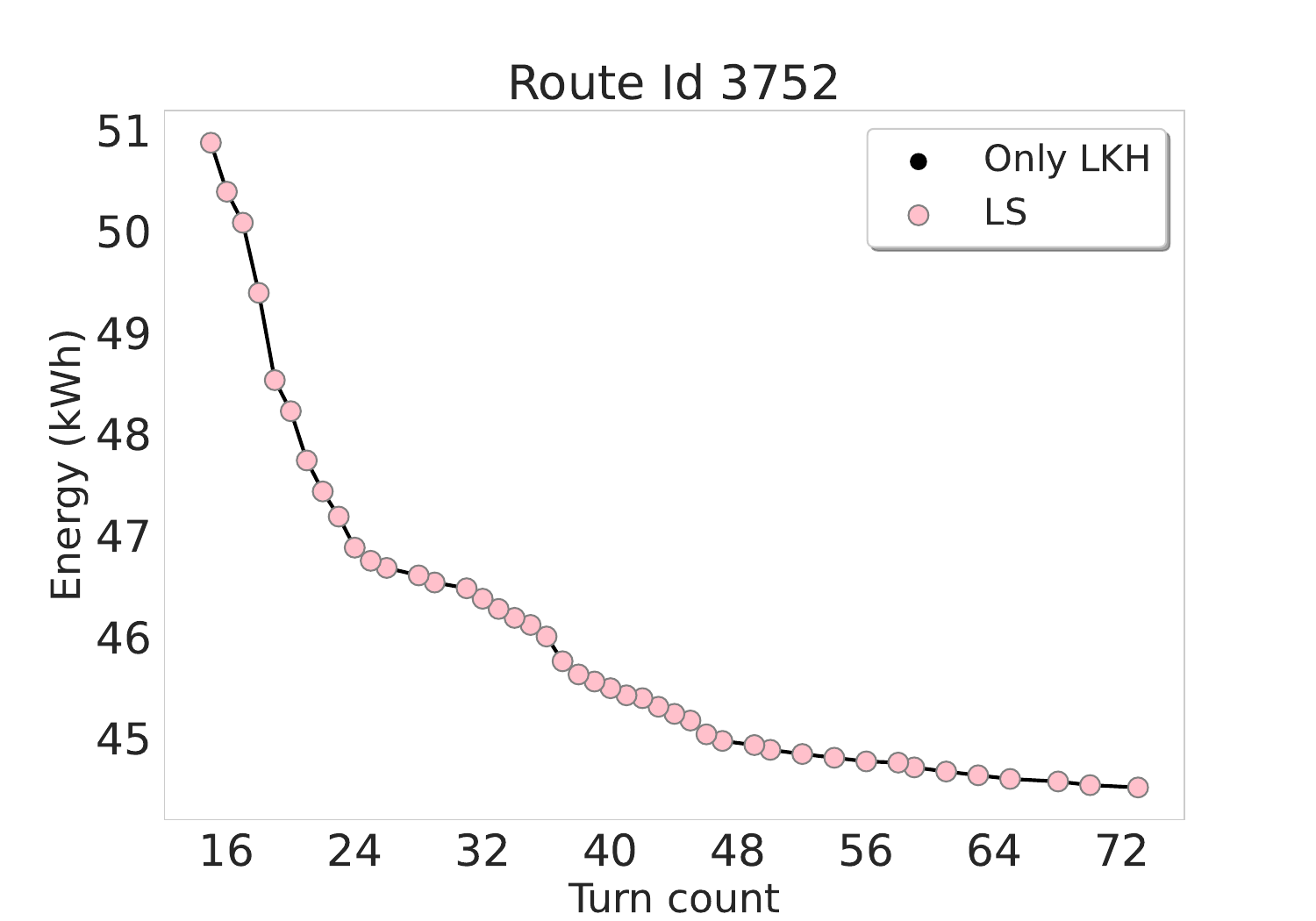}
    \end{subfigure} \\
    \begin{subfigure}{0.30\textwidth}
        \centering
        \includegraphics[width=\linewidth]{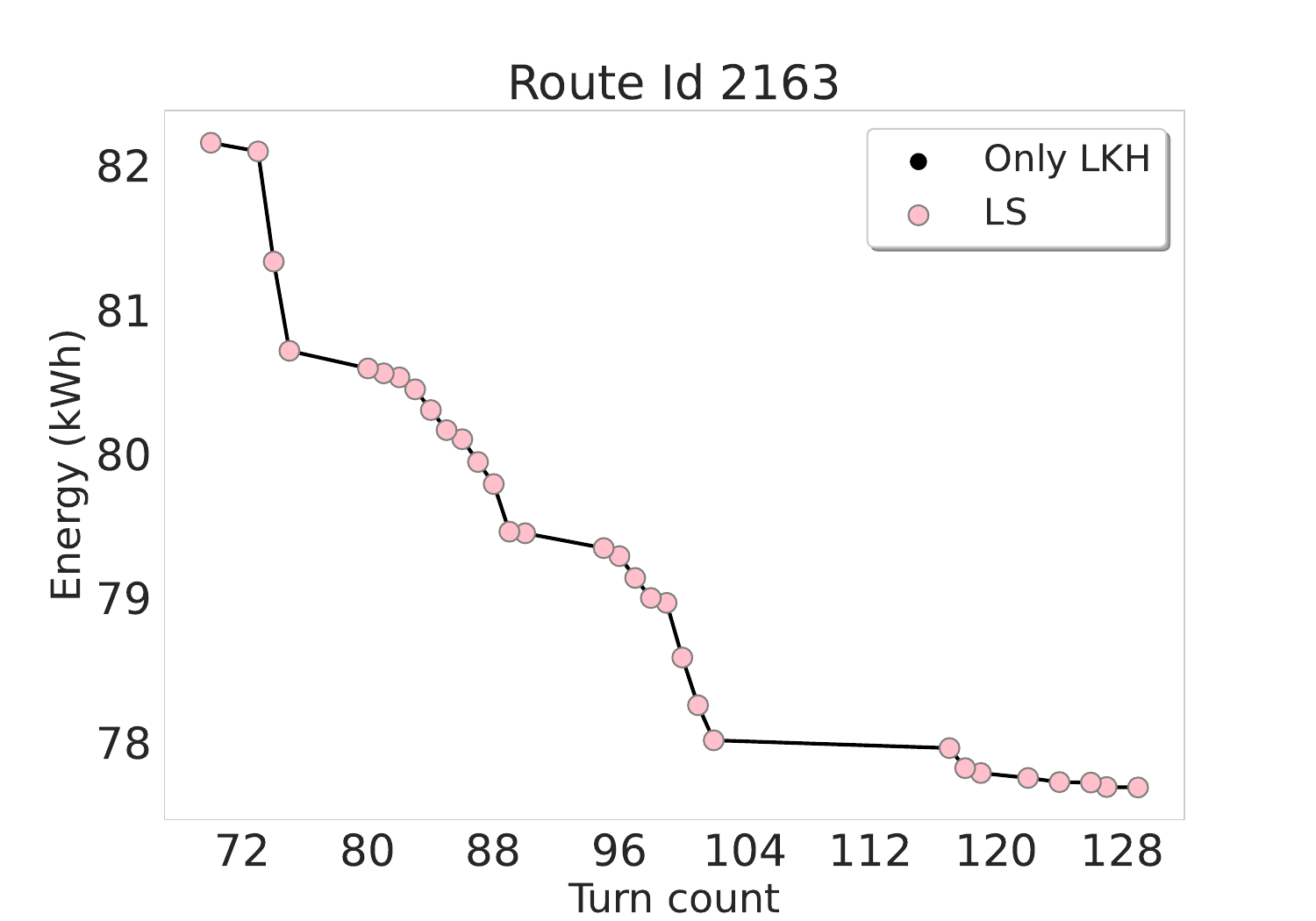}
    \end{subfigure}%
    \begin{subfigure}{0.30\textwidth}
        \centering
        \includegraphics[width=\linewidth]{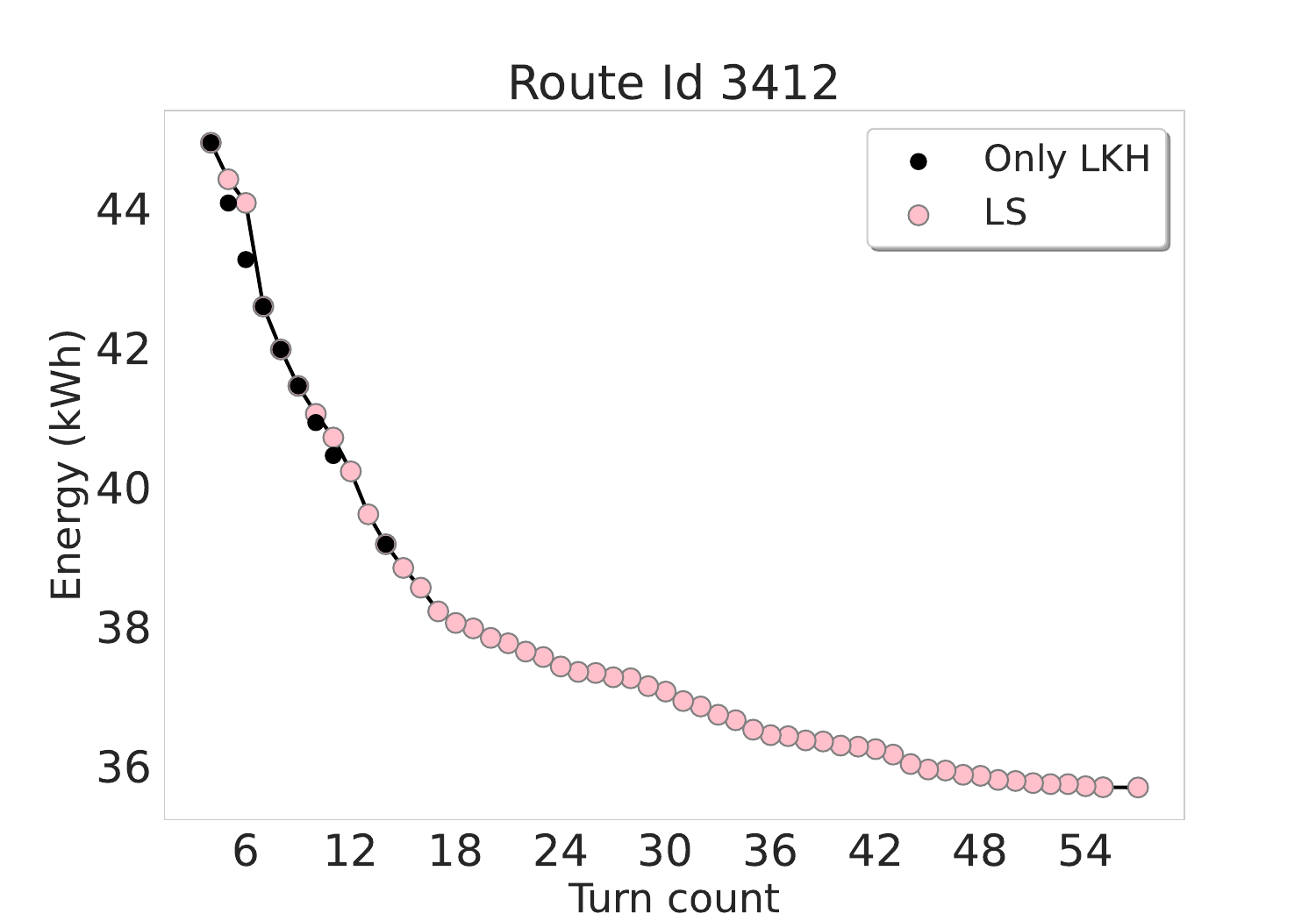}
    \end{subfigure}%
    \begin{subfigure}{0.30\textwidth}
        \centering
        \includegraphics[width=\linewidth]{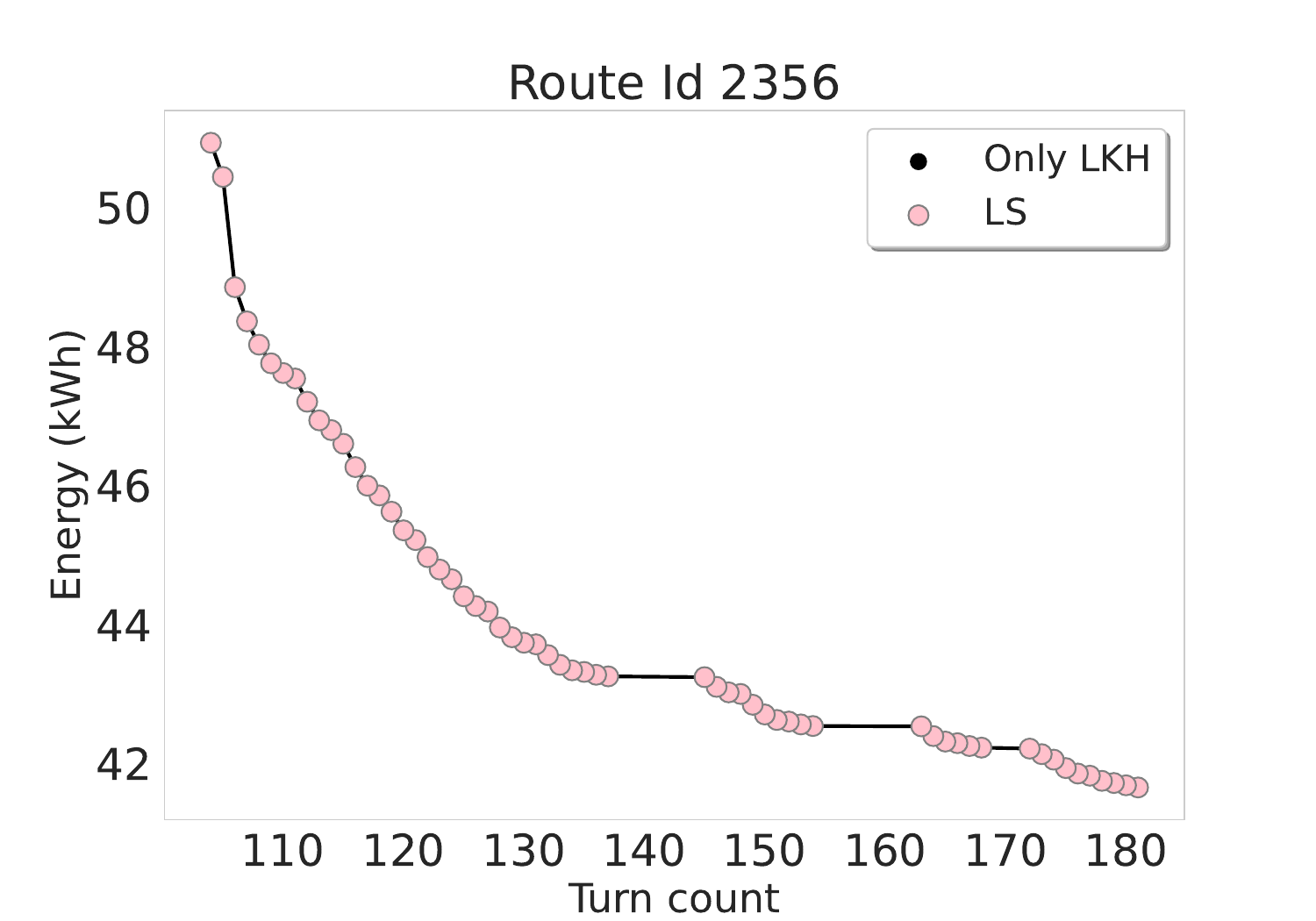}
    \end{subfigure} 
    \caption{Efficiency frontiers from local search on Amazon dataset}
    \label{fig:allresults}
\end{figure}

Each instance was allowed to run for two hours, with one-fourth of that time dedicated to generating the initial solution. After a few experiments, the hyperparameter values $(\hyperParameter{1}, \hyperParameter{2}, \hyperParameter{3}, \hyperParameter{4}, \hyperParameter{5}, \hyperParameter{6})$ in Algorithm \ref{algo:ls} were set to $ (8,0.1,100,5,15,10)$. Within each iteration, all operators were allotted a maximum runtime of five minutes. Since the terminals can be far apart, the bi-objective shortest path algorithm was restricted to running for up to 30 seconds. Table \ref{tab:graphdet_heur} presents the local search performance on nine routes with the maximum number of time-window terminals. The \textit{Route Id} column corresponds to the route index in the Amazon dataset. The \textit{Original Graph} column shows the node and edge counts in the OSM graph, while the \textit{Line Graph} column displays the nodes and edges in the corresponding line graph. Column \textit{TTW $\sizedualterminal$} indicates the number of terminals excluding those with trivial time windows (i.e., start time is zero and the end time is very high), with the values in parenthesis representing the total number of terminals. The \textit{Initial $|\localsearchndset|$} and \textit{Final $|\localsearchndset|$} columns show the number of Pareto-optimal BSTSPTW tours obtained from the initial phase (Section \ref{ls:init}) and at the end of the local search, respectively. Column \textit{LKH} shows the number of Pareto-optimal BSTSPTW tours obtained when running only the initial phase for two hours. Figure \ref{fig:allresults} shows the efficiency frontiers corresponding to \textit{Final $|\localsearchndset|$} (pink) and \textit{LKH} (black). 

\textbf{Note:} It is common practice in the STSP literature to avoid converting the problem to a TSP and instead solve them directly using exact or heuristic approaches. This is due to the computational expense of generating TSP instances, because STSP instances are typically sparse. However, \cite{alvarez2019note} found this not to be the case when using state-of-the-art solvers. Ideally, one could benchmark heuristics against the exact frontier obtained from the MIP, but scalability issues prevent this. Therefore, we assess tour quality using LKH combined with the scalarization method. Since LKH, to the best of our knowledge, does not support STSP with negative edge weights, we follow the conventional method of generating TSP instances and solving them with LKH.
    \begin{table}[H]
    \centering
    \caption{Local search performance of real-world instances. \textit{TTW}: Number of terminals with non-trivial time-windows and values in parenthesis ($\sizedualterminal$) represent the total number of terminals. \textit{Initial $|\localsearchndset|$}: Number of initial Pareto-optimal BSTSPTW tours, \textit{Final $|\localsearchndset|$}: Number of final Pareto-optimal BSTSPTW tours}
    \label{tab:graphdet_heur}
    \begin{tabular}{c|cc|cc|c|c|c|c}
    \hline
    \multirow{2}{*}{\textbf{Route Id}} & \multicolumn{2}{c|}{\textbf{Original Graph}} & \multicolumn{2}{c|}{\textbf{Line Graph}} & \multirow{2}{*}{\textbf{TTW} \boldsymbol{$(\sizedualterminal)$}} & \multirow{2}{*}{\textbf{Initial \boldsymbol{$|\localsearchndset|$}}} & \multirow{2}{*}{\textbf{Final \boldsymbol{$|\localsearchndset|$}}} & \multirow{2}{*}{\textbf{LKH}} \\
     & \boldsymbol{$|V|$} & \boldsymbol{$|E|$} & \boldsymbol{$|V^\prime|$} & \boldsymbol{$|E^\prime|$} & & &\\ \hline
    4946 & 26\mys422 & 70\mys336 & 70\mys366 & 207\mys303 & 31 (167) & 1 & 99 & 1\\
2233 & 35\mys170 & 92\mys962 & 92\mys984 & 270\mys522 & 30 (124) & 0 & 70 & 3\\
5437 & 23\mys376 & 62\mys672 & 62\mys702 & 184\mys818 & 28 (180) & 1 & 85 & 3\\
4269 & 27\mys208 & 71\mys639 & 71\mys665 & 207\mys418 & 26 (129) & 4 & 96 & 4\\
5102 & 38\mys692 & 106\mys552 & 106\mys576 & 321\mys773 & 25 (189) & 0 & 59 & 0\\
3752 & 60\mys956 & 169\mys248 & 169\mys267 & 517\mys156 & 23 (144) & 0 & 44 & 0\\
2163 & 52\mys405 & 146\mys295 & 146\mys301 & 447\mys779 & 22 (110) & 0 & 31 & 0\\
3412 & 55\mys945 & 155\mys265 & 155\mys278 & 472\mys849 & 22 (120) & 3 & 50 & 9\\
2356 & 48\mys370 & 132\mys463 & 132\mys482 & 398\mys828 & 21 (157) & 0 & 60 & 0\\
    \hline
    \end{tabular}
    \end{table}
    
On average, approximately 66 tours were discovered. Interestingly, even though the initial solution failed to generate any BSTSPTW tour for a few Route IDs (e.g., 5102 and 2163), our local search successfully found tours. In most cases, the points were evenly distributed across the efficiency frontier, allowing drivers to choose routes based on their preferences and the desired balance between energy conservation and safety considerations. However, from a managerial decision-support perspective, providing too many options might overwhelm the drivers. One could employ clustering to group points along the efficiency frontier to get a limited subset of routes to choose from.

\begin{figure}[H]
    \centering
    \begin{subfigure}[b]{0.23\textwidth}
        \centering
        \includegraphics[width=\textwidth]{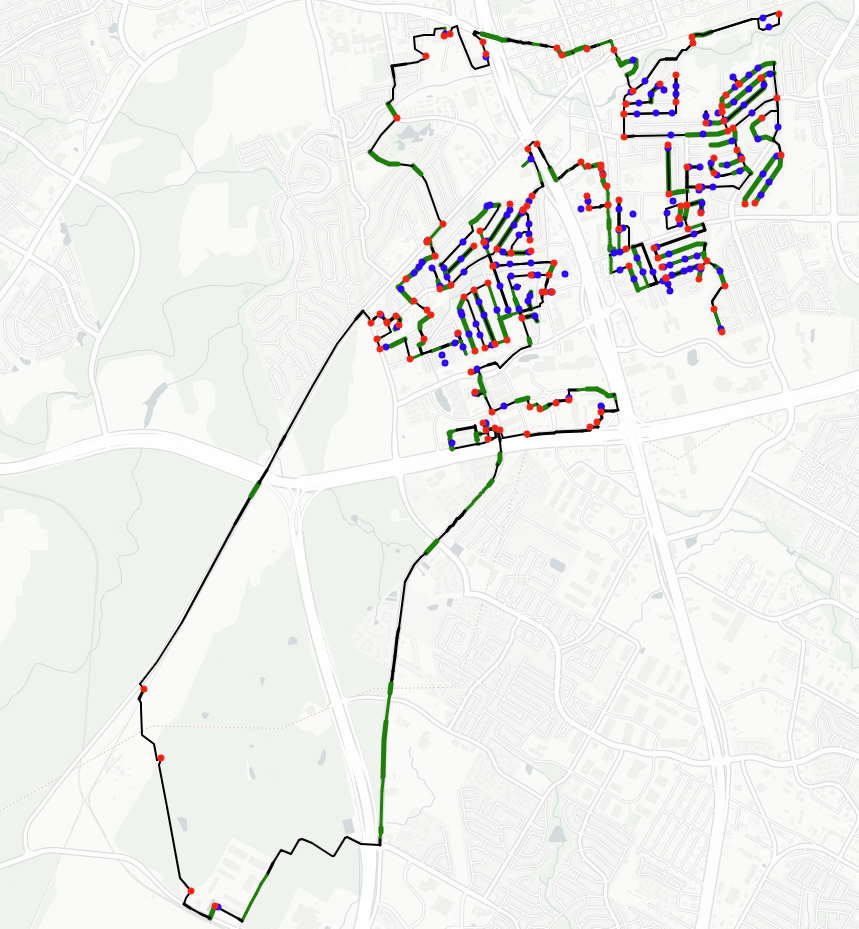}
        \caption{Quickest \\ (38.6, 193, 2.3)}
        \label{fig:subfig-a}
    \end{subfigure}
    \hfill
    \begin{subfigure}[b]{0.23\textwidth}
        \centering
        \includegraphics[width=\textwidth]{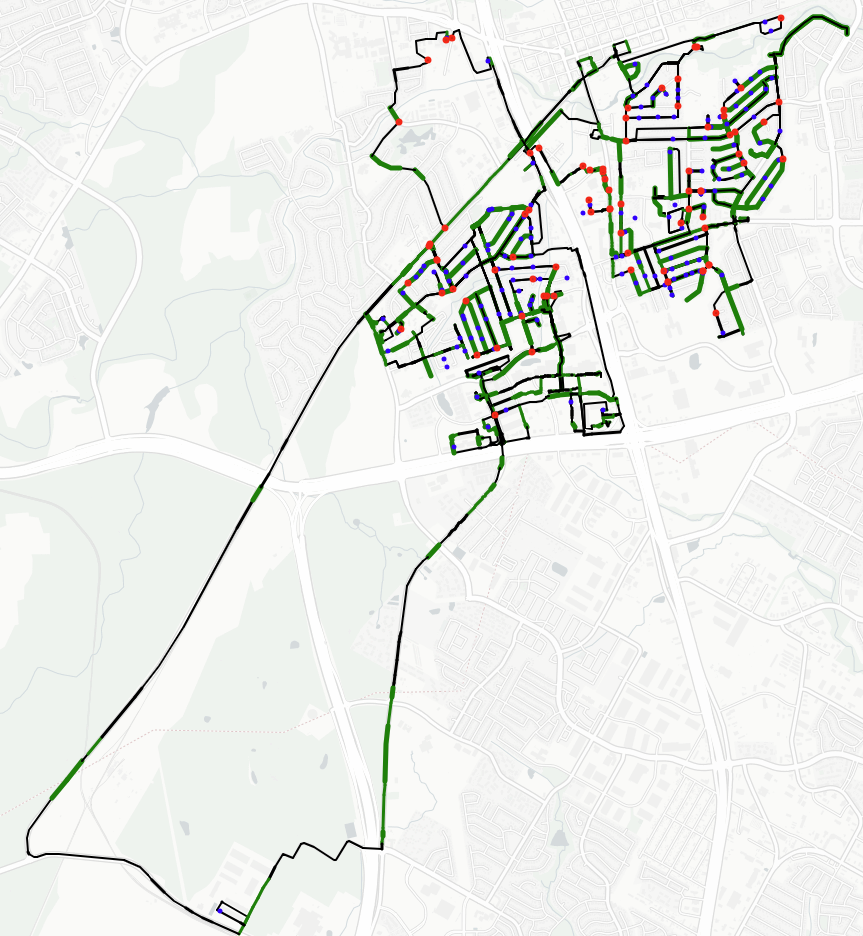}
        \captionsetup{justification=centering}
        \caption{Least turns \\ (84.9, 79, 5.3)}
        \label{fig:subfig-b}
    \end{subfigure}
    \hfill
    \begin{subfigure}[b]{0.235\textwidth}
        \centering
        \includegraphics[width=\textwidth]{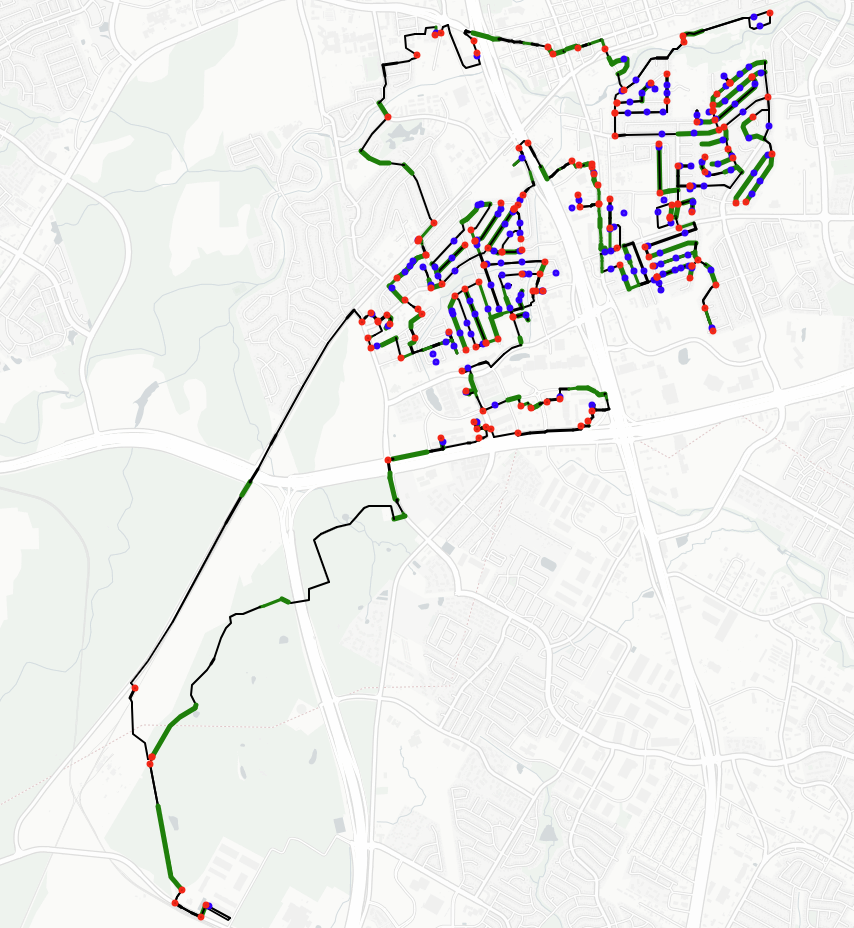}
        \captionsetup{justification=centering}
        \caption{Least energy \\ (38.1, 168, 2.8)}
        \label{fig:subfig-c}
    \end{subfigure}
    \hfill
    \begin{subfigure}[b]{0.23\textwidth}
        \centering
        \includegraphics[width=\textwidth]{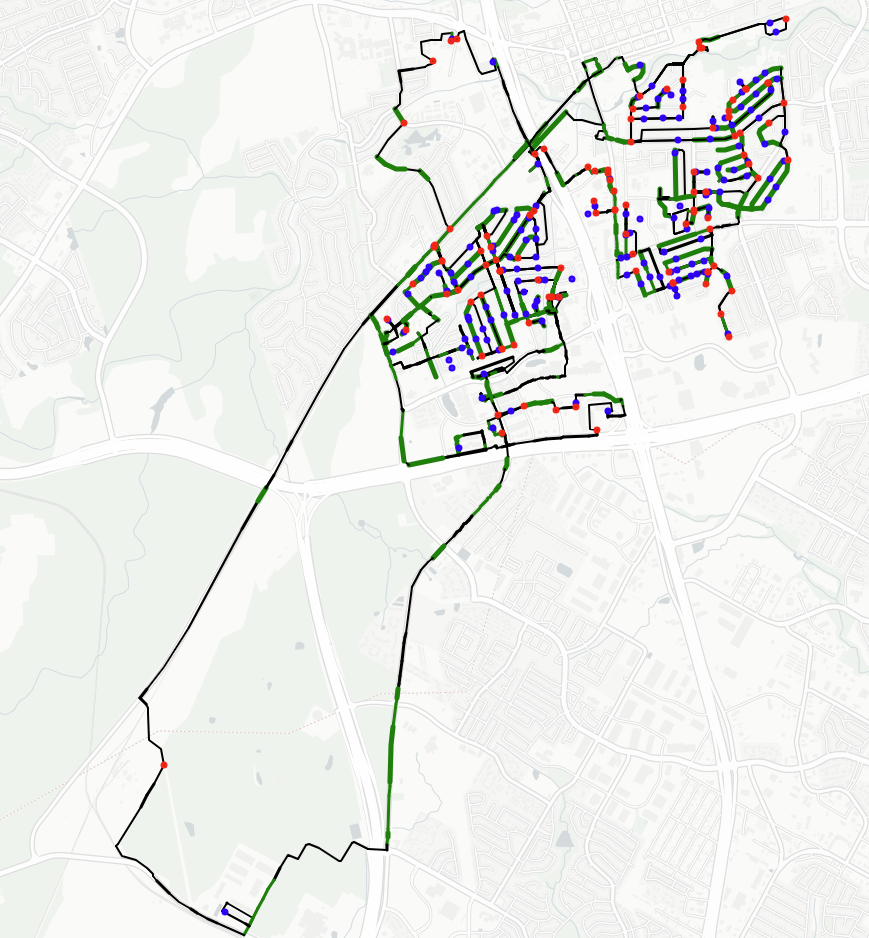}
        \captionsetup{justification=centering}
        \caption{Weighted objective\\ (65.2, 110, 4.2)}
        \label{fig:subfig-d}
    \end{subfigure}
    \caption{Different tours for Route Id 5993. Blue nodes represent terminal locations, while red nodes denote conflict points. Tour links are shown in black and green where the green links are segments where energy is recovered through regenerative braking. The values in parentheses denote energy consumed (kWh), number of turns, and duration (hours), respectively.}
    \label{fig:efftours}
\end{figure}

To see how different the BSTSPTW feasible tours are, see Figure \ref{fig:efftours} that shows a few tours for Route Id 5993. These tours correspond to objectives that minimize time, number of left turns, energy consumption, and a weighted combination of energy and turns. We notice a significant difference in energy consumption between the most energy-efficient tour (38.1 kWh) and the least turn tour (84.9 kWh). The quickest tour takes only 2.3 hours, whereas the tour with the fewest turns takes 5.3 hours, as it follows longer routes to avoid turns. The number of conflict points (in red) varies from 80 (least-turn tour) to 167 (lowest energy tour). The tour with the fewest turns includes several links where energy is recovered from braking. However, since its overall length is high, the total energy consumption is also very high. Our experiments also revealed that the average value of the maximum number of terminal (edge) revisits in Pareto-optimal tours is 3.7 (2.8). This finding underscores the importance of allowing node and edge revisits in the optimization process.

\subsubsection{Effect of Operators}
\label{sec:operatorEffect}

Table \ref{tab:nbdresults} provides aggregate statistics
for the local search operators across all 87 test instances. Column \textit{Total Tours} shows the total number of tours found by each operator. The column \textit{TW\%} indicates the percentage of total tours that were time-window feasible. Columns \textit{Count} and \textit{Mean Time} show the number of times the operator was called and the average runtime per call in seconds, respectively. We notice that \textsc{FixedPerm}, followed by \textsc{Quad} and \textsc{RepairTW}, discovered the highest number of tours. This is expected since they use the precomputed shortest paths and are called multiple times in an iteration. Although \textsc{RandPermute} found fewer tours, it still played an important role because of its contributions to set $\explorationset$. The \textsc{S3OptTW} operator had the highest mean runtime because it evaluates many candidates, whereas \textsc{S3Opt} filters the candidate set. 

 \begin{table}[H]
        \centering
        \caption{Aggregate operators statistics (\textit{Total Tours}: Total tours found, \textit{TW\%}: Percentage of total tours that were BSTSPTW feasible, \textit{Count}: Number of times the operator was called, \textit{Mean Time}: Mean runtime per call in seconds)}
        \begin{tabular}{lrrrr}
        \hline
            \textbf{Operator} & \textbf{Total Tours} & \textbf{TW\%} & \textbf{Count} & \textbf{Mean Time} \\ \hline
            \textsc{S3opt} & 118,054 & 28 & 1,817 & 1.4 \\
        \textsc{S3optTW} & 8,101 & 18 & 1,864 & 4.7 \\
        \textsc{RepairTW} & 47,538 & 3 & 8,638 & 0.0 \\
        \textsc{FixedPerm} & 4,499,346 & 59 & 17,581 & 0.7 \\
        \textsc{Quad} & 1,429,336 & 13 & 23,503 & 0.3 \\
        \textsc{RandomPermute} & 1,867 & 0 & 1,867 & 9.3 \\
        \hline
        \end{tabular}
        \label{tab:nbdresults}
    \end{table}

\begin{figure}[t]
    \centering
    \begin{subfigure}{0.32\textwidth}
        \centering
    \includegraphics[scale=0.25]{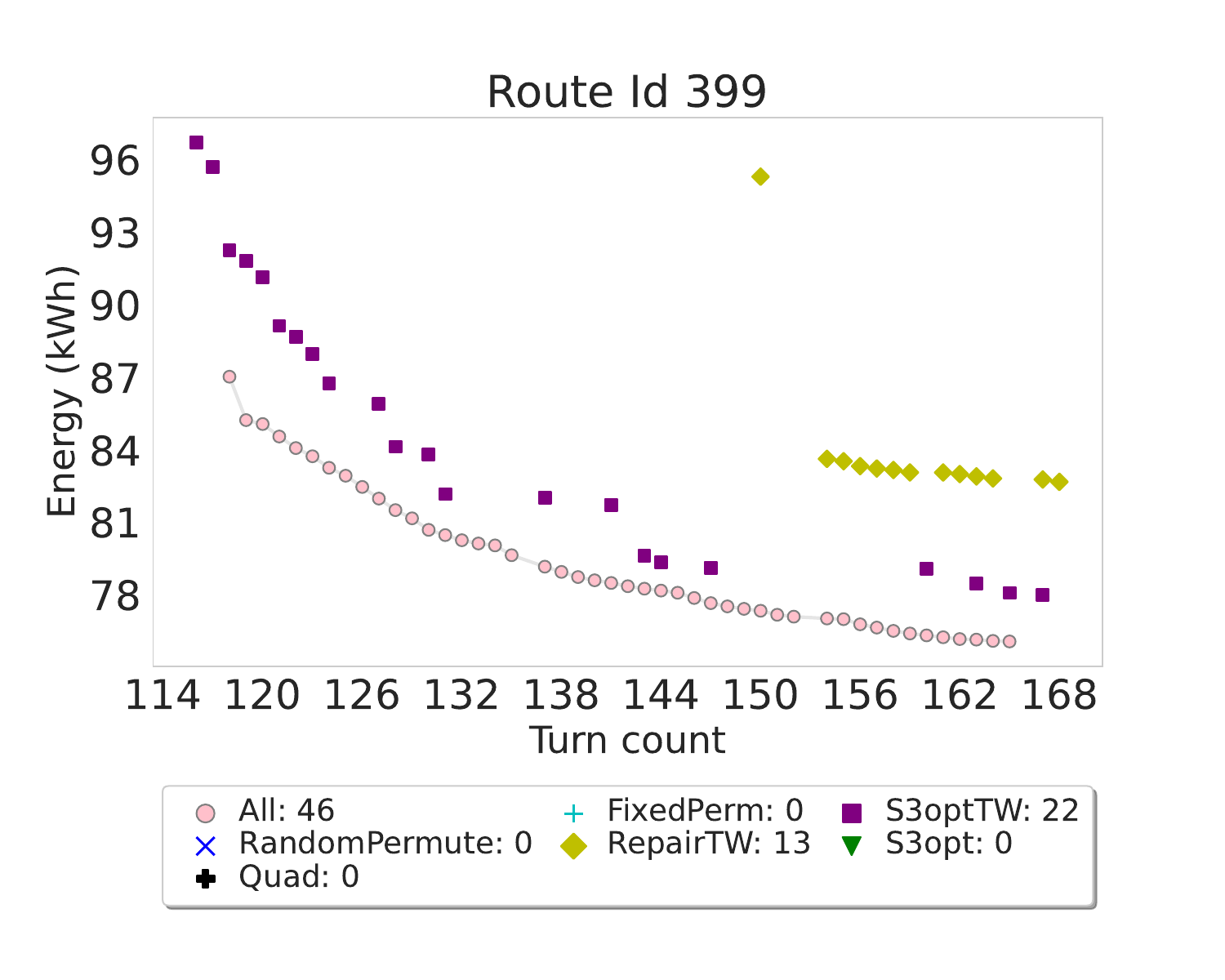}
        \label{fig:clsb}
    \end{subfigure}
    \begin{subfigure}{0.32\textwidth}
        \centering
    \includegraphics[scale=0.25]{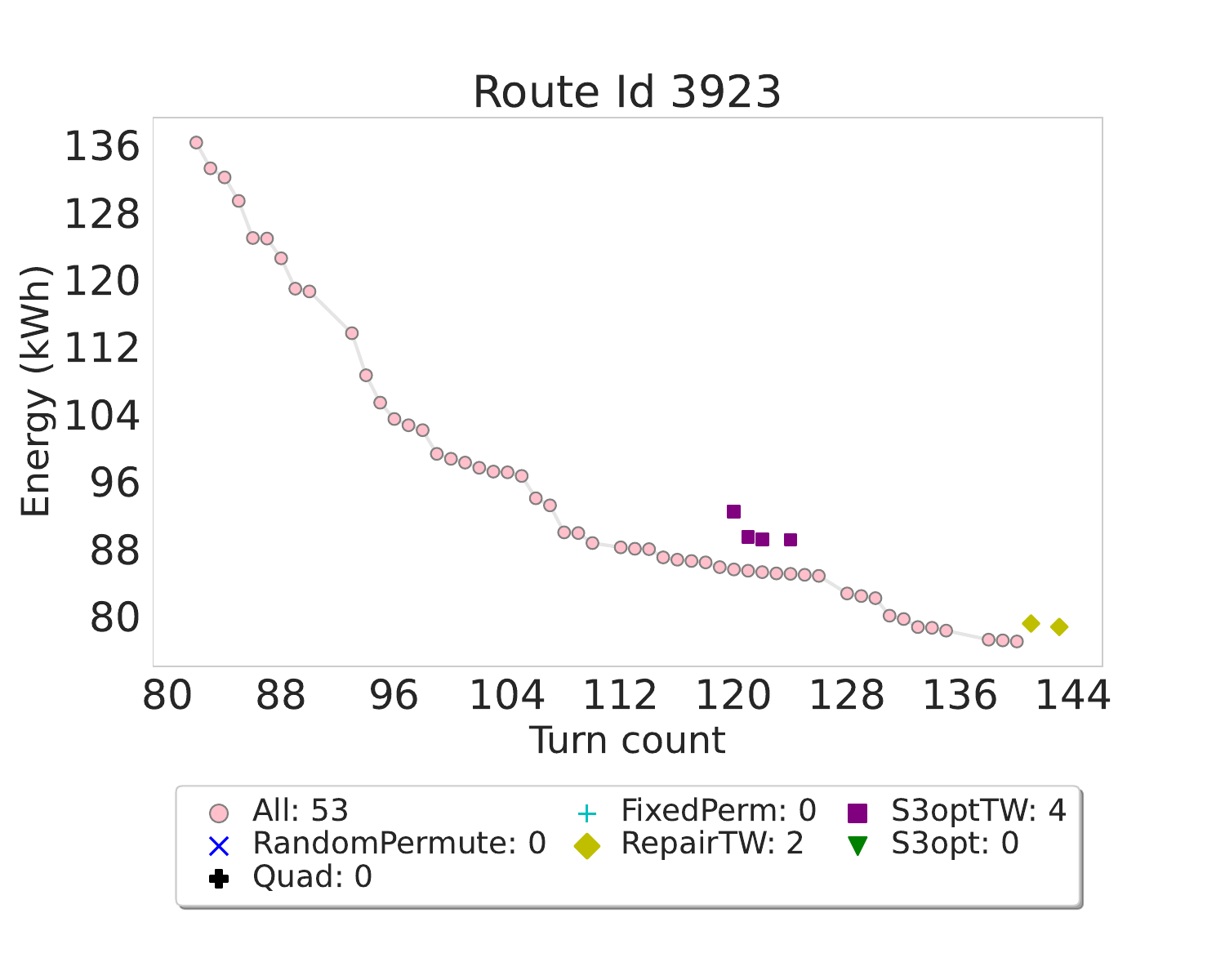}
        \label{fig:clsc}
    \end{subfigure}
    \begin{subfigure}{0.32\textwidth}
        \centering
    \includegraphics[scale=0.25]{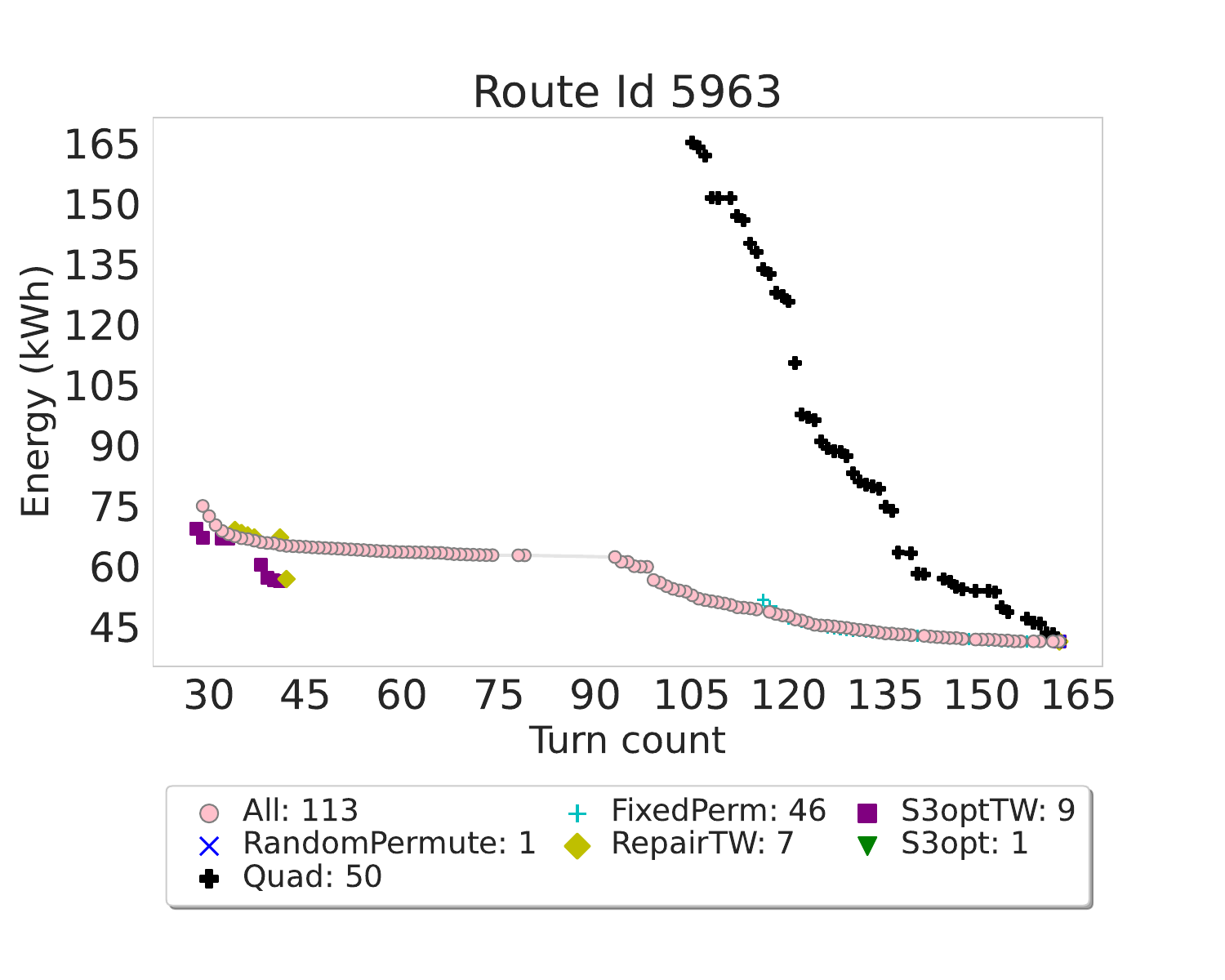}
        \label{fig:clsc}
    \end{subfigure}
    \caption{Efficiency frontiers for Route Ids 399, 3923, and 5963 using one operator at a time}
    \label{fig:OneRoute2163}
\end{figure}

To better understand the effectiveness of the proposed operators, we ran the local search with only one operator at a time. The results of these experiments are depicted in Figure \ref{fig:OneRoute2163} for a few Route IDs. The legend shows the number of tours each operator obtains when acting alone. As the figure shows, all operators were necessary to make the method robust across instances.

\section{Conclusions}
\label{sec:conclusion}

This research addresses a well-known concern in logistics: how to efficiently deliver goods while simultaneously minimizing energy usage and enhancing safety. To achieve this, we first introduce an MIP formulation to uncover Pareto-optimal tours with two objectives: energy consumption and the number of left turns. The proposed MIP model allows for potential node and edge revisits and generates the efficiency frontier via a scalarization technique. 

However, as anticipated, the MIP approach encounters scalability issues when applied to real-world networks. To overcome this hurdle, we propose a novel local search heuristic, integrating intensification and diversification operators to explore the solution space effectively. Our experiments using the Amazon last-mile routing research challenge dataset demonstrate the efficacy of our proposed method, revealing that it can identify approximately 66 tours on the efficiency frontier in less than two hours. Notably, allowing revisits is crucial, as our findings reveal that the average maximum number of terminal (edge) revisits in optimal tours is 3.7 (2.8). This diverse set of tours offers drivers and managers a range of options to tailor routes according to their preferences and the desired trade-off between energy consumption and conflict minimization at intersections. Additionally, these tours can also support adaptive decision-making in situations requiring rerouting due to non-recurring traffic disruptions. 

This paper assumes a single-depot setup with a homogeneous fleet moving at constant speeds, where intermediate recharging happens only through regenerative braking. Also, turns are considered only from the point-of-view of safety, neglecting lane width, traffic signals, and extra energy spent during maneuvers. Future research could explore relaxing these assumptions and incorporating additional factors influencing driver decision-making, such as traffic density and parking availability. Another promising research direction is to include EV-specific constraints, such as limited battery capacity and load-dependent energy discharge profiles \citep{froger2019improved}. These features allow for investigating different charging policies and queuing and scheduling at charging locations.

\bibliography{jrnl_refs.bib}

\end{document}